\documentclass[twoside, 10pt]{article}

\usepackage{amsmath}
\usepackage{paralist}
\usepackage{amssymb,graphicx}
\usepackage{hyperref}
\def\v{\varepsilon}

\def\x{\xi}
\def\t{\theta}
\def\T{\Theta}
\def\k{\kappa}

\def\mb{\mathbf}
\def\a{\alpha}
\def\b{\beta}

\def\d{\delta}
\def\l{\lambda}
\def\f{\frac}
\def\p{\phi}
\def\r{\rho}

\def\s{\sigma}
\def\z{\zeta}
\def\o{\omega}
\def\di{\displaystyle}
\def\i{\infty}
\def\wt{\widetilde}

\newtheorem{theorem}{Theorem}[section]

\newtheorem{lemma}[theorem]{Lemma}
\newtheorem{proposition}{Proposition}[section]

\newtheorem{remark}{Remark}

\begin{document}

\title{\bf  The  Limit of the Boltzmann
Equation to the  Euler Equations for Riemann Problems} \vskip 0.5cm
\author{
Feimin Huang,\thanks{Academy of Mathematics and Systems
   Science, Chinese Academy of Sciences, Beijing, 100190, China ({\tt fhuang@amt.ac.cn}),
   supported in part by NSFC Grant No.
10825102 for distinguished youth scholar, and National Basic
Research Program of China (973 Program) under Grant No.
2011CB808002.} \quad Yi Wang,\thanks{ Academy of Mathematics and
Systems
   Science, Chinese Academy of Sciences, Beijing, 100190, China ({\tt
wangyi@amss.ac.cn}), supported by the NSFC Grant No. 10801128.}
\quad Yong Wang,\thanks{Academy of Mathematics and Systems
   Science, Chinese Academy of Sciences, Beijing, 100190, China ({\tt yongwang@amss.ac.cn}).}
\quad Tong Yang\thanks{Department of Mathematics, City University of
HongKong, Hong Kong ({\tt matyang@cityu.edu.hk}), supported by the
General Research Fund of Hong Kong, CityU 104310 and the Croucher Foundation.} }
\date{}
\maketitle
\begin{abstract}
The convergence of the Boltzmann equaiton to the compressible Euler
equations when the Knudsen number tends to zero has been a long
standing open problem in the kinetic theory. In the setting of
Riemann solution that contains the generic superposition of shock,
rarefaction wave and contact discontinuity to the Euler equations,
we succeed in justifying this limit by introducing hyperbolic waves
with different solution backgrounds to capture the extra masses
carried by the hyperbolic approximation of the rarefaction wave and
the diffusion approximation of contact discontinuity.\\[3mm]
 \noindent{\it 2010 Mathematics Subject Classification:}  76P05, 35Q20,
74J40, 82B40, 82C40.

 \noindent{\it Key
words and phrases:} Hydrodynamic limit, Boltzmann equation, Euler
equations, Riemann solution
\end{abstract}
 \tableofcontents
\section{Introduction}
\renewcommand{\theequation}{\arabic{section}.\arabic{equation}}
\setcounter{equation}{0}

As the fundamental equation in statistical mechanics, the Boltzmann
equation takes the form of
\begin{equation*}
f_t + \xi\cdot \nabla_Xf= \f{1}{\v}Q(f,f),
\end{equation*}
where  $f(t,X,\xi)$ is the
 density distribution function of  particles at time $t$ with location
 $X$  and velocity $\xi$.
In this equation, the physical parameter $\v>0$ called   Knudsen
number
 is proportional to the mean free path of the interacting particles.

It was known since its derivation that the Boltzmann equation is
closely related to the systems of fluid dynamics, in particular, the
systems of Euler and Navier-Stokes equations. In fact, the first
derivation of the fluid dynamical components and systems
 from the kinetic equations can be
traced back to the dates of  Maxwell  and Boltzmann.  Their early
derivations rest on some arguments as  how the various terms in a
kinetic equation balance each other.  These balance arguments seem
arbitrary to some extent. For this,
 Hilbert proposed  a
systematic  expansion in 1912, and  Enskog and Chapman independently
proposed another expansion in 1916 and 1917 respectively.

Either the Hilbert expansion or Chapman-Enskog expansion yields the
compressible Euler equations in the leading order with respect to
the Knudsen number $\v$, and the compressible Navier-Stokes
equations, Burnett equations in the subsequent orders.   To justify
these formal approximations in rigorous mathematics, that is,
hydrodynamic limits,
 has been proved to be extremely challenging and most remains open,
 in part because the basic
well-posedness and regularity questions are still mostly unsolved
 for these fluid
equations.

The justification of the fluid limits of the Boltzmann equation is also related to  the
Hilbert's sixth problem, "Mathematical treatment of the axioms of
physics", in which it says that "The investigations on the
foundations of geometry suggest the problem: To treat in the same
manner, by means of axioms, those physical sciences in which
mathematics plays an important part; in the first rank are the
theory of probabilities and mechanics. ...... Thus Boltzmann's work
on the principles of mechanics suggests the problem of developing
mathematically the limiting processes, there merely indicated, which
lead from the atomistic view to the laws of motion of continua.
...... Further, the mathematician has the duty to test exactly in
each instance whether the new axioms are compatible with the
previous ones."

The goal of this paper is to justify the limiting process of the
Boltzmann equation to the system of the compressible Euler equations
in the setting of Riemann solutions. The Riemann problem was first
formulated and studied by Riemann in 1860s when he studied the one
space dimensional isentropic gas dynamics with initial data being
two constant states. The solution to this problem turns out to be
fundamental in the theory of hyperbolic conservation laws because it
not only captures  the local and global behavior of solutions, but
also fully represents the effect of the nonlinearity in the
structure of the solutions. It is now well known that for the system
of Euler equations, there are three basic wave patterns, that is,
shock wave, rarefaction wave and contact discontinuity. These three
types of waves have essential differences, that is, shock is
compressive, rarefaction is expansive, and contact discontinuity has
some diffusive structure. Therefore, how to study the hydrodynamic
limit of Boltzmann equation for the full Riemann solution that
consists of the superposition of these three typical waves is still
very challenging in mathematics.

In this paper, by introducing two types of  hyperbolic waves that carry
the extra masses in the hyperbolic approximation of the rarefaction
wave profile and the diffusive approximation of the contact
discontinuity, we succeed in proving rigorously that there exists a
family of solutions to the Boltzmann equation that converges to a
Maxwellian determined by a  Riemann solution consisting of three
basic wave patterns when the Knudsen number tends to zero.
Furthermore, a convergence rate is obtained in term of the Knudsen
number.

By coping with the essential properties of individual wave pattern,
the hydrodynamic limit for a single wave was justified in the previous
works separately. More precisely, by using the compressibility of
shock wave profile, Yu \cite{Yu} showed that when the solution of
the Euler equations (\ref{(1.5)}) contains only non-interacting
shocks, there exists a sequence of solutions to
 the Boltzmann equation that converge to a local Maxwellian defined by the
 solution of the Euler equations (\ref{(1.5)}) uniformly
away from the shock in any fixed time interval. In this work,
a generalized Hilbert
expansion was introduced, and the analytic technique of matching the
inner and outer expansions developed by Goodman-Xin
\cite{Goodman-Xin} for conservation laws was
 used.  On the other hand, by using the time decay
properties of the rarefaction wave, similar problem was studied by
Xin-Zeng \cite{Xin-Zeng}. Moreover, by using the diffusive structure
in the contact discontinuity as for the Navier-Stokes equations, the
hydrodynamic limit to the contact discontinuity was proved by
Huang-Wang-Yang in  \cite{Huang-Wang-Yang}.

However, up to now, how to deal with the general Riemann solution
that consists of all three basic waves is still a challenging open
problem. This is mainly due to the difficulty in handling the wave
interactions and also unifying the different approaches in the
analysis used for each single wave pattern. In order to overcome
these difficulties in justifying the limit, our main idea in this
paper is to introduce two families of hyperbolic waves, called
hyperbolic wave I and II, that capture the propagation of the extra
mass created by the approximate hyperbolic rarefaction wave profile
in the viscous setting and the diffusion approximation of contact
discontinuity.

We now briefly explain why the two families of hyperbolic waves we
introduced are essential for the proof. As in the previous works on
the rarefaction wave in the setting of either Navier-Stokes
equations or the Boltzmann equation, the approximate rarefaction
wave is constructed as a hyperbolic wave profile. Therefore, we need
to precisely capture the error in the second order of the
approximation for the Boltzmann equation in term of the Knudsen
number, that is, in the Navier-Stokes level. And this reduces to
study the propagation of the extra mass induced by the viscosity and
heat conductivity. For this, we introduce the hyperbolic wave I as a
solution to the linearized system around the approximate rarefaction
wave profile with source terms given by the viscosity and heat
conductivity induced by the rarefaction wave profile to recover the
viscous terms. We can show that the hyperbolic wave I decays like
the first-order derivative of the rarefaction wave profile so that
the decay properties given in Lemma \ref{hw1}  are good enough to
carry out the analysis.

The main difficulty comes from the approximation of the contact
discontinuity. First of all, such an approximation, that is,
2-viscous contact wave, behaves like a diffusion wave profile as for
the Navier-Stokes equations. Due to the lack of sufficient decay in
$\v$ and the non-conservative error terms when taking the
anti-derivative of the perturbation, we need to remove the leading
error terms and non-conservative terms  in such approximation before
taking the anti-derivative. The hyperbolic wave II is constructed to
remove these error terms due to the viscous contact wave
approximation. Note that the construction of the hyperbolic wave II
can not be done simply around the 2-viscous contact wave as the
hyperbolic wave I for the rarefaction wave profile. Otherwise, the
wave interaction terms thus induced will lead to insufficiently
decay in $\v$ due to 2-viscous contact wave and it seems that these
terms are essential in wave interactions. Instead, it is constructed
around the superposition of the approximate 1-rarefaction wave, the
hyperbolic wave I, the 2-viscous contact wave and  the 3-shock
profile as a whole. Thus some wave interaction terms can be absorbed
in the hyperbolic wave II and the other wave interaction terms can
be handled by some subtle and careful calculations. Due to the
non-conservative terms and insufficient decay rates of $\varepsilon$
of error terms induced by the 2-contact wave, we can not use the
anti-derivative technique to analyze the hyperbolic wave II. Since
the derivative of 3-shock profile is negative and tends to infinity
as the Knudsen number $\varepsilon\to 0+$, we have to impose the
condition at the time $t=T$ (see \eqref{A2.8} below) for the
linearized hyperbolic system of hyperbolic wave II so that the
monotonicity of 3-shock wave is fully utilized. This idea is
different from the previous stability analysis on the shock profile
which  is based on the anti-derivative technique.

With the help of these two hyperbolic waves and the corresponding
new estimates, we can justify the limiting process from the
Boltzmann equation to compressible Euler equations for the generic
Riemann problems by elaborate analysis after a hyperbolic scaling.

We now formulate the problem. Consider the Boltzmann
equation with slab symmetry
\begin{equation}
f_t + \xi_1 f_x= \f{1}{\v}Q(f,f), \label{(B)}
\end{equation}
where $\xi=(\xi_1,\xi_2,\xi_3)\in {\mb R}^3$ and $x\in {\mb R}^1$.
Here,  the collision operator takes the form of
$$
 Q(f,g)(\xi) \equiv \f{1}{2}
\int_{{\mb R}^3}\!\!\int_{{\mb S}_+^2} \Big(f(\xi')g(\xi_*')+
f(\xi_*')g(\xi')-f(\xi) g(\xi_*) - f(\xi_*)g(\xi) \Big)
B(|\xi-\xi_*|, \hat\t)
 \; d \xi_* d\Omega,
$$
where $\xi',\xi_*'$ are the velocities after an elastic collision of
two particles with velocities  $\xi,\xi_*$ before the collision.
Here, $\hat\t$ is the angle between the relative velocity
$\xi-\xi_*$ and the unit vector $\Omega$ in ${\mb S}^2_+=\{\Omega\in
{\mb S}^2:\ (\xi-\xi_*)\cdot \Omega\geq 0\}$. The conservations of
momentum and energy yield the following relations between the
velocities before and after collision:
$$
 \xi'= \xi -[(\xi-\xi_*)\cdot \Omega] \; \Omega, \qquad
 \xi_*'= \xi_* + [(\xi-\xi_*)\cdot \Omega] \; \Omega.
$$

We will concentrate on
 the hard sphere model where  the cross-section takes the form of
$$
B(|\xi-\xi_*|,\hat\t)=|(\xi-\xi_*, \Omega)|=|\x-\x_*|\cos\hat\t.
$$
On the other hand, it is noted that  the analysis can be
applied to at least hard potential if we can assume Lemma
\ref{Lemma-shock} on the shock wave profile holds true.

As we mentioned earlier, formally,  when the Knudsen number $\v$
tends to zero, the limit  of the Boltzmann equation \eqref{(B)} is
 the system of compressible Euler equations
that consists of conservations of mass, momentum and energy:
\begin{equation}
\left\{
\begin{array}{l}
\di\rho_t+(\rho u_1)_x=0,\\
\di(\rho u_1)_t+(\rho u_1^2+p)_x=0,\\
\di(\rho u_i)_t+(\rho u_1u_i)_x=0,~i=2,3,\\
\di[\rho(e+\f{|u|^2}{2})]_t+[\rho u_1(E+\f{|u|^2}{2})+pu_1]_x=0,
\end{array}
\right. \label{(1.5)}
\end{equation}
where
\begin{equation}
\left\{
\begin{array}{l}
\di\rho(t,x)=\int_{\mb{R}^3}\varphi_0(\xi)f(t,x,\xi)d\xi,\\
\di\rho u_i(t,x)=\int_{\mb{R}^3}\varphi_i(\xi)f(t,x,\xi)d\xi,~i=1,2,3,\\
\di\rho(e+\f{|u|^2}{2})(t,x)=\int_{\mb{R}^3}\varphi_4(\xi)f(t,x,\xi)d\xi.
\end{array}
\right. \label{(1.6+)}
\end{equation}
Here, $\rho$ is the density, $u=(u_1,u_2,u_3)$ is the macroscopic
velocity, $e$ is the internal energy, and $p=R\rho\t$
with $R$ being the gas constant is  the pressure. Note that the
temperature $\t$ is related to the internal energy by
$e=\f{3}{2}R\t$, and $\varphi_i(\xi)(i=0,1,2,3,4)$ are the collision
invariants given by
\begin{equation}
 \varphi_0(\xi) = 1,~~~
 \varphi_i(\xi) = \xi_i~  (i=1,2,3),~~~
 \varphi_4(\xi) = \f{1}{2} |\xi|^2,
\label{collision-invar}
\end{equation}
that satisfy
$$
\int_{{\mb R}^3} \varphi_i(\xi) Q(g_1,g_2) d \xi =0,\quad {\textrm
{for} } \ \  i=0,1,2,3,4.
$$

Instead of using either Hilbert expansion or Chapman-Enskog
expansion, we will apply  the macro-micro decomposition introduced
in  \cite{Liu-Yang-Yu}.
 For a solution $f(t,x,\xi)$ of (\ref{(B)}), set
$$
f(t,x,\xi)=\mb{M}(t,x,\xi)+\mb{G}(t,x,\xi),
$$
where the local Maxwellian
$\mb{M}(t,x,\xi)=\mb{M}_{[\rho,u,\t]}(\xi)$ represents the
macroscopic  component of the solution defined by the five conserved
quantities, i.e., the mass density $\rho(t,x)$, the momentum $\rho
u(t,x)$, and the total energy $\rho(e+\f{1}{2}|u|^2)(t,x)$ given in
(\ref{(1.6+)}), through
\begin{equation}
\mb{M}=\mb{M}_{[\rho,u,\t]} (t,x,\xi) = \f{\rho(t,x)}{\sqrt{ (2 \pi
R \t(t,x))^3}} e^{-\f{|\xi-u(t,x)|^2}{2R\t(t,x)}}. \label{(1.7)}
\end{equation}
And $\mb{G}(t,x,\xi)$  represents the microscopic component.

From now on, the inner product of $g_1$ and $g_2$ in
$L^2_{\xi}({\mb R}^3)$ with respect to a given Maxwellian
$\tilde{\mb{M}}$ is denoted by:
\begin{equation}
 \langle g_1,g_2\rangle_{\tilde{\mb{M}}}\equiv \int_{{\mb R}^3}
 \f{1}{\tilde{\mb{M}}}g_1(\xi)g_2(\xi)d \xi.\label{product}
\end{equation}
 If $\tilde{\mb{M}}$ is the local
Maxwellian $\mb{M}$ defined in (\ref{(1.7)}),
  the macroscopic space is spanned by the
following five pairwise orthogonal base,
\begin{equation}
\left\{
\begin{array}{l}
 \chi_0(\xi) \equiv {\di\f1{\sqrt{\rho}}\mb{M}}, \\[2mm]
 \chi_i(\xi) \equiv {\di\f{\xi_i-u_i}{\sqrt{R\t\rho}}\mb{M}} \ \ {\textrm {for} }\ \  i=1,2,3, \\[2mm]
 \chi_4(\xi) \equiv
 {\di\f{1}{\sqrt{6\rho}}(\f{|\xi-u|^2}{R\t}-3)\mb{M}},\\
 \langle\chi_i,\chi_j\rangle=\delta_{ij}, ~i,j=0,1,2,3,4.
 \end{array}
\right.\label{orthogonal-base}
\end{equation}
For brevity, if $\tilde{\mb{M}}$ is the local Maxwellian $\mb{M}$,
we will simply use $\langle\cdot,\cdot\rangle$ to denote
  $\langle\cdot,\cdot\rangle_{\mb{M}}$. By using the above base,
the macroscopic projection $\mb{P}_0$ and microscopic projection
$\mb{P}_1$ can be defined as
\begin{equation*}
 \mb{P}_0g = {\di\sum_{j=0}^4\langle g,\chi_j\rangle\chi_j},\qquad
 \mb{P}_1g= g-\mb{P}_0g.
\end{equation*}
Note that a function $g(\xi)$ is called microscopic  if
$$
\int g(\xi)\varphi_i(\xi)d\xi=0,~i=0,1,2,3,4,
$$
where again $\varphi_i(\xi)$ represents the collision invariants.

Notice that the solution $f(t,x,\xi)$ to the Boltzmann equation
(\ref{(B)}) satisfies
$$
\mb{P}_0f=\mb{M},~~~\mb{P}_1f=\mb{G},
$$
and the Boltzmann equation (\ref{(B)}) becomes
\begin{equation}
(\mb{M}+\mb{G})_t+\xi_1(\mb{M}+\mb{G})_x
=\f{1}{\v}[2Q(\mb{M},\mb{G})+Q(\mb{G},\mb{G})]. \label{(1.9)}
\end{equation}
By integrating the product of the equation (\ref{(1.9)}) and the
collision invariants $\varphi_i(\xi)(i=0,1,2,3,4)$  with respect to
$\xi$ over ${\mb R}^3$, one has the following  system for the fluid
variables $(\rho, u, \theta)$:
\begin{equation}
\left\{
\begin{array}{lll}
\di \rho_{t}+(\rho u_1)_x=0, \\
\di (\rho u_1)_t+(\rho u_1^2
+p)_x=-\int\xi_1^2\mb{G}_xd\xi,  \\
\di (\rho u_i)_t+(\rho u_1u_i)_x=-\int\xi_1\xi_i\mb{G}_xd\xi,~ i=2,3,\\
\di [\rho(e+\f{|u|^2}{2})]_t+[\rho
u_1(e+\f{|u|^2}{2})+pu_1]_x=-\int\f12\xi_1|\xi|^2\mb{G}_xd\xi.
\end{array}
\right. \label{(1.10)}
\end{equation}

 Note that the above fluid-type system is not
self-contained and one more  equation for the microscopic
 component ${\mb{G}}$
is needed and it
 can be obtained by applying the projection operator
$\mb{P}_1$ to  (\ref{(1.9)}):
\begin{equation}
\mb{G}_t+\mb{P}_1(\xi_1\mb{M}_x)+\mb{P}_1(\xi_1\mb{G}_x)
=\f{1}{\v}\left[\mb{L}_\mb{M}\mb{G}+Q(\mb{G}, \mb{G})\right].
\label{(1.11)}
\end{equation}
Here $\mb{L}_\mb{M}$ is the linearized collision operator of
$Q(f,f)$
 with respect to the local Maxwellian $\mb{M}$ given by
$$
\mb{L}_\mb{M} g=2Q(\mb{M}, g)=Q(\mb{M}, g)+ Q(g,\mb{M}).
$$
Note that  the null space $\mathfrak{N}$ of $\mb{L}_\mb{M}$ is
spanned by the macroscopic variables:
$$
\chi_j(\xi), ~j=0,1,2,3,4.
$$
Furthermore, there exists a positive constant $\wt\sigma>0$ such
that for any function $g(\xi)\in \mathfrak{N}^\bot$, cf.
\cite{Grad},
$$
\langle g,\mb{L}_\mb{M}g\rangle \le -\wt\sigma\langle
\nu(|\xi|)g,g\rangle ,
$$
where $\nu(|\xi|)=O(1)(1+|\xi|)$ is the collision frequency for the
hard sphere model.

Consequently, the linearized collision operator $\mb{L}_\mb{M}$ is a
dissipative operator on $L^2({\mb R}^3)$, and its inverse
$\mb{L}_\mb{M}^{-1}$  is a bounded operator on $\mathfrak{N}^\bot$.
It follows from (\ref{(1.11)}) that
\begin{equation}
\mb{G}=\v \mb{L}_\mb{M}^{-1}[\mb{P}_1(\xi_1\mb{M}_x)] +\Pi,
\label{(1.12)}
\end{equation}
with
\begin{equation}
\Pi=\mb{L}_\mb{M}^{-1}[\v(\mb{G}_t+\mb
{P}_1(\xi_1\mb{G}_x))-Q(\mb{G}, \mb{G})].\label{(1.13)}
\end{equation}
Plugging  (\ref{(1.12)}) into (\ref{(1.10)}) gives
\begin{equation}
\left\{
\begin{array}{l}
\di \rho_{t}+(\rho u_1)_x=0,\\
\di (\rho u_1)_t+(\rho u_1^2 +p)_x=\f{4\v}{3}(\mu(\t)
u_{1x})_x-\int\xi_1^2\Pi_xd\xi,  \\
\di (\rho u_i)_t+(\rho u_1u_i)_x=\v(\mu(\t)
u_{ix})_x-\int\xi_1\xi_i\Pi_xd\xi,~ i=2,3,\\
\di [\rho(\t+\f{|u|^2}{2})]_t+[\rho
u_1(\t+\f{|u|^2}{2})+pu_1]_x=\v(\k(\t)\t_x)_x+\f{4\v}{3}(\mu(\t)u_1u_{1x})_x\\
\di\qquad\qquad +\v\sum_{i=2}^3(\mu(\t)u_iu_{ix})_x
-\int\f12\xi_1|\xi|^2\Pi_xd\xi,
\end{array}
\right. \label{(1.14)}
\end{equation}
where the viscosity coefficient $\mu(\t)>0$ and the heat
conductivity coefficient $\k(\t)>0$ are smooth functions of the
temperature $\t$. Here,  we normalize the gas constant $R$ to be
$\f{2}{3}$ so that $e=\t$ and $p=\f23\rho\t$.

Since the problem considered in this paper is  one dimensional in
the space variable $x\in {\bf R}$, in the macroscopic level,
 it is
more convenient to rewrite the equation (\ref{(B)}) and the system
(\ref{(1.5)}) in the {\it Lagrangian} coordinates. For this,
 set the coordinate transformation:
\begin{equation}
  (t,x)\Rightarrow \Big(t,\int_{(0,0)}^{(t,x)} \rho(\tau,y)dy-(\r
 u_1)(\tau,y)d\tau\Big),\label{Lag}
\end{equation}
where $\int_A^B fdy+gd\tau$ represents a line integraton from point $A$ to point $B$ on
${\mb R}^+\times {\mb R}$. Here, the value of the integration is unique because of the conservation of mass.

We will still denote the {\it Lagrangian} coordinates by $(t,x)$ for
the simplicity of notations. Then (\ref{(B)}) and (\ref{(1.5)}) in
the Lagrangian coordinates become, respectively,
\begin{equation}
f_t-\f{u_1}{v}f_x+\f{\xi_1}{v}f_x=\f{1}{\v}Q(f,f),\label{Lag-B}
\end{equation}
and
\begin{equation}
\left\{
\begin{array}{llll}
\di v_{t}-u_{1x}=0,\\
\di u_{1t}+p_x=0,\\
\di u_{it}=0, ~i=2,3,\\
\di (\t+\f{|u|^{2}}{2}\bigr)_t+ (pu_1)_x=0.\\
\end{array}
\right.\label{(1.16)}
\end{equation}
Moreover, (\ref{(1.10)})-(\ref{(1.14)}) take the form of
\begin{equation}
\left\{
\begin{array}{llll}
\di v_t-u_{1x}=0,\\
\di u_{1t}+p_x=-\int\xi_1^2\mb{G}_xd\xi,\\
\di u_{it}=-\int\xi_1\xi_i\mb{G}_xd\xi,
~i=2,3,\\
\di\bigl(\t+\f{|u|^{2}}{2}\bigr)_{t}+ (pu_1)_x=-\int\f12\xi_1|\xi|^2\mb{G}_xd\xi,\\
\end{array}
\right.\label{(1.17)}
\end{equation}
\begin{equation}
\mb{G}_t-\f{u_1}{v}\mb{G}_x+\f{1}{v}\mb{P}_1(\xi_1\mb{M}_x)+\f{1}{v}\mb{P}_1(\xi_1\mb{G}_x)=\f{1}{\v}(\mb{L}_\mb{M}\mb{G}+Q(\mb{G},\mb{G})),\label{(1.18)}
\end{equation}
with
\begin{equation}
\mb{G}=\v \mb{L}^{-1}_\mb{M}(\f{1}{v} \mb{P}_1(\xi_1
\mb{M}_x))+\Pi_1,\label{(1.19)}
\end{equation}
\begin{equation}
\Pi_1=\mb{L}_\mb{M}^{-1}[\v(\mb{G}_t-\f{u_1}v\mb{G}_x+\f{1}{v}\mb{P}_1(\xi_1\mb{G}_x))-Q(\mb{G},\mb{G})],\label{(1.21)}
\end{equation}
and
\begin{equation}
\left\{
\begin{array}{llll}
\di v_t-u_{1x}=0,\\
\di u_{1t}+p_x=\f{4\v}{3}(\f {\mu(\t)}vu_{1x})_{x}-\int\xi_1^2\Pi_{1x}d\xi,\\
\di u_{it}=\v(\f{\mu(\t)}{v}u_{ix})_x-\int\xi_1\xi_i\Pi_{1x}d\xi,
~i=2,3,\\
\di\bigl(\t+\f{|u|^{2}}{2}\bigr)_{t}+
(pu_1)_{x}=\v(\f{\k(\t)}{v}\t_x)_x+\f{4\v}{3}(\f{\mu(\t)}{v}u_1u_{1x})_x\\
\di\qquad+\v\sum_{i=2}^3(\f{\mu(\t)}{v} u_iu_{ix})_x
-\int\f12\xi_1|\xi|^2\Pi_{1x}d\xi.
\end{array}
\right. \label{(1.22)}
\end{equation}

The Riemann problem for  the Euler system \eqref{(1.16)} is an initial value problem with initial data
\begin{equation*}
(v,u,\t)(t=0,x)=\left\{
\begin{array}{l}
(v_-,u_{-},\t_-),~~~x<0,\\
(v_+,u_{+},\t_+),~~~x>0,
\end{array}
\right.
\end{equation*}
where, $u=(u_1,u_2,u_3)$,  $u_{\pm}=(u_{1\pm},0,0)$ and $v_\pm>0,u_{1\pm},\t_\pm>0$ are
 constants. It is known that the generic solution
to the Riemann problem consists of three waves that propagates at
different speeds, that is, shock, rarefaction wave and contact
discontinuity, cf. \cite{C-F, Lax}. We denote this solution  by
$(\wt{V},\wt{U},\wt{\Theta})(t,x)$. Note that
$\tilde{U}=(\tilde{U}_1,0,0)$.

Given the right end state $(v_+,u_{1+},\t_+)$,  the following wave
curves for the left end state  $(v,u_1,\t)$ in the phase space
are defined with $v>0$ and
$\t>0$ for the Euler equations \eqref{(1.16)}.

$\bullet$ Contact discontinuity curve:
\begin{equation}
CD(v_+,u_{1+},\t_+)= \{(v,u_1,\t)  |  u_1=u_{1+}, p=p_+, v
\not\equiv v_+
 \}. \label{(2.1)}
\end{equation}

$\bullet$ $i$-Rarefaction wave curve $(i=1,3)$:
\begin{equation}
 R_i (v_+, u_{1+}, \theta_+):=\Bigg{ \{} (v, u_1, \theta)\Bigg{ |}v<v_+ ,~u_1=u_{1+}-\int^v_{v_+}
 \lambda_i(\eta,
s_+) \,d\eta,~ s(v, \theta)=s_+\Bigg{ \}},\label{(2.2)}
\end{equation}
where $s_+=s(v_+,\t_+)$ and $\l_i=\l_i(v,s)$ is the $i$-th
characteristic speed of  \eqref{(1.16)}.

$\bullet$ $i$-Shock wave curve $(i=1,3)$:
\begin{equation}
 S_i (v_+, u_{1+}, \theta_+):=\Bigg{ \{} (v, u_1, \theta)\Bigg{ |}
\begin{array}{ll}
-s_i(v_+-v)-(u_{1+}-u_1)=0,\\
-s_i(u_{1+}-u_1)+(p_+-p)=0,\\
-s_i(E_+-E)+(p_+u_{1+}-pu_1)=0,
\end{array}
{\rm and}~~\l_{i+}<s_i<\l_{i-}\Bigg{ \}},\label{3-shock-curve}
\end{equation}
where $E=\t+\f{|u|^2}{2}, p=\f{2\t}{3v},E_+=\t_++\f{|u_+|^2}{2},
p_+=\f{2\t_+}{3v_+}$ , $\l_{i\pm}=\l_i(v_\pm,\t_\pm)$ and $s_i$ is
the $i-$shock speed.

For definiteness, we consider the case when the solution to the Riemann problem
is a superposition of a 1-rarefaction and a 3-shock
wave with a contact discontinuity in between, that is,
$(v_-,u_{1-},\t_-) \in$ $R_1$-$CD$-$S_3(v_+,u_{1+},\t_+)$. Then
there exist uniquely two intermediate states $(v_*, u_{1*},\t_*)$
and $(v^*, u_1^*,\t^*)$ such that $(v_-, u_{1-},\t_-)\in R_1(v_*,
u_{1*},\t_*)$, $(v_*, u_{1*},\t_*)\in CD(v^*, u_1^*,\t^*)$  and
$(v^*, u_1^*,\t^*)\in S_3(v_+,u_{1+},\t_+)$.

Hence,  the wave pattern $(\wt V,\wt U,\wt{\mathcal{E}})(t,x)$
can be written as
\begin{equation}
\begin{array}{l}
 \left(\begin{array}{cc} \wt V\\ \wt U_1 \\
\wt{\mathcal{E}}
\end{array}
\right)(t, x)= \left(\begin{array}{cc}v^{r_1}+ v^{cd}+ v^{s_3}\\ u_1^{r_1}+ u_1^{cd}+ u_1^{s_3} \\
E^{r_1}+ E^{cd}+ E^{s_3}
\end{array}
\right)(t, x) -\left(\begin{array}{cc} v_*+v^*\\ u_{1*}+u_1^*\\
E_*+E^*
\end{array}
\right), \quad \di \wt U_i=0,(i=2,3),
\end{array}
\label{(2.38)}
\end{equation}
where $(v^{r_1}, u_1^{r_1}, \t^{r_1} )(t,x)$ is the 1-rarefaction
wave defined in \eqref{(2.2)} with the right state $(v_+, u_{1+},
\t_+)$ given by $(v_*, u_{1*}, \theta_* )$, $(v^{cd}, u_1^{cd},
\t^{cd} )(t,x)$ is the contact discontinuity defined in
\eqref{(2.1)} with the states $(v_-, u_{1-}, \t_-)$ and $(v_+,
u_{1+}, \t_+)$ given by $(v_*, u_{1*}, \theta_* )$ and $(v^*,
u_{1}^*, \theta^* )$ respectively, and $(v^{s_3}, u_1^{s_3},
\t^{s_3})(t,x)$ is the 3-shock wave defined in \eqref{3-shock-curve}
with the left state $(v_-, u_{1-}, \t_-)$ given by $(v^*, u_1^*,
\theta^* )$.

Consequently, we can define
\begin{equation}
\wt\T(t,x)=(\wt{\mathcal{E}}(t,x)-\f{\wt U(t,x)^2}{2}).
\label{superposition-1}
\end{equation}

Due to the singularity of the rarefaction wave at $t=0$, in this
paper, we consider the problem
 in the time interval $[h,T]$ for any small fixed $h>0$ up to any
arbitrarily  fixed time $T>0$. To investigate the
interaction between the waves and the initial layer is another
interesting topic that will not be discussed here. With the
above preparation, the main result can be stated as follows.

\begin{theorem}\label{Theorem 2.2}
Let  $(\wt V,\wt U,\wt\T)(t,x)$  be a Riemann solution to the Euler
equations which is a superposition of a 1-rarefaction wave, a
2-contact discontinuity and a 3-shock wave, and
$\d=|(v_+-v_-,u_+-u_-,\t_+-\t_-)|$ be the wave strength. There exist
a small positive constant $\delta_0$, and a global Maxwellian
$\mb{M}_\star=\mb{M}_{[v_\star,u_\star,\t_\star]}$ such that if the
wave strength satisfies $\d\leq\d_0$, then in any time interval
$[h,T]$ with $0<h<T$, there exists a positive constant
$\v_0=\v_0(\d,h,T)$, such that if the Knudsen number $\v\leq \v_0$,
then the Boltzmann equation admits a family of smooth solutions
$f^{\v,h}(t,x,\xi)$ satisfying
\begin{equation*}
\sup_{(t,x)\in\Sigma_{h,T}}\|f^{\v,h}(t,x,\xi)-\mb{M}_{[\wt V,\wt
U,\wt\T]}(t,x,\xi)\|_{L_\xi^2(\f{1}{\sqrt{\mb{M}_\star}})}\leq
C_{h,T}~ \v^\f15|\ln\v|,
\end{equation*}
where $\Sigma_{h,T}=\{(t,x)| h\leq t\leq T, |x|\geq h, |x-s_3t|\geq
h\}$, the norm $\|\cdot\|_{L_\xi^2(\f{1}{\sqrt{\mb{M}_\star}})}$ is
 $\|\f{\cdot}{\sqrt{\mb{M}_\star}}\|_{L_\xi^2(\mb{R}^3)}$ and the positive constant $C_{h,T}$  depends on  $h$ and $T$  but
is independent of $\v$. Consequently,  when $ \v\rightarrow 0+$ and
then $h\rightarrow 0+, T\rightarrow+\i$, we have
\begin{equation*}
\|f^{\v,h}(\xi)-\mb{M}_{[\wt V,\wt
U,\wt\T]}(\xi)\|_{L_\xi^2(\f{1}{\sqrt{\mb{M}_\star}})}(t,x)\rightarrow
0,~~a.e.~{\rm in}~~\mathbf{R}^+\times \mathbf{R}.
\end{equation*}
\end{theorem}

\begin{remark} Theorem \ref{Theorem 2.2}
shows that  away from the initial time
$t=0$, the contact discontinuity  at $x=0$ and the shock
discontinuity at $x=s_3t$, for small total wave strength $\d\leq \delta_0$
and
Knudsen number $\v\leq\v_0$,  there exists a
family of smooth solutions $f^{\v,h}(t,x,\xi)$ of the Boltzmann
equation  which tends to the Maxwellian $\mb{M}_{[\wt
V,\wt U,\wt\T]}(t,x,\xi)$ with $(\wt V,\wt U,\wt\T)(t,x)$ being the
Riemann solution to the Euler equations  as
 a superposition of a 1-rarefaction wave, a 2-contact discontinuity and
a 3-shock wave when $\v\rightarrow 0$ with a convergence
rate $\v^{\f{1}{5}}|\ln\v|$. Note that this superposition of waves is the most generic case
for the Riemann problem. Similar results hold for any other superpositions of
waves by using the same analysis.
\end{remark}

\begin{remark}
The proof of the above theorem crucially depends on the introduction
of two kinds of hyperbolic waves. The hyperbolic wave I was
constructed by Huang-Wang-Yang \cite{Huang-Wang-Yang-3} for the
compressible Navier-Stokes equations to recover the viscous terms to
the inviscid approximation of rarefaction wave pattern where the
rarefaction wave structure plays an important role in the
construction.

 The
hyperbolic wave II is constructed to remove the error terms due to the
viscous contact wave approximation. Note that
the construction of the hyperbolic wave II can not be done simply around
the contact wave approximation as the hyperbolic wave I for the rarefaction wave.
Otherwise, the wave interaction terms thus induced will lead to insufficiently
decay in term of the Kundsen number.
Instead, it is constructed around the
superposition of the approximate 1-rarefaction wave, the hyperbolic
wave I, the 2-viscous contact wave and  the 3-shock profile as a whole. Moreover, it also takes care of  the non-conservative  terms in the previous
reduced system so that   energy
estimates can be taken for  anti-derivative of the perturbation.
\end{remark}

\begin{remark}
Note that the analysis  can also be applied to the
vanishing viscosity limit of the one dimensional compressible Navier-Stokes equations. In fact, the vanishing viscosity
limit of the one dimensional  compressible Navier-Stokes equations in
some sense can be
viewed as a special case of hydrodynamic limit of Boltzmann
equation to the Euler equations by neglecting the
microscopic effect.
\end{remark}

\begin{remark}
If the total wave strength $\d=|(v_+-v_-,u_+-u_-,\t_+-\t_-)|\leq
\d_0,$ then from the wave curves defined in \eqref{(2.1)},
\eqref{(2.2)} and \eqref{3-shock-curve},  we know that  $\d^{R_1},
\d^{CD},\d^{S_3}\leq C\d_0$ where $\d^{R_1}, \d^{CD},\d^{S_3}$ are
the wave strengths of rarefaction wave, contact discontinuity and
shock wave, respectively.
\end{remark}

Let us now review some previous works on the hydrodynamic limits to
the Boltzmann equation.  For the case when the Euler equations have
smooth solutions, the vanishing Knudsen number limit of the Boltzmann
equation has been studied even in the case with an initial layer,
cf. Caflisch \cite{Caflisch}, Lachowicz \cite{Lachowicz}, Nishida
\cite{Nishida} and Ukai-Asona \cite{Ukai-Asona} etc. However, as
well-known, solutions of the Euler equations
 in general develop  singularities, such as shock waves and contact
discontinuities. Therefore,  how to verify the hydrodynamic limit
from the Boltzmann equation to the Euler equations with basic wave
patterns becomes a natural problem in the process to the general
setting. In this direction,  with slab symmetry, as mentioned
earlier, there were studies on each individual wave pattern. For
superposition of different types of waves, to our knowledge, there
is only one result given in \cite{Huang-Wang-Yang-2} about the
superposition of two rarefaction waves and one contact
discontinuity.

On the other hand, for the incompressible equations, there are
works, such as those by Bardos-Golse-Levermore,
 Bardos-Levermore-Ukai-Yang, Bardos-Ukai,
Golse-Saint Raymond, Levermore-Masmoudi and Sone which studied
direct derivations of the incompressible Navier-Stokes equations in
the long time scaling, about which more is known, cf. \cite{BGL, BU,
BGLY,GS, LM, sone, sone-2} and the references therein. In
particular, Golse and Saint-Raymond
 showed that the limits of suitably rescaled
sequences of the  DiPerna-Lions renormalized solutions to the
Boltzmann equation are the Leray solutions to the incompressible
Navier-Stokes equations. However, even in this aspect, the
uniqueness and regularity of the solution are still big issues.
Since we will concentrate on the compressible Euler limit in this
paper, we will not go into details about the incompressible limits.

Furthermore, the Boltzmann equation provides more information than
the classical fluid dynamical systems
 so that it describes some phenomena
which can not be modeled by using the classical systems, such as
  Euler and Navier-Stokes equations. This kind of interesting
phenomena, such as the thermal creep flow in a rarefied gas was
known since the time of Maxwell. Some mathematical formulations and
numerical computations
 on the basis of kinetic equations
were studied since 1960s, cf. the works by Sone \cite{sone,sone-2}.
 However, the justification of this kind of
fluid dynamics is almost open with rigorous mathematical
theory.

Finally, we briefly outline the proof of the theorem. Firstly, we
define the individual wave profile. Then we
introduce the first family of hyperbolic wave by linearizing around
the approximate rarefaction profile and by adding the viscosity and heat
conductivity terms induced by the  profile. Then we
define the first approximation of the superposition of this hyperbolic wave
together with the three basic wave patterns so that it   takes care
of the hyperbolicity of the rarefaction wave in the viscous setting.

Based on this,  we linearize the fluid system around this profile and consider
the propagation of the extra error due to the contact discontinuity
approximation and then define a second set of hyperbolic wave. By
adding these two sets of hyperbolic waves to the superposition of
the three basic wave profiles, we will perform the energy estimate on the
Boltzmann equation with suitable initial data through  the
macro-micro decomposition. Precisely, for the macroscopic component,
we will consider the anti-derivative of the perturbation after
applying a hyperbolic scaling. By using the dissipation in the
fluid-type system and the linearized Boltzmann operator
 on the microscopic component, we
can close the energy estimate through a suitable chosen a priori
assumption. Then the statements in the theorem follow.

The rest of the paper will be arranged as follows. In Section 2, we will
construct the approximate solutions to the Boltzmann equation
 corresponding to the basic wave patterns to the Euler
system. Then we obtain the detailed information on the difference
between the Riemann solution to Euler system and
 the
approximate solution  to the Boltzmann equation by the construction.
In Section 3, we will construct a family of solutions to the
Boltzmann equation around the approximate solution by using energy
method to close the a priori estimate. Since the proofs of two
Propositions \ref{Prop3.1} and \ref{Prop3.2} about
 the lower and higher order energy estimates respectively are very technical and long, we put them to the
Appendices.

\

\noindent\textbf{Notations:} Throughout this paper, the positive
generic constants which are independent of $\v, T, h$ are denoted by
$c,C,C_i(i=1,2,3,\cdots)$, while $C_{h,T}$ represents a generic
positive constant depending on $h$ and $T$ but  independent
of $\v$. And we will use $\|\cdot\|$ to denote the standard
$L_2(\mathbf{R};dy)$ norm, and $\|\cdot\|_{H^i}~(i=1,2,3,\cdots)$ to
denote the Sobolev $H^i(\mathbf{R};dy)$ norm. Sometimes, we also use
$O(1)$ to denote a uniform bounded constant which is independent of
$\v, T, h$.

\section{Approximate Wave Patterns}
\setcounter{equation}{0}

In this section, we will construct the approximate wave profile that
consists of three basic wave patterns and two hyperbolic waves.
For this,
 we will firstly recall the construction of the approximate
rarefaction wave for the Boltzmann equation. Then we will introduce the
hyperbolic waves I to correct the error terms coming from the
hyperbolic approximation.
 Then we will construct the viscous contact wave to Boltzmann
equation and study the non-conservative error terms.  The
 viscous shock profile  to the Boltzmann
equation will then be recalled. With the above wave patterns,
 we will introduce the hyperbolic wave II to take care of
the error terms due to the viscous contact wave by avoiding the
interaction between the viscous contact wave with the wave patterns
defined earlier.

\subsection{Rarefaction Wave}

For the rarefaction wave, since there is no exact
rarefaction wave profile for either the Navier-Stokes equations or
the Boltzmann equation, the following approximate rarefaction wave
profile satisfying the Euler equations was introduced in \cite{MN-86,Xin1}. For the completeness of the presentation, we
include its definition and the properties obtained in the above
two papers as follows.

If $(v_-, u_{1-}, \theta_-) \in R_1 (v_+, u_{1+}, \theta_+)$, then
there exists  a $1$-rarefaction wave $(v^{r_1}, u_1^{r_1},
E^{r_1})(x/t)$ which is a global  solution to the following Riemann
problem
\begin{eqnarray}
\left\{
\begin{array}{l}
\di  v_{t}- u_{1x}= 0,\\
\di u_{1t} +  p_{x} = 0 ,
\\
\di E_t + (pu_1)_x =0,\\
\di  (v, u_1, \t)(t=0,x)=\left\{
\begin{array}{l}
\di (v_-, u_{1-}, \t_-),   x<
0 ,\\
\di  (v_+, u_{1+}, \t_+), x> 0 .
\end{array}
\right.
\end{array} \right.\label{(2.8)}
\end{eqnarray}
Consider the following inviscid Burgers equation with Riemann data
\begin{equation}
\left\{
\begin{array}{l}
w_t+ww_x=0,\\
w(t=0,x)=\left\{
\begin{array}{ll}
w_-,&x<0,\\
w_+,&x>0.
\end{array}
 \right.
  \end{array}
 \right.\label{(2.9)}
\end{equation}
If $w_-<w_+$, then the above Riemann problem admits a rarefaction
wave solution
\begin{equation}
w^r(t,x)=w^r(\f xt)=\left\{
\begin{array}{ll}
w_-,&\f xt\leq w_-,\\
\f xt,&w_-\leq \f xt\leq w_+,\\
w_+,&\f xt\geq w_+.
\end{array}
\right.\label{(2.10)}
\end{equation}

As in \cite{Xin1},  the approximate rarefaction wave $(V^{R_1},
U^{R_1}, \Theta^{R_1})(t,x)$ to the problem (\ref{(2.8)}) can be
constructed by the solution of the Burgers equation
\begin{eqnarray}
\left\{
\begin{array}{l}
\di w_{t}+ww_{x}=0,\\
\di w( 0,x
)=w_\s(x)=w(\f{x}{\s})=\f{w_++w_-}{2}+\f{w_+-w_-}{2}\tanh\f{x}{\s},
\end{array}
\right.\label{(2.11)}
\end{eqnarray}
where $\s>0$ is a small parameter
 to be determined later to be $\v^\f15$. Note that the solution $w^r_\s(t,x)$ of the
problem (\ref{(2.11)}) is given by
$$
w^r_\s(t,x)=w_\s(x_0(t,x)),\qquad x=x_0(t,x)+w_\s(x_0(t,x))t.
$$




The smooth approximate rarefaction wave profile denoted by
$(V^{R_1}, U^{R_1}, \T^{R_1})(t,x)$ can be defined  by
\begin{eqnarray}
\left\{
\begin{array}{l}
\di  S^{R_1}(t,x)=s(V^{R_1}(t,x),\T^{R_1}(t,x))=s_+,\\
\di w_\pm=\l_{1\pm}:=\l_1(v_\pm,\t_\pm), \\
\di w_\s^r(t,x)= \l_1(V^{R_1}(t,x),s_+),\\
\di U^{R_1}_1(t,x)=u_{1+}-\int^{V^{R_1}(t,x)}_{v_+}  \l_1(v,s_+)
dv,\\
\di U^{R_1}_i(t,x)\equiv0,~i=2,3.
\end{array} \right.\label{(2.12)}
\end{eqnarray}
Note that $(V^{R_1}, U^{R_1}, \Theta^{R_1})(t,x)$ defined above
satisfies
\begin{eqnarray}
\begin{cases}
  V^{R_1}_t-U^{R_1} _{1x} = 0,     \cr
    U^{R_1}_{1t}+P^{R_1}_x
    =0,\cr
U^{R_1}_{it}=0,~i=2,3,\cr
    \mathcal{E}^{R_1}_t
        + (P^{R_1} U^{R_1}_1)_x
    =0,
\end{cases}\label{rarefaction-equ}
\end{eqnarray}
where $ P^{R_1}=p( V^{R_1}, \T ^{R_1})=\f{2\T^{R_1}}{3V^{R_1}}$ and
$\mathcal{E}^{R_1}=\T^{R_1}+\f{|U^{R_1}|^2}{2}$. The properties of
the rarefaction wave profile can be summarized as follows.

\begin{lemma}\label{Lemma 2.3}(\cite{Xin1}) The approximate rarefaction waves $(V^{R_1},
U^{R_1}, \Theta^{R_1})(t,x)$ constructed in \eqref{(2.12)} have the
following properties:
\begin{enumerate}
\item[(1)] $U^{R_1}_{1x}(t,x)>0$ for $x\in \mathbf{R}$, $t>0$;
\item[(2)] For any $1\leq p\leq +\i,$ the following estimates holds,
$$
\begin{array}{ll}
\|(V^{R_1},U^{R_1}_1, \Theta^{R_1})_x\|_{L^p(dx)} \leq
C\min\big{\{}\d^{R_1}\s^{-1+ 1/p},~
(\d^{R_1})^{1/p}t^{-1+1/p}\big{\}},\\
\|(V^{R_1},U^{R_1}_1, \Theta^{R_1})_{xx}\|_{L^p(dx)} \leq
C\min\big{\{}\d^{R_1}\s^{-2+ 1/p},~ \s^{-1+1/p}t^{-1}\big{\}},\\
\end{array}
$$
where the positive constant $C$  depends only on $p$ and the wave
strength;
\item[(3)] If $x\geq \l^{R_1}_{1+}t$, then
$$
\begin{array}{l}
 |(V^{R_1},U^{R_1},\T^{R_1})(t,x)-(v_+,u_+,\t_+)|\leq Ce^{-\f{2|x-\l_{1+}t|}{\s}},\\[2mm]
 |\partial^k_x(V^{R_1},U^{R_1},\T^{R_1})(t,x)|\leq
\f{C}{\s^k}e^{-\f{2|x-\l_{1+}t|}{\s}},~k=1,2;
\end{array}
$$

\item[(4)] There exist positive constants $C$ and $\s_0$ such that
for $\s\in(0,\s_0)$ and $t>0,$
$$
\sup_{x\in\mathbf{R}}|(V^{R_1},U^{R_1},
\mathcal{E}^{R_1})(t,x)-(v^{r_1},u^{r_1}, E^{r_1})(\f xt)|\leq
\f{C}{t}[\s\ln(1+t)+\s|\ln\s|].
$$
\end{enumerate}
\end{lemma}

\subsection{Hyperbolic Wave I}

Since the whole wave profile consisits of a shock wave whose rate of
change in the shock region is of the order of $\v^{-1}$, we have to consider
the anti-derivative of the perturbation in order to cope with the correct sign
as in the stability analysis. From
\eqref{rarefaction-equ}, we know that the approximate rarefaction
wave $(V^{R_1},U^{R_1},\T^{R_1})(t,x)$ satisfies the compressible
Euler equations exactly without viscous terms. Thus if we carry out
the energy estimates to the anti-derivative variables, the error
terms due to the viscous terms from the approximate rarefaction wave
are not good enough to get the desired estimates. In order to
overcome this difficulty, we introduce the
hyperbolic wave I to recover these viscous terms.

 This hyperbolic wave denoted by  $(d_1,d_2,d_3)(t,x)$ can be defined as follows.
Consider a
linear system
\begin{equation}
\left\{
\begin{array}{l}
\di  d_{1t}-d_{2x}=0,\\
\di
d_{2t}+(p^{R_1}_vd_1+p^{R_1}_{u_1}d_2+p^{R_1}_Ed_3)_x=\f43\v(\f{\mu(\T^{R_1})U^{R_1}_{1x}}{V^{R_1}})_x,
\\\di d_{3t}+[(pu_1)^{R_1}_vd_1+(pu_1)^{R_1}_{u_1}d_2+(pu_1)^{R_1}_Ed_3]_x=
\v(\f{\k(\T^{R_1})\T^{R_1}_x}{V^{R_1}})_x+\f43\v(\f{\mu(\T^{R_1})U^{R_1}_1U^{R_1}_{1x}}{V^{R_1}})_x,
\end{array} \right.
\label{hyp}
\end{equation}
where $p=\f{R\t}{v}=p(v,u,E)=\f{2E-u^2}{3v}$ and
$p^{R_1}_v=p_v(V^{R_1}, U^{R_1}, \mathcal{E}^{R_1})$ etc.
Note that the left hand side of the above system is the linearization
of the Euler equation around the rarefaction wave approximation.  We
want to solve this linear hyperbolic system \eqref{hyp} on the time
interval $[h,T]$. Firstly, we diagonalize the above system by rewriting it as
\begin{equation}
\left(
\begin{array}{l}
\di d_1\\
\di d_2\\
\di d_3
\end{array}
\right)_t +\left[A^{R_1}\left(
\begin{array}{l}
\di d_1\\
\di d_2\\
\di d_3
\end{array}
\right)\right]_x=\left(
\begin{array}{c}
\di 0\\
\di H^{R_1}_1\\
\di H^{R_1}_2
\end{array}
\right),
 \label{hyp-1}
\end{equation}
where $H^{R_1}_1=\v(\f{\mu(\T^{R_1})U^{R_1}_{1x}}{V^{R_1}})_x,
H^{R_1}_2=\v(\f{\k(\T^{R_1})\T^{R_1}_x}{V^{R_1}})_x+\v(\f{\mu(\T^{R_1})U^{R_1}_1U^{R_1}_{1x}}{V^{R_1}})_x$.
Here,
 the matrix
$$
A^{R_1}=\left(
\begin{array}{ccc}
\di 0&-1&0\\
\di p^{R_1}_v&p^{R_1}_{u_1}&p^{R_1}_E\\
\di (pu_1)^{R_1}_v&(pu_1)^{R_1}_{u_1}&(pu_1)^{R_1}_E
\end{array}
\right)
$$
has three distinct eigenvalues
$\l_1^{R_1}:=\l_1(V^{R_1},s_\pm)<0\equiv\l^{R_1}_2<\l_3(V^{R_1},s_\pm):=\l_3^{R_1}$
and the corresponding left and right eigenvectors denoted
$l_j^{R_1},r_j^{R_1}~(j=1,2,3)$ respectively, satisfy
$$
L^{R_1}A^{R_1}R^{R_1}={\rm diag(\l_1^{R_1},0,\l_3^{R_1})}\equiv
\Lambda^{R_1},\qquad
L^{R_1}R^{R_1}={\rm Id.},
$$
Here $L^{R_1}=(l_1^{R_1},l_2^{R_1},l_3^{R_1})^t,
R^{R_1}=(r_1^{R_1},r_2^{R_1},r_3^{R_1})$ with
$l_i^{R_1}=l_i(V^{R_1}, U^{R_1}_1, s_+)$ and $
r_i^{R_1}=r_i(V^{R_1}, U^{R_1}_1, s_+)~(i=1,2,3)$ and ${\rm Id.}$ is
the $3\times3$ identity matrix. Now we set
\begin{equation}
(D_1,D_2,D_3)^t=L^{R_1}(d_1,d_2,d_3)^t.\label{D}
\end{equation}
Then
\begin{equation}
(d_1,d_2,d_3)^t=R^{R_1}(D_1,D_2,D_3)^t,\label{D-d}
\end{equation}
and $(D_1,D_2,D_3)$ satisfies the system
\begin{equation}
\begin{array}{rr}
 \left(
\begin{array}{l}
\di D_1\\
\di D_2\\
\di D_3
\end{array}
\right)_t +\left[\Lambda^{R_1}\left(
\begin{array}{l}
\di D_1\\
\di D_2\\
\di D_3
\end{array}
\right)\right]_x=L^{R_1}\left(
\begin{array}{c}
\di 0\\
\di H^{R_1}_1\\
\di H^{R_1}_2
\end{array}
\right) +L^{R_1}_tR^{R_1}\left(
\begin{array}{l}
\di D_1\\
\di D_2\\
\di D_3
\end{array}
\right)+L^{R_1}_xR^{R_1}\Lambda^{R_1}\left(
\begin{array}{l}
\di D_1\\
\di D_2\\
\di D_3
\end{array}
\right).
\end{array}
 \label{hyp-1-1}
\end{equation}
Due to the fact that the $1-$Riemann invariant is constant along the approximate
rarefaction wave curve, we have
$$
L^{R_1}_t=-\l_1^{R_1}L^{R_1}_x.
$$
Substituting the above equation into \eqref{hyp-1-1}, we obtain the
diagonalized system
\begin{equation}
\left\{
\begin{array}{l}
\di  D_{1t}+(\l_1^{R_1}D_1)_x=b^{R_1}_{12}H^{R_1}_1+b^{R_1}_{13}H^{R_1}_2+a^{R_1}_{12}V_x^{R_1}D_2+a^{R_1}_{13}V_x^{R_1}D_3,\\
\di
D_{2t}=b^{R_1}_{22}H^{R_1}_1+b^{R_1}_{23}H^{R_1}_2+a^{R_1}_{22}V_x^{R_1}D_2+a^{R_1}_{23}V_x^{R_1}D_3,
\\\di D_{3t}+(\l_3^{R_1}D_3)_x=b^{R_1}_{32}H^{R_1}_1+b^{R_1}_{33}H^{R_1}_2+a^{R_1}_{32}V_x^{R_1}D_2+a^{R_1}_{33}V_x^{R_1}D_3,
\end{array} \right.
\label{hyp-2}
\end{equation}
where $a^{R_1}_{ij}, b^{R_1}_{ij}$ are some given functions of
$V^{R_1},U_1^{R_1}$ and $S^{R_1}=s_\pm$. Note that in the
diagonalized system \eqref{hyp-2}, the equations of $D_2,D_3$ are
decoupled from $D_1$ because of the property of the rarefaction wave.

Now we impose the following boundary condition to the above linear
hyperbolic system \eqref{hyp-2} in the domain $(t,x)\in
[h,T]\times\mathbf{R}$:
\begin{equation}
D_1(t=h,x)=0,\qquad D_2(t=T,x)=D_3(t=T,x)=0.\label{cond}
\end{equation}
With this boundary condition,  we can solve the linear diagonalized hyperbolic system
(\ref{hyp-2}) under the conditions (\ref{cond}). Moreover, we have
the following estimates on the solution.

\

\begin{lemma}\label{hw1} There exists a positive constant $C_{h,T}$
independent of $\v$ such that
\begin{itemize}
\item[(1)] $$\|\f{\partial^k}{\partial x^k}d_i(t,\cdot)\|^2_{L^2(dx)} \leq
C_{h,T}~ \f{\v^2}{\s^{2k+1}},\quad i=1,2,3,~k=0,1,2,3.$$
\item[(2)] If $x>\l_{1+}t$, then we have
$$
|d_i(x,t)|\leq C_{h,T}~ \f{1}{\s}e^{- \f {|x-\l_{1+}t|}\s}, \qquad
|d_{ix}(x,t)|\leq C_{h,T}~ \f{1}{\s^2}e^{- \f
{|x-\l_{1+}t|}\s},~i=1,2,3.
$$
\end{itemize}
\end{lemma}
The proof of Lemma \ref{hw1} can be done similarly as in
\cite{Huang-Wang-Yang-3} for the compressible Navier-Stokes
equations.

\subsection{Viscous Contact Wave}
In this subsection, we construct the contact wave
$(V^{CD},U^{CD},\T^{CD})(t,x)$ for the Boltzmann equation motivated
by \cite{Huang-Yang}. Consider the Euler system \eqref{(1.16)} with
a Riemann initial data
\begin{equation}
(v,u,\t)(t=0,x)=\left\{
\begin{array}{l}
(v_-,u_{-},\t_-),~~~x<0,\\
(v_+,u_{+},\t_+),~~~x>0,
\end{array}
\right. \label{(2.26)}
\end{equation}
where $u_{\pm}=(u_{1\pm},0,0)$ and $v_\pm>0,\t_\pm>0,u_{1\pm}$ are
given constants. It is known (cf. \cite{Smoller}) that the Riemann
problem (\ref{(1.16)}), (\ref{(2.26)}) admits a contact
discontinuity
\begin{equation}
(v^{cd},u^{cd},\t^{cd})(t,x)=\left\{
\begin{array}{l}
(v_-,u_{-},\t_-),~~~x<0,\\
(v_+,u_{+},\t_+),~~~x>0,
\end{array}
\right. \label{(2.27)}
\end{equation}
provided that $(v_-,u_{1-},\t_-)\in CD(v_+,u_{1+},\t_+)$, that is,
\begin{equation}
u_{1+}=u_{1-},\qquad p_-:=\f{2\t_-}{3v_-}=p_+:=\f{2\t_+}{3v_+}.
\label{(2.28)}
\end{equation}
Then for the Navier-Stokes equations, by the energy equation
$(\ref{(1.22)})_4$
and the mass equation $(\ref{(1.22)})_1$ with
$v\approx\f{2\t}{3p_+}$, cf. \cite{Huang-Wang-Yang-2}, we can
obtain the following nonlinear diffusion equation
\begin{equation}
\t_t=\v(a(\t)\t_x)_x,~~~a(\t)=\f{9 p_+\k(\t)}{10\t}. \label{(2.30)}
\end{equation}
From \cite{Atkinson-Peletier} and \cite{Duyn-Peletier}, we know that
the nonlinear diffusion equation (\ref{(2.30)}) admits a unique
self-similar solution $\hat{\T}(\eta),~\eta=\f{x}{\sqrt{\v(1+t)}}$
satisfying the  boundary conditions $
\hat{\T}(\pm\i,t)=\t_\pm. $ Let $\delta^{CD}=|\t_+-\t_-|$, then
$\hat{\T}(t,x)$ has the property that
\begin{equation}
\hat{\T}_x(t,x)=\f{O(1)\delta^{CD}}{\sqrt{\v(1+t)}}e^{-\f{cx^2}{\v(1+t)}}, \label{(2.31)}
\end{equation}
with some positive constant $c$  depending only on  $\t_{\pm}$.

Correspondingly, we can define the Navier-Stokes profile by
\begin{equation}
\begin{array}{ll}
\di \hat V=\f{2}{3p_+}\hat\T,\\
\di \hat{U}_1=u_{1+}+\f{2\v a(\hat{\T})}{3p_+}\hat{\T}_x,~~~\hat
U_i=0,i=2,3.
\end{array}
 \label{NS-CD}
\end{equation}
For the Boltzmann equation, if we still use the above Navier-Stokes
profile $(\hat V,\hat U,\hat \T)$, we can not get any decay
with respect to the Knudsen number $\v$ due to  the non-fluid component.
Hence, we construct a Boltzmann
contact wave as follows. Set
\begin{equation}
\mb{G}^{CD}(t,x,\xi)=\f{3\v}{2v\t}\mb{L}^{-1}_\mb{M}\big\{\mb{P}_1[\xi_1(\f{|\xi-u|^2}{2\t}\T^{CD}_x+\xi\cdot{U}^{CD}_{x})\mb{M}]\big\},
\label{G-CD}
\end{equation}
and
\begin{equation}
\Pi^{CD}_{11}=\mb{L}_\mb{M}^{-1}\big[\v(-\f{u_1}v\mb{G}^{CD}_x+\f{1}{v}\mb{P}_1(\xi_1\mb{G}^{CD}_x))-Q(\mb{G}^{CD},\mb{G}^{CD})\big],
 \label{Pi-CD-1}
\end{equation}
where $(V^{CD},U^{CD},\T^{CD})(t,x)$ is the viscous contact wave for
the Boltzmann equation to be constructed later.

Note that for the Boltzmann equation, the leading terms in
 the energy equation $(\ref{(1.22)})_4$ can be written as
\begin{equation}
\t_t+p_+u_{1x}=\v(\f{\k(\t)\t_x}{v})_x-\int\f12\x_1|\x|^2\Pi^{CD}_{11x}d\x+u_{1+}\int\f12\x_1^2\Pi^{CD}_{11x}d\x.\label{B-CD}
\end{equation}
By the definition of $\Pi^{CD}_{11}$ in \eqref{Pi-CD-1}, we have
\begin{equation}
-\int\f12\x_1|\x|^2\Pi^{CD}_{11}d\x+u_{1+}\int\f12\x_1^2\Pi^{CD}_{11}d\x=\Delta_{11}+\Delta_{12},
\end{equation}
where
\begin{equation}
\Delta_{11}=\v^2\Big[g_{11}\t_x\T^{CD}_x+g_{12}v_x\T^{CD}_x+g_{13}(\T^{CD}_x)^2+g_{14}\T^{CD}_{xx}\Big],\label{delta-41}
\end{equation}
with $g_{1i}=g_{1i}(v,u,\t),(i=1,2,3,4)$ being smooth functions of
$(v,u,\t)$, and
\begin{equation}
\Delta_{12}=O(1)\v^2\Big[\big(|v_x|+|u_x|+|\t_x|+|\T^{CD}_x|+|U^{CD}_x|\big)|U^{CD}_x|+|u_x||\T^{CD}_x|+|U^{CD}_{xx}|\Big].
 \label{delta-42}
\end{equation}
Thus, by choosing the leading term and dropping the higher order
term $\Delta_{12x}$  in \eqref{B-CD}, we have
\begin{equation}
\t_t=\v(a(\t)\t_x)_x+\f35\Delta_{11x},
 \label{B-Diffu}
\end{equation}
where $a(\t)$ is defined in \eqref{(2.30)} and $\Delta_{11}$ is
defined in \eqref{delta-41}. To represent
 the microscopic effect on the wave profile, we want to define
$\T^{CD}$ to be close to $ \hat\T(\f{x}{\sqrt{\v(1+t)}})+\hat\T^{nf}(t,x)$ with
$\hat\T$ being determined by \eqref{(2.30)}, \eqref{(2.31)} and
$\hat\T^{nf}$ represents the part of the nonlinear diffusion wave
coming from the non-fluid component not appearing in the
Navier-Stokes level. Moreover, the term $\hat\T^{nf}$ decays faster
than $\hat\T$ so that it can be viewed as the perturbation around
the Navier-Stokes profile $\hat\T$. To construct $\hat\T^{nf}$, we
linearize the equation \eqref{B-Diffu} around the Navier-Stokes
profile $\hat\T$ and drop all the higher order terms. This leads to
a linear diffusion equation for $\hat\T^{nf}$
\begin{equation}
\hat\T^{nf}_t=\v(a(\hat\T)\hat\T^{nf}_x)_x+\v(a^\prime(\hat\T)\hat\T_x\hat\T^{nf})_x+\f35\wt\Delta_{11x},\label{B-Diffu2}
\end{equation}
where $\wt\Delta_{11}=\v^2(\tilde g_{11}+\f{2}{3p_+} \tilde
g_{12}+\tilde g_{13})(\hat\T_x)^2+\v^2\tilde g_{14}\hat\T_{xx}$ with
$\tilde g_{1i}=\tilde g_{1i}(\hat V,\hat U,\hat \T)~(i=1,2,3,4)$.
Integrating \eqref{B-Diffu2} with respect to $x$ yields that
\begin{equation}
\Xi_{1t}=\v a(\hat\T)\Xi_{1xx}+\v
a^\prime(\hat\T)\hat\T_x\Xi_{1x}+\f35\wt\Delta_{11},
\end{equation}
where
\begin{equation}
\Xi_1(t,x)=\int_{-\i}^x\hat\T^{nf}(t,x)dx.
\end{equation}
Note that $\wt\Delta_{11}$ takes the form of
$\f{\v}{1+t}A^1(\f{x}{\sqrt{\v(1+t)}})$ and satisfies that
$$
|\wt\Delta_{11}|=O(\d^{CD})\v(1+t)^{-1}e^{-\f{x^2}{4a(\t_\pm)\v
(1+t)}}, ~~~{\rm as}~~x\rightarrow\pm\i.
$$
We can check that there exists a self-similar solution
$\Xi_1(\f{x}{\sqrt{\v(1+t)}})$ for \eqref{B-Diffu2} with the
boundary conditions $\Xi_1(-\i)=0, \Xi_1(+\i)=\Xi_{1+}$. Here
$\Xi_{1+}$ can be any given constant satisfying
$|\Xi_{1+}|<\d^{CD}$. It is worthy to point out that even though the
function $\Xi_1(t,x)$ depends on the constant $\Xi_{1+}$,
$\hat\T^{nf}(t,x)=\Xi_{1x}(t,x)\rightarrow 0$ as
$x\rightarrow\pm\i.$ That is, the choice of the constant $\Xi_{1+}$
has no influence on the ansantz as long as $|\Xi_{1+}|<\d^{CD}$.
From now on, we fix $\Xi_{1+}$ so that the function $\Xi_1(t,x)$ is
uniquely determined and its derivative $\Xi_{1x}=\hat\T^{nf}$
has the property
\begin{equation}
|\hat\T^{nf}|=|\Xi_{1x}|=O(\d^{CD})\v^{\f12}(1+t)^{-\f12}e^{-\f{x^2}{4a(\t_\pm)\v(1+t)}},~~~{\rm
as}~~x\rightarrow\pm\i.\label{T-nf}
\end{equation}
Then we apply the similar procedure to construct the second and the
third components of the velocity of the contact wave denoted by
$U^{CD}_i~(i=2,3)$ as follows. The leading part of the equation for $u_i$ in
$\eqref{(1.22)}_i~(i=2,3)$ is
\begin{equation}
u_{it}=\v\big(\f{3p_+\mu(\t)}{2\t}u_{ix}\big)_x-\int\x_1\x_i\Pi^{CD}_{11x}d\x.
\end{equation}
Firstly,  we have
\begin{equation}
-\int\x_1\x_i\Pi^{CD}_{11}d\x=\Delta_{i1}+\Delta_{i2},
\end{equation}
where
\begin{equation}
\Delta_{i1}=\v^2\Big[g_{i1}\t_x\T^{CD}_x+g_{i2}v_x\T^{CD}_x+g_{i3}(\T^{CD}_x)^2+g_{i4}\T^{CD}_{xx}\Big],\label{delta-i1}
\end{equation}
with $g_{ij},(i=2,3,~j=1,2,3,4)$ being the smooth functions of
$(v,u,\t)$ and
\begin{equation}
\Delta_{i2}=O(1)\v^2\Big[\big(|v_x|+|u_x|+|\t_x|+|\T^{CD}_x|+|U^{CD}_x|\big)|U^{CD}_x|+|u_x||\T^{CD}_x|+|U^{CD}_{xx}|\Big].
 \label{delta-i2}
\end{equation}
Thus we expect the viscous contact wave $U^{CD}_i~(i=2,3)$ to satisfy
the following linear equation
\begin{equation}
U^{CD}_{it}=\v\big(\f{3p_+\mu(\hat\T)}{2\hat\T}U^{CD}_{ix}\big)_x+\wt\Delta_{i1x},~~~i=2,3,\label{U-CD-i}
\end{equation}
where where $\wt\Delta_{i1}=\v^2(\tilde g_{i1}+\f{2}{3p_+} \tilde
g_{i2}+\tilde g_{i3})(\hat\T_x)^2+\v^2\tilde g_{i4}\hat\T_{xx}$ with
$\tilde g_{ij}=\tilde g_{ij}(\hat V,\hat U,\hat
\T)~(i=2,3,~j=1,2,3,4)$. Integrating \eqref{U-CD-i} with respect to
$x$ yields that
\begin{equation}
\Xi_{it}=\v \f{3p_+\mu(\hat\T)}{2\hat\T}\Xi_{ixx}+\wt\Delta_{i1},
\end{equation}
where
\begin{equation}
\Xi_i(t,x)=\int_{-\i}^xU^{CD}_i(t,x)dx.
\end{equation}
Note that $\wt\Delta_{i1}$ takes the form
$\f{\v}{1+t}A^i(\f{x}{\sqrt{\v(1+t)}}), i=2,3$ and satisfies that
$$
|\wt\Delta_{i1}|=O(\d^{CD})\v(1+t)^{-1}e^{-\f{x^2}{4a(\t_\pm)\v
(1+t)}}, ~~~{\rm as}~~x\rightarrow\pm\i,
$$

We can check that there exists a self-similar solution
$\Xi_i(\f{x}{\sqrt{\v(1+t)}})$ for \eqref{U-CD-i} with the boundary
conditions $\Xi_i(-\i)=0, \Xi_i(+\i)=\Xi_{i+},~(i=2,3)$. Again, here
$\Xi_{i+}$ can be any given constant satisfying
$|\Xi_{i+}|<\d^{CD}$. As we explained before, the choice of the
constant $\Xi_{i+}$ has no influence on the ansantz as long as
$|\Xi_{1+}|<\d^{CD}$. We fix $\Xi_{i+}$ so that the function
$\Xi_i(t,x)$ is uniquely determined and the derivative
$\Xi_{ix}=U^{CD}_i$ has the property
\begin{equation}
|U^{CD}_i|=|\Xi_{ix}|=O(\d^{CD})\v^{\f12}(1+t)^{-\f12}e^{-\f{x^2}{4b(\t_\pm)\v(1+t)}},~~~{\rm
as}~~x\rightarrow\pm\i,
\end{equation}
with $b(\t_\pm)=\max\{a(\t_\pm),\f{3p_+\mu(\t_\pm)}{2\t_\pm})\}$.

In summary,   the viscous contact wave $(V^{CD},U^{CD},\T^{CD})(t,x)$ can be
defined by
\begin{equation}
\begin{array}{ll}
\di V^{CD}=\f{2}{3p_+}(\hat{\T}+\hat\T^{nf}),\\
\di U^{CD}_1=u_{1+}+\f{2
}{3p_+}\Big[\v a(\hat{\T})\hat{\T}_x+\v(a(\hat\T)\hat\T^{nf})_x+\f35\wt\Delta_{11}\Big],\\
U^{CD}_i=\Xi_{ix},~~(i=2,3),\\[3mm]
\di\T^{CD}=\hat{\T}+\hat\T^{nf}+H,
\end{array}
 \label{contact-wave}
\end{equation}
where
\begin{equation}\label{H11}
H=-\f{\v
V^{CD}}{p_+}\Big[a(\hat\T)\hat\T_t+(a(\hat\T)\hat\T^{nf})_t\Big]+\f{4\v^2\mu(\hat\T)}{3p_+}\Big[(a(\hat\T)\hat\T_x)_x
+(a(\hat\T)\hat\T^{nf})_{xx}\Big]-\f32V^{CD}\int\x_1^2\wt\Pi_{11}d\x,
\end{equation}
is chosen such that the momentum equation that the viscous contact
wave satisfies has an error term with sufficient decay in $\v$.
Without $H,$ the error term in the momentum equation decay like
$\v^{\f12}$. In order to get $\v$ order decay in the error term, we
should introduce the higher order approximate term $H$ in the
definition of $\T^{CD}$ in the contact wave. Here,
$\di\int\x_1^2\wt\Pi_{11}d\x$ in \eqref{H11} is the corresponding
function defined in \eqref{Pi-CD-1} by replacing both the variables
$(v,u,\t)$ and $(V^{CD},U^{CD},\T^{CD})$ by $(\hat V,\hat
U,\hat\T)$, and it satisfies
\begin{equation}
\int\x_1^2\wt\Pi_{11}d\x=O(1)\v^2|(\hat\T_x^2,\hat\T_{xx})|.\label{hat-Pi-11}
\end{equation}
Hence, by \eqref{(2.31)}, \eqref{T-nf} and \eqref{hat-Pi-11}, we
have
\begin{equation}
H=O(\d^{CD})\v(1+t)^{-2}e^{-\f{cx^2}{\v(1+t)}},~~~{\rm
as}~~x\rightarrow\pm\i. \label{H}
\end{equation}

Now the contact wave $(V^{CD},U^{CD},\T^{CD})(t,x)$ defined in
\eqref{contact-wave} satisfies the following system
\begin{equation}
\left\{\begin{array}{llll}
\di V^{CD}_t-U^{CD}_{1x}=0,\\
\di U^{CD}_{1t}+P^{CD}_x=\f{4\v}{3}(\f{\mu(\T^{CD})}{V^{CD}}U^{CD}_{1x})_x-\int\x_1^2\Pi^{CD}_{11x}d\x+Q^{CD}_1,\\
\di
U^{CD}_{it}=\v(\f{\mu(\T^{CD})}{V^{CD}}U^{CD}_{ix})_x-\int\x_1\x_i\Pi^{CD}_{11x}d\x+Q^{CD}_i,
i=2,3,\\
\di\mathcal{E}^{CD}_{t}+
(P^{CD}U^{CD}_{1})_{x}=\v(\f{\k(\T^{CD})}{V^{CD}}\T^{CD}_x)_x+\f{4\v}{3}(\f{\mu(\T^{CD})U^{CD}_{1}U^{CD}_{1x}}{V^{CD}})_x\\
\di\qquad+\sum_{i=2}^3\v(\f{\mu(\T^{CD})U^{CD}_iU^{CD}_{ix}}{V^{CD}})_x-\int\x_1\f{|\x|^2}{2}\Pi^{CD}_{11x}d\x+Q^{CD}_4,
\end{array}
\right. \label{(2.33)}
\end{equation}
where $P^{CD}=\f{2\T^{CD}}{3V^{CD}}$,
$\mathcal{E}^{CD}=\T^{CD}+\f{|U^{CD}|^2}{2}$,
$\di-\int\x_1\x_i\Pi^{CD}_{11x}d\x~(i=1,2,3)$ and $\di
-\int\x_1\f{|\x|^2}{2}\Pi^{CD}_{11x}d\x$ are the corresponding
functions defined in \eqref{Pi-CD-1} by replacing the variables
$(v,u,\t)$ by $(V^{CD},U^{CD},\T^{CD})$, respectively. Moreover,
\begin{equation}
\begin{array}{ll}
\di
Q^{CD}_1=\f{2\v}{5p_+}\wt\Delta_{11t}-\f{8\v^2}{15p_+}\big(\f{\mu(\hat\T)\wt\Delta_{11x}}{V^{CD}}\big)_x
-\f{4\v}{3}\Big(\f{\mu(\T^{CD})-\mu(\hat\T)}{V^{CD}}U^{CD}_{1x}\Big)_x+\int\x_1^2(\Pi_{11}^{CD}-\wt\Pi_{11})_xd\x\\
\qquad~\di =
O(1)\delta^{CD}\v(1+t)^{-2}e^{-\f{cx^2}{\v(1+t)}},\qquad{\rm as}~~
x\rightarrow\pm\i,
\end{array}
\label{Q-CD-1}
\end{equation}
\begin{equation}
\begin{array}{ll}
\di Q^{CD}_i&\di=
-\v\Big(\f{\mu(\T^{CD})-\mu(\hat\T)}{V^{CD}}U^{CD}_{ix}\Big)_x+\int\x_1\x_i(\Pi_{i1}^{CD}-\wt\Pi_{i1})_xd\x+\wt\Delta_{i2x}\\
&\di = O(1)\delta^{CD}\v(1+t)^{-2}e^{-\f{cx^2}{\v(1+t)}},\qquad{\rm
as}~~ x\rightarrow\pm\i,~~~~i=2,3,
\end{array}
\label{Q-CD-i}
\end{equation}
and
\begin{equation}
\begin{array}{ll}
 Q^{CD}_4&\di=-\f{5\v}{3}\Big[\Big(a(\T^{CD})-a(\hat\T)-a^\prime(\hat\T)(\hat\T^{CD}-\hat\T)\Big)(\hat\T_x+\hat\T^{nf}_x)+a(\T^{CD})H_x\\
 &\di~~+a^\prime(\hat\T)(\hat\T^{CD}-\hat\T)\hat\T^{nf}_x\Big]_x+\f{3p_+\v}{2}\Big(\f{\k(\T^{CD})\T^{CD}_xH}{\T^{CD}(\T^{CD}-H)}\Big)_x
 +H_t\\
 &\di~~+\f{2U^{CD}_{1x}H}{3V^{CD}}-\f{4\v\mu(\T^{CD})}{3V^{CD}}
(U^{CD}_{1x})^2-\sum_{i=2}^3\f{\v\mu(\T^{CD})}{V^{CD}}(U^{CD}_{ix})^2\\
&\di\quad+\int\x_1\f{|\x|^2}{2}(\Pi_{11}^{CD}-\wt\Pi_{11})_xd\x+(U^{CD}_1-u_{1+})\int\x_1^2(\Pi_{11}^{CD}-\wt\Pi_{11})_xd\x\\
&\di\quad+\sum_{i=2}^3U^{CD}_i\int\x_1\x_i(\Pi_{11}^{CD}-\wt\Pi_{11})_x d\x+\sum_{i=1}^3U^{CD}_iQ^{CD}_i+\wt\Delta_{12x}\\
&\di=O(1)\delta^{CD}\v(1+t)^{-2}e^{-\f{cx^2}{\v(1+t)}},\qquad{\rm
as}~~ x\rightarrow\pm\i,
\end{array}
\label{Q-CD-4}
\end{equation}
with some positive constant $c>0$  depending only on $\t_\pm$ and
$\wt\Delta_{i2x}~(i=1,2,3)$ being the corresponding functions defined in
\eqref{delta-42} and \eqref{delta-i2} by replacing both $(v,u,\t)$
and $(V^{CD},U^{CD},\T^{CD})$ by $(\hat V,\hat U,\hat\T)$.

Note that from (\ref{(2.31)}), we have
\begin{equation}
|(V^{CD},U^{CD},\T^{CD})(t,x)-(v^{cd},u^{cd},\t^{cd})(t,x)|= O(1)\delta^{CD} e^{-\f{cx^2}{2\v(1+t)}}.\\
\label{(2.37)}
\end{equation}

\subsection{Shock Profile}
In this subsection,  we will firstly recall the shock profile $F^{S_3}(x-\bar s_3t,\xi)$ of the
Boltzmann equation (\ref{(B)}) in Eulerian coordinates with its existence and
properties given in the papers by Caflisch-Nicolaenko \cite{Caflish-Nicolaenko} and Liu-Yu
\cite{Liu-Yu}, \cite{Liu-Yu-1}.  And then we will state the corresponding
properties in the Lagrangian coordinates used in this paper.

First of all,
$F^{S_3}(x-\bar s_3t,\xi)$ satisfies
\begin{equation}
\left\{
\begin{array}{ll}
\di -\bar s_3(F^{S_3})^\prime+\x_1(F^{S_3})^\prime=\f{1}{\v}Q(F^{S_3},F^{S_3}),\\[5mm]
F^{S_3}(\pm\i,\xi)=\mathbf{M}_{\pm}(\xi):=\mathbf{M}_{[\r_\pm,u_\pm,\t_\pm]}(\xi),
\end{array}
\right.\label{BS}
\end{equation}
 where $^\prime=\f{d}{d\vartheta}$, $\vartheta=x-\bar s_3t,$ $u_\pm=(u_{1\pm},0,0)$
 and
$(\r_{\pm},u_{\pm},\t_\pm)$ satisfy Rankine-Hugoniot condition
\begin{equation}
\left\{
\begin{array}{ll}
\di-\bar s_3(\r_+-\r_-)+(\r_+u_{1+}-\r_-u_{1-})=0,\\
\di -\bar
s_3(\r_+u_{1+}-\r_-u_{1-})+(\r_+u_{1+}^2+p_+-\r_-u_{1-}^2-p_-)=0,\\
\di -\bar
s_3(\r_+E_+-\r_-E_-)+(\r_+u_{1+}E_++p_+u_{1+}-\r_-u_{1-}E_--p_-u_{1-})=0,
\end{array}
\right. \label{RH-E}
\end{equation}
and Lax entropy condition
\begin{equation}
\l_{3+}^E<\bar s_3<\l_{3-}^{E},
 \label{Lax-E}
\end{equation}
with $\bar s_3$ being 3-shock wave speed and
$\l_3^E=u_1+\f{\sqrt{10\t}}{3}$ being the third characteristic
eigenvalue of the Euler equations in the Eulerian coordinate and
$\l_{3\pm}^E=u_{1\pm}+\f{\sqrt{10\t_\pm}}{3}$.

 By the macro-micro decomposition  around the local Maxwellian
$\mb{M}^{S_3}$, set
$$
F^{S_3}(x,t,\x)=\mathbf{M}^{S_3}(x,t,\x)+\mathbf{G}^{S_3}(x,t,\x),
$$
where
$$
\mathbf{M}^{S_3}(x,t,\x)=\mathbf{M}_{[\r^{S_3},u^{S_3},\t^{S_3}]}(x,t,\x)=\f{\r^{S_3}(x,t)}{\sqrt{(2\pi
R\t^{S_3}(x,t))^3}}e^{-\f{|\x-u^{S_3}(x,t)|^2}{2R\t^{S_3}(x,t)}},
$$
with
\begin{equation}
\left(
\begin{array}{c}
 \r^{S_3}\\
 \rho^{S_3} u^{S_3}_i\\
 \rho^{S_3}(\t^{S_3}+\frac{|u^{S_3}|^2}{2})
\end{array}
\right)=\int_{\mathbf{R}^3}\left(
\begin{array}{c}
 1\\
 \x_i\\
\frac{|\x|^2}{2}
\end{array}
\right)F^{S_3}(x,t,\xi)d\xi,~~i=1,2,3.\label{shock-fluid}
\end{equation}

With respect to the inner product
$\langle\cdot,\cdot\rangle_{\mathbf{M}^{S_3}}$ defined in
\eqref{product},  we can now define the macroscopic projection
$\mb{P}^{S_3}_0$ and microscopic projection $\mb{P}^{S_3}_1$ by
\begin{equation}
 \mb{P}^{S_3}_0g = {\di\sum_{j=0}^4\langle g,\chi^{S_3}_j\rangle_{\mathbf{M}^{S_3}}\chi^{S_3}_j,}
 \qquad
 \mb{P}^{S_3}_1g= g-\mb{P}^{S_3}_0g,
\label{PS3}
\end{equation}
where $\chi^{S_3}_j~(j=0,1,2,3,4)$ are the corresponding pairwise
orthogonal base defined in \eqref{orthogonal-base} by replacing
$(\r,u,\t, \mb{M})$ by $(\r^{S_3},u^{S_3},\t^{S_3}, \mb{M}^{S_3})$.

Under the above macro-micro decomposition, the solution
$F^{S_3}=F^{S_3}(x-\bar s_3t,\xi)$ satisfies
$$
\mb{P}^{S_3}_0F^{S_3}=\mb{M}^{S_3},~\mb{P}^{S_3}_1F^{S_3}=\mb{G}^{S_3},
$$
and the Boltzmann equation (\ref{BS}) becomes
\begin{equation}
(\mb{M}^{S_3}+\mb{G}^{S_3})_t+\xi_1(\mb{M}^{S_3}+\mb{G}^{S_3})_x
=\f{1}{\v}[2Q(\mb{M}^{S_3},\mb{G}^{S_3})+Q(\mb{G}^{S_3},\mb{G}^{S_3})].
\label{(1.6)}
\end{equation}
Correspondingly, we have the following fluid-type system for the
fluid components of shock profile:
\begin{equation}
\begin{array}{l}
\left\{
\begin{array}{l}
\di \r^{S_3}_{t}+(\r^{S_3} u^{S_3}_1)_x=0,\\
\di (\r^{S_3} u^{S_3}_1)_t+[\r^{S_3} (u^{S_3}_1)^2
+p^{S_3}]_x=\f{4\v}{3}\big(\mu(\t^{S_3})
u^{S_3}_{1x}\big)_x-\int\xi_1^2\Pi^{S_3}_xd\xi,  \\
\di (\r^{S_3} u^{S_3}_i)_t+(\r^{S_3}
u^{S_3}_1u^{S_3}_i)_x=\v\big(\mu(\t^{S_3})
u^{S_3}_{ix}\big)_x-\int\xi_1\xi_i\Pi^{S_3}_xd\xi,~ i=2,3,\\
\di [\rho^{S_3}(\t^{S_3}+\f{|u^{S_3}|^2}{2})]_t+[\r^{S_3}
u^{S_3}_1(\t^{S_3}+\f{|u^{S_3}|^2}{2})+p^{S_3}u^{S_3}_1]_x=\v\big(\k(\t^{S_3})\t^{S_3}_x\big)_x\\
\di\qquad
+\f{4\v}{3}\big(\mu(\t^{S_3})u^{S_3}_1u^{S_3}_{1x}\big)_x+\v\sum_{i=2}^3\big(\mu(\t^{S_3})u^{S_3}_iu^{S_3}_{ix}\big)_x
-\int\f12\xi_1|\xi|^2\Pi^{S_3}_xd\xi.
\end{array}
\right.
\end{array}\label{shock-Euler}
\end{equation}
In fact, from the invariance of the equation \eqref{BS} by changing
$\x_i$ with $-\x_i$ and the fact that $u_{i\pm}=0,$ we have $\di
u^{S_3}_i=\int\xi_1\xi_i\Pi^{S_3}_xd\xi\equiv0$ for $i=2,3.$

And the  equation for the non-fluid component
${\mb{G}}^{S_3}$ is
$$
\mb{G}^{S_3}_t+\mb{P}^{S_3}_1(\xi_1\mb{M}^{S_3}_x)+\mb{P}^{S_3}_1(\xi_1\mb{G}^{S_3}_x)
=\f{1}{\v}\left[\mb{L}_{\mb{M}^{S_3}}\mb{G}^{S_3}+Q(\mb{G}^{S_3}.
\mb{G}^{S_3})\right],
$$
Here $\mb{L}_{\mb{M}^{S_3}}$ is the linearized collision operator of
$Q(F^{S_3},F^{S_3})$
 with respect to the local Maxwellian $\mb{M}^{S_3}$:
$$
\mb{L}_{\mb{M}^{S_3}} g=2Q(\mb{M}^{S_3}, g)=Q(\mb{M}^{S_3}, g)+
Q(g,\mb{M}^{S_3}).
$$
Thus
$$
\begin{array}{l}
\di \mb{G}^{S_3}=\v \mb{L}_{\mb{M}^{S_3}}^{-1}[\mb{P}^{S_3}_1(\xi_1\mb{M}^{S_3}_x)] +\Pi^{S_3},\\[6mm]
\di \Pi^{S_3}=\mb{L}_{\mb{M}^{S_3}}^{-1}[\v(\mb{G}^{S_3}_t+\mb
{P}^{S_3}_1(\xi_1\mb{G}^{S_3}_x))-Q(\mb{G}^{S_3}, \mb{G}^{S_3})].
\end{array}
$$
Now we recall  the properties of the shock profile $F^{S_3}(x-\bar
s_3t,\x)$ that are given or can be induced by  Liu-Yu in Theorem
6.8, \cite{Liu-Yu-1}.

\begin{lemma}\label{Lemma-shock} (\cite{Liu-Yu-1})  If the shock wave strength $\d^{S_3}$ is small enough,
then the Boltzmann equation \eqref{(B)} admits a 3-shock profile
solution $F^{S_3}(x-\bar s_3t,\x)$ uniquely up to a shift satisfying
the following properties:
\begin{itemize}
\item [(1)] The shock profile converges to its far fields
exponentially fast with an exponent proportional to the magnitude of
the shock wave strength, that is
$$
\left\{
\begin{array}{ll}
\di |(\r^{S_3}-\r_{\pm},u^{S_3}_1-u_{1\pm},\t^{S_3}-\t_\pm)|\leq
C\d^{S_3} e^{-c_\pm\f{\d^{S_3}|\vartheta|}{\v}},\quad{\rm
as}~~\vartheta\rightarrow \pm\i,\\
\di\big(\int
\f{\nu(|\x|)|\mb{G}^{S_3}|^2}{\mb{M}_0}d\x\big)^{\f12}\leq
C(\d^{S_3})^2 e^{-c_\pm\f{\d^{S_3}|\vartheta|}{\v}},\quad{\rm
as}~~\vartheta\rightarrow \pm\i,
\end{array}
\right.
$$
with  $\d^{S_3}$ being the 3-shock strength and  $\mb{M}_0$ being
the global Maxwellian which is close to the shock profile
with its precise definition given in Theorem 6.8, \cite{Liu-Yu-1}.

\item[(2)] Compressibility of 3-shock profile:
$$
(\l^{E}_3)_\vartheta< 0,\qquad
\l^{E}_3=u^{S_3}_1+\f{\sqrt{10\t^{S_3}}}{3}.
$$

\item[(3)] The following properties hold:
$$
\begin{array}{ll}
\di \r^{S_3}_\vartheta\sim u^{S_3}_{1\vartheta}\sim
\t^{S_3}_\vartheta\sim(\l^{E}_3)_\vartheta\sim \f{1}{\v}\big(\int
\f{\nu(|\x|)|\mb{G}^{S_3}|^2}{\mb{M}_0}d\x\big)^{\f12},
\end{array}
$$
where $A\sim B$ denotes the equivalence of the quantities $A$ and
$B$, and
$$
\left\{
\begin{array}{ll}
\di u^{S_3}_i\equiv 0,~~\int\x_1\x_i\Pi^{S_3}d\x\equiv0, ~~i=2,3,\\
\di \di |\partial^k_\vartheta(\r^{S_3},u^{S_3}_{1},\t^{S_3})|\leq
C\f{(\d^{S_3})^{k-1}}{\v^{k-1}}|(\r^{S_3}_\vartheta,
u^{S_3}_{1\vartheta},
\t^{S_3}_\vartheta)|,~k\ge 2,\\[3mm]
\di \big(\int
\f{\nu(|\x|)|\partial^k_\vartheta\mb{G}^{S_3}|^2}{\mb{M}_0}d\x\big)^{\f12}\leq
C\f{(\d^{S_3})^k}{\v^k}\big(\int
\f{\nu(|\x|)|\mb{G}^{S_3}|^2}{\mb{M}_0}d\x\big)^{\f12},~~k\ge 1,\\
\\ \di
|\int\x_1\varphi_i(\x)\Pi^{S_3}_{\vartheta}d\x|\leq
C\d^{S_3}|u^{S_3}_{1\vartheta}|,~~i=1,2,3,4,
\end{array}
\right.
$$
where $\varphi_i(\x)~(i=1,2,3,4)$ are the collision invariants
defined in \eqref{collision-invar}.
\end{itemize}
\end{lemma}

Now we rewrite this shock profile in Lagrangian coordinate by using
the transformation \eqref{Lag} and use $(\tilde t,\tilde x)$
for the Lagrangian coordinate to distingish it from the Eulerian
coordinate $(t,x)$ at this moment. Then the shock profile  in Lagrangian coordinate can be written as $\tilde
F^{S_3}(\tilde x-s_3\tilde t,\x)$ with $s_3$  determined by the
3-shock wave curve given in \eqref{3-shock-curve}. First, from the
Rankine-Hugoniot condition in Eulerian and Lagrangian coordinates, we have
$$
-\bar s_3(\r_+-\r_-)+(\r_+u_{1+}-\r_-u_{1-})=0,
$$
and
$$
-s_3(v_+-v_-)-(u_{1+}-u_{1-})=0,
$$
with $v_\pm=\f{1}{\r_\pm},$ respectively. Thus we have the following
relation between $\bar s_3$ and $s_3$
\begin{equation}
s_3=\r_\pm(\bar s_3-u_{1\pm}). \label{s-3}
\end{equation}
On the other hand, we have from $\eqref{shock-Euler}_1$ that
\begin{equation}
\r^{S_3}(x-s_3 t)[\bar s_3-u^{S_3}_1(x-s_3 t)]={\rm
const.}=\r_\pm(\bar s_3-u_{1\pm}).
\label{SF}
\end{equation}

Note that
$$
\begin{array}{ll}
\di \tilde x-s_3\tilde t &\di=-\int_0^t\r^{S_3}u^{S_3}_1(-\bar
s_3\tau)d\tau+\int_0^x\r^{S_3}(y-\bar s_3\tau)dy-\int_0^t\r_\pm(\bar
s_3-u_{1\pm}) d\tau\\
&\di=-\int_0^t\r^{S_3}u^{S_3}_1(-\bar
s_3\tau)d\tau+\int_0^x\r^{S_3}(y-\bar s_3\tau)dy-\int_0^t\r^{S_3}(-s_3 \tau)[\bar s_3-u^{S_3}_1(-s_3 \tau)] d\tau\\
&\di=\int_0^x\r^{S_3}(y-\bar s_3\tau)dy-\bar
s_3\int_0^t\r^{S_3}(-s_3\tau)d\tau=\int_0^{x-\bar s_3t}\r^{S_3}(\eta)d\eta,
\end{array}
$$
where in the second equality we have used  \eqref{SF}.

This shows that under the Lagrangian transformation \eqref{Lag}, the
shock profile $F^{S_3}(x-\bar s_3t,\x)$ in Eulerian coordinate can
be exactly transformed to the shock profile $\tilde F^{S_3}(\tilde
x-s_3\tilde t,\x)$ in Lagrangian coordinate. Moreover, we have
$$
\tilde F^{S_3}_\eta(\eta,\x)=\r^{S_3}F^{S_3}_\vartheta(\vartheta,\x)
$$
with $\eta=\tilde x-s_3\tilde t$.

For simplicity of the notations, from now on, we  use $(t,x)$ to denote the
Lagrangian coordinate and $F^{S_3}(\eta,\x)$
with $\eta=x-s_3t$ to denote the 3-shock profile of Boltzmann
equation in Lagrangian coordinate. And in the Lagrangian
coordinate,  we have the following Lemma.

\begin{lemma}\label{Lemma-shock-L} Assume that $(v_-,u_-,\t_-)\in
S_3(v_+,u_+,\t_+)$,
 then there exists a unique shock profile
$F^{S_3}(\eta,\x)$ with $\eta=x-s_3t$ up to a shift, to the Boltzmann
equation \eqref{Lag-B} in Lagrangian coordinate. Moreover, there are
positive constants $c_\pm$ and $C$ such that for
$\eta\in\mathbf{R}$,
$$
\left\{
\begin{array}{lll}
\displaystyle s_3 V^{S_3}_{\eta}=-U^{S_3}_{1\eta}>0,\\
\di U^{S_3}_i\equiv 0,~~\int\x_1\x_i\Pi^{S_3}_1d\x\equiv0, ~~i=2,3,\\
\di (|V^{S_3}-v_{\pm}|,|U^{S_3}_1-u_{1\pm}|,|\T^{S_3}-\t_\pm|)\le C\d^{S_3} e^{-\f{c_\pm\d^{S_3}|\eta|}{\v}},\quad{\rm as}~~\eta\rightarrow \pm\i,\\
\di\big(\int
\f{\nu(|\x|)|\mb{G}^{S_3}|^2}{\mb{M}_0}d\x\big)^{\f12}\leq
C(\d^{S_3})^2 e^{-c_\pm\f{\d^{S_3}|\eta|}{\v}},\quad{\rm
as}~~\eta\rightarrow \pm\i.
\end{array}
\right.$$
Furthermore, we have
$$
\begin{array}{ll}
\di V^{S_3}_\eta\sim U^{S_3}_{1\eta}\sim \T^{S_3}_\eta\sim
\f{1}{\v}\big(\int
\f{\nu(|\x|)|\mb{G}^{S_3}|^2}{\mb{M}_0}d\x\big)^{\f12},
\end{array}
$$
and
$$
\begin{array}{ll}
\di \di |\partial^k_\eta(V^{S_3},U^{S_3}_{1},\T^{S_3})|\leq
C\f{(\d^{S_3})^{k-1}}{\v^{k-1}}|(V^{S_3}_\eta, U^{S_3}_{1\eta},
\T^{S_3}_\eta)|,~~k\ge 2,\\[3mm]
\di \big(\int
\f{\nu(|\x|)|\partial^k_\eta\mb{G}^{S_3}|^2}{\mb{M}_0}d\x\big)^{\f12}\leq
C\f{(\d^{S_3})^k}{\v^k}\big(\int
\f{\nu(|\x|)|\mb{G}^{S_3}|^2}{\mb{M}_0}d\x\big)^{\f12},~~k\ge 1,
\end{array}
$$
and
$$
|\int\x_1\varphi_i(\x)\Pi^{S_3}_{1\eta}d\x|\leq
C\d^{S_3}|U^{S_3}_{1\eta}|,~~i=1,2,3,4,
$$
with $\varphi_i(\x)$ being the collision invariants.
\end{lemma}

Furthermore, we have
\begin{equation}
 \left\{
\begin{array}{lll}
\displaystyle V^{S_3}_t-U_{1x}^{S_3}=0,\\
\displaystyle
U_{1t}^{S_3}+P^{S_3}_x=\f43\v(\frac{\mu(\T^{S_3})U_{1x}^{S_3}}{V^{S_3}})_x-\int\x_1^2\Pi^{S_3}_{1x}d\x,\\
\di U_{it}^{S_3}=\v(\frac{\mu(\T^{S_3})U_{ix}^{S_3}}{V^{S_3}})_x-\int\x_1\x_i\Pi^{S_3}_{1x}d\x,~i=2,3,\\
 \displaystyle
\mathcal{E}^{S_3}_t+(P^{S_3}U_1^{S_3})_x=\v(\f{\k(\T^{S_3})\T^{S_3}_x}{V^{S_3}})_x+\f43\v(\frac{\mu(\T^{S_3})U_1^{S_3}U_{1x}^{S_3}}{V^{S_3}})_x\\
\di\qquad\qquad\qquad\quad+\v\sum_{i=2}^3(\frac{\mu(\T^{S_3})U_i^{S_3}U_{ix}^{S_3}}{V^{S_3}})_x
-\int\x_1\f{|\x|^2}{2}\Pi^{S_3}_{1x}d\x,\\
\end{array}
\right.\label{shock-profile}
\end{equation}
 where $\mathcal{E}^{S_3}=\T^{S_3}+\f{|U^{S_3}|^2}{2}$ and
$(v_{\pm},u_{\pm},\t_\pm)$ satisfy Rankine-Hugoniot condition and Lax entropy
condition  and $s_3$ is 3-shock wave speed.

Correspondingly, we have the following equation for the
non-fluid part  of
3-shock profile.
$$
\mb{G}^{S_3}_t-\f{U_1^{S_3}}{V^{S_3}}\mb{G}^{S_3}_x+\f{1}{V^{S_3}}\mb{P}^{S_3}_1(\xi_1\mb{M}^{S_3}_x)
+\f{1}{V^{S_3}}\mb{P}^{S_3}_1(\xi_1\mb{G}^{S_3}_x)
=\f{1}{\v}\left[\mb{L}_{\mb{M}^{S_3}}\mb{G}^{S_3}+Q(\mb{G}^{S_3},
\mb{G}^{S_3})\right].
$$
Here, $\mb{L}_{\mb{M}^{S_3}}$ is the linearized collision operator of
$Q(F^{S_3},F^{S_3})$
 with respect to the local Maxwellian $\mb{M}^{S_3}$:
$$
\mb{L}_{\mb{M}^{S_3}} g=2Q(\mb{M}^{S_3}, g)=Q(\mb{M}^{S_3}, g)+
Q(g,\mb{M}^{S_3}).
$$
Thus
\begin{equation}
\begin{array}{l}
\di \mb{G}^{S_3}=\v \mb{L}_{\mb{M}^{S_3}}^{-1}\Big[\f{1}{V^{S_3}}\mb{P}^{S_3}_1(\xi_1\mb{M}^{S_3}_x)\Big] +\Pi^{S_3}_1,\\[6mm]
\di
\Pi^{S_3}_1=\mb{L}_{\mb{M}^{S_3}}^{-1}\Big[\v\Big(\mb{G}^{S_3}_t-\f{U_1^{S_3}}{V^{S_3}}\mb{G}^{S_3}_x
+\f{1}{V^{S_3}}\mb{P}^{S_3}_1(\xi_1\mb{G}^{S_3}_x)\Big)-Q(\mb{G}^{S_3},\mb{G}^{S_3})\Big].
\end{array}
\label{Pi-S-1}
\end{equation}

\subsection{Hyperbolic Wave II}

The purpose of this subsection is to construct the second
hyperbolic wave. Up to now, we can define the following
approximate composite wave profile $(\bar V,\bar
U,\bar{\mathcal{E}})(t,x)$
\begin{equation}
\begin{array}{ll}
 \left(\begin{array}{cc} \bar V\\ \bar U_1 \\
  \bar{\mathcal{E}}
\end{array}
\right)(t, x)= \left(\begin{array}{cc}  V^{R_1}+d_1+ V^{CD}+ V^{S_3}\\ U_1^{R_1}+d_2+ U_1^{CD}+ U_1^{S_3} \\
\mathcal{E}^{R_1}+d_3+ \mathcal{E}^{CD}+ \mathcal{E}^{S_3}
\end{array}
\right)(t, x) -\left(\begin{array}{cc} v_*+v^*\\ u_{1*}+u_1^{*}\\
E_*+E^*
\end{array}
\right),\\[7mm]
\di\bar U_i=U^{CD}_i,i=2,3,
\end{array}
 \label{(2.39)}
\end{equation}
where $\bar{\mathcal{E}}=\bar\T+\f{|\bar U|^2}{2}$, $(V^{R_1},
U_1^{R_1}, \mathcal{E}^{R_1} )(t,x)$ is the 1-rarefaction wave
defined in \eqref{(2.12)} with the right state $(v_+, u_{1+}, E_+)$
replaced by $(v_*, u_{1*}, E_* )$, $(V^{CD}, U_1^{CD},
\mathcal{E}^{CD} )(t,x)$ is the viscous contact wave defined in
\eqref{contact-wave} with the states $(v_-, u_{1-}, E_-)$ and $(v_+,
u_{1+}, E_+)$ replaced by $(v_*, u_{1*}, E_* )$ and $(v^*, u^*_1,
E^* )$ respectively, and $(V^{S_3}, U_1^{S_3},
\mathcal{E}^{S_3})(t,x)$ is the fluid part of 3-shock profile of
Boltzmann equation defined in \eqref{shock-profile} with the left
state $(v_-, u_{1-}, E_-)$ replaced by $(v^*, u_1^*, E^* )$.

Moreover, we can check that this profile satisfies
\begin{equation}
\left\{
\begin{array}{l}
  \bar V_t- \bar U _{1x} = 0,  \\
   \di   \bar U_{1t}+ \bar P_x
    = \f43\v (\frac{\mu( \bar \T) \bar U_{1x}}{ \bar V}) _{x}-\int\x_1^2\Pi^{CD}_{11x}d\x-\int\x_1^2\Pi^{S_3}_{1x}d\x+\bar{Q}_{1x}+Q^{CD}_1,\\
    \di \bar U_{it}=\v(\frac{\mu( \bar \T) \bar U_{1x}}{ \bar V}) _{x}-\int\x_1\x_i\Pi^{CD}_{11x}d\x-\int\x_1\x_i\Pi^{S_3}_{1x}d\x+\bar Q_{ix}+Q_i^{CD},~~i=2,3,\\
 \di  \bar{\mathcal{E}}_t+ (\bar P\bar U_{1})_x
    =\v( \frac{\k( \bar \T) \bar \T_{x}}{  \bar V})_x + \f43\v(\frac{ \mu( \bar \T)\bar U_{1} \bar U_{1x}}{  \bar
    V})_x+\sum_{i=2}^3\v(\frac{ \mu( \bar \T)\bar U_{i} \bar U_{ix}}{  \bar
    V})_x\\
    \di\hspace{3cm}-\int\x_1\f{|\x|^2}{2}\Pi^{CD}_{11x}d\x-\int\x_1\f{|\x|^2}{2}\Pi^{S_3}_{1x}d\x
     +\bar Q_{4x}+Q^{CD}_4,
\end{array}
\right.\label{super-1}
\end{equation}
where $\bar P =p(\bar V ,\bar\T)$, $Q_i^{CD}~(i=1,2,3,4)$ are
defined in \eqref{Q-CD-1}, \eqref{Q-CD-i} and \eqref{Q-CD-4},
respectively, and
$$\begin{array}{ll}
   \di \bar Q_1=&\di\Big[\bar
   P-(p^{R_1}_vd_1+p^{R_1}_{u_1}d_2+p^{R_1}_Ed_3)-P^{R_1}-P^{CD}-P^{S_3}+p_*+p^*\Big]\\[3mm]
   &\di \di
   -\f43\v\Big[\f{\mu(\bar \T)\bar U_{1x}}{\bar V}-\f{\mu(\T^{R_1})U^{R_1}_{1x}}{V^{R_1}}-\f{\mu(\T^{CD})U^{CD}_{1x}}{V^{CD}}-\f{\mu(\T^{S_3})U^{S_3}_{1x}}{V^{S_3}}\Big],\\[3mm]
   \di \bar Q_i=&\di-\v\Big[\f{\mu(\bar\T)}{\bar
   V}-\f{\mu(\T^{CD})}{V^{CD}}\Big]U^{CD}_{ix},~i=2,3,\\[3mm]
\end{array}$$
$$\begin{array}{ll}
  \di \bar Q_4=&\di
  \Big[\bar P\bar U_{1}-\big((pu_1)^{R_1}_vd_1+(pu_1)^{R_1}_{u_1}d_2+(pu_1)^{R_1}_Ed_3\big)-P^{R_1}U^{R_1}_{1}-P^{CD}U^{CD}_{1}-P^{S_3}U^{S_3}_{1}\\[3mm]
  &\di+p_*u_{1*}+p^*u_1^*\Big]-\v\Big[\f{\k(\bar \T)\bar \T_x}{\bar V}-\f{\k(\T^{R_1})\T^{R_1}_x}{V^{R_1}}-\f{\k(\T^{CD})\T^{CD}_x}{V^{CD}}-\f{\k(\T^{S_3})\T^{S_3}_x}{V^{S_3}}\Big]\\[3mm]
  &\di-\f43\v\Big[\f{\mu(\bar \T)\bar U_{1}\bar U_{1x}}{\bar V}-\f{\mu(\T^{R_1})U^{R_1}_{1}U^{R_1}_{1x}}{ V^{R_1}}-\f{\mu(\T^{CD})U^{CD}_{1}U^{CD}_{1x}}{ V^{CD}}-\f{\mu(\T^{S_3})U^{S_3}_{1}U^{S_3}_{1x}}{
  V^{S_3}}\Big]\\
&\di -\sum_{i=2}^3\v\Big[\f{\mu(\bar \T)}{\bar
V}-\f{\mu(\T^{CD})}{V^{CD}}\Big]U^{CD}_{i} U^{CD}_{ix}  .
\end{array}$$
From \eqref{super-1}, we have that
\begin{equation}
\begin{array}{ll}
 \di \bar \T=&\di \T^{R_1}+d_3+\T^{CD}+\T^{S_3}-(\t_*+\t^*)\\
 &\di-\f12[|\bar U|^2-|U^{R_1}|^2-|U^{CD}|^2-|U^{S_3}|^2+u_{1*}^2+(u^*_{1})^2].
\end{array}
\label{bar-T}
\end{equation}

Direct computation yields
\begin{equation}
\begin{array}{lll}
   \di \bar Q_1&=&\di O(1)\Big[|(V^{R_1}-v_*,\T^{R_1}-\t_*,U^{R_1}_1-u_{1*},d_1,d_2,d_3)||(V^{CD}-v_*,U^{CD}_1-u_{1*},\\
   &&\qquad\quad \di \T^{CD}-\t_*,V^{S_3}-v^*,U^{S_3}_1-u_1^*,\T^{S_3}-\t^*)|\\
&&\di
+|(V^{S_3}-v^*,U^{S_3}_1-u_1^*,\T^{S_3}-\t^*)||(V^{CD}-v^*,U^{CD}_1-u_{1*},U^{CD}_1-u_1^*,\T^{CD}-\t^*)|\\
&&\di
+\v|U^{R_1}_{1x}||(V^{CD}-v_*,\T^{CD}-\t_*,V^{S_3}-v^*,\T^{S_3}-\t^*)|\\
&&\di
+\v|U^{CD}_{1x}||(V^{R_1}-v_*,\T^{R_1}-\T_*,d_1,d_2,d_3,V^{S_3}-v^*,\T^{S_3}-\t^*)|
\\
&&\di+\v|U^{S_3}_{1x}||(V^{R_1}-v_*,\T^{R_1}-\T_*,d_1,d_2,d_3,V^{CD}-v^*,\T^{CD}-\t^*)|\Big]\\
&&\di+O(1)\Big[|(d_1,d_2,d_3)|^2+\v|d_{2x}|+\v|U^{R_1}_{1x}||(d_1,d_2,d_3)|\Big]\\[2mm]
& :=&\di\bar Q_{11}+\bar Q_{12},
\end{array}
\label{barQ1}
\end{equation}
\begin{equation}
\begin{array}{lll}
   \di \bar Q_i= O(1)\v\Big[|(V^{R_1}-v_*,\T^{R_1}-\t_*,U^{R_1}_1-u_{1*},d_1,d_2,d_3, V^{S_3}-v^*,\T^{S_3}-\t^*,U^{S_3}_1-u_{1}^*)||U^{CD}_{ix}|\Big],~i=2,3,\\
  \end{array}
\label{barQi}
\end{equation}
and
\begin{equation}\begin{array}{lll}
   \di \bar Q_4&=&\di O(1)\Big[|(V^{R_1}-v_*,U^{R_1}_1-u_{1*},\T^{R_1}-\t_*,d_1,d_2,d_3)||(V^{CD}-v_*,U^{CD}_1-u_{1*},\\
   &&\qquad\quad \di \T^{CD}-\t_*,V^{S_3}-v^*,U^{S_3}_1-u_1^*,\T^{S_3}-\t^*)|\\
&&\di
+|(V^{S_3}-v^*,U^{S_3}_1-u_1^*,\T^{S_3}-\t^*)||(V^{CD}-v^*,U^{CD}_1-u_1^*,\T^{CD}-\t^*)|\\
&&\di
+\v|(U^{R_1}_{1x},\T^{R_1}_x)||(V^{CD}-v_*,U^{CD}_1-u_{1*},\T^{CD}-\t_*,V^{S_3}-v^*,\T^{S_3}-\t^*)|\\
&&\di
+\v|(U^{CD}_{1x},\T^{CD}_x)||(V^{R_1}-v_*,U^{R_1}_1-u_{1*},\T^{R_1}-\T_*,d_1,d_2,d_3,V^{S_3}-v^*,\T^{S_3}-\t^*)|
\\
&&\di+\v|(U^{S_3}_{1x},\T^{S_3}_x)||(V^{R_1}-v_*,\T^{R_1}-\T_*,d_1,d_2,d_3,V^{CD}-v^*,U^{CD}_1-u_1^*,\T^{CD}-\t^*)|\Big]\\
&&\di+O(1)\Big[|(d_1,d_2,d_3)|^2+\v|(d_{2x},d_{3x})|+\v|(U^{R_1}_{1x},\T^{R_1}_x)||(d_1,d_2,d_3)|\Big]\\[2mm]
& :=&\di\bar Q_{41}+\bar Q_{42}.
\end{array}\label{barQ4}
\end{equation}
Here, $\bar Q_{11}$, $\bar Q_i~(i=2,3)$  and $\bar Q_{41}$ represent
the interaction of waves in different families, $\bar Q_{12}$ and
$\bar Q_{42}$ represent the error terms coming from the approximate
rarefaction wave and the hyperbolic wave I.

Firstly, we estimate the interaction terms $\bar Q_{11}$, $\bar
Q_i~(i=2,3)$ and $\bar Q_{41}$ by dividing the whole domain
$\Omega=\{(t,x)|(t,x)\in[h,T]\times\mathbf{R}\}$ into three regions:
$$
\begin{array}{ll}
\Omega_{R_1}&\di =\{(t,x)\in\Omega|2x\leq \l_{1*}t\},\\
\Omega_{CD}&\di =\{(t,x)\in\Omega|\l_{1*}t<2x<s_3t\},\\
\Omega_{S_3}&\di =\{(t,x)\in\Omega|2x\geq s_3t\},
\end{array}
$$
where $\l_{1*}=\l_1(v_*,\t_*)$ and $s_3$ is the 3-shock speed.

From Lemma \ref{Lemma 2.3}, we have the following estimates in each
region:
\begin{itemize}
\item In $\Omega_{R_1}$,

$\begin{array}{ll}|V^{S_3}-v^*|=O(1)\d^{S_3}e^{-\f{c\d^{S_3}|x-s_3t|}{\v}}=O(1)e^{-\f{C|x|}{\v}}e^{-\f{C_h}{\v}},
\end{array}$

$\begin{array}{ll}|(V^{CD}-v_*,V^{CD}-v^*)|=O(1)\d^{CD}e^{-\f{C|x||\l_{1*}|t}{8\v(1+t)}}e^{-\f{C(\l_{1*}t)^2}{8\v(1+t)}}=O(1)e^{-\f{C_h
|x|}{\v}}e^{-\f{C_h}{\v}};
\end{array}
$
\item In $\Omega_{CD}$,

$\begin{array}{ll}|(V^{R_1}-v_*,d_1,d_2,d_3)|&=O(1)e^{-\f{2|x-\l_{1*}t|}{\s}}=O(1)e^{-\f{C_h
|x|}{\s}}e^{-\f{C_h}{\s}},
\end{array}$

$\begin{array}{ll}|V^{S_3}-v^*|&=O(1)e^{-\f{c\d^{S_3}|x-s_3t|}{\v}}=O(1)e^{-\f{C_h
|x|}{\v}}e^{-\f{C_h}{\v}};
\end{array}$
\item In $\Omega_{S_3}$,

$\begin{array}{ll}|(V^{R_1}-v_*,d_1,d_2,d_3)|&=O(1)e^{-\f{2|x-\l_{1*}t|}{\s}}=O(1)e^{-\f{2
|x|}{\s}}e^{-\f{C_h}{\s}},
\end{array}$

$\begin{array}{ll}|(V^{CD}-v_*,V^{CD}-v^*)|&=O(1)\d^{CD}e^{-\f{C|x||s_3t|}{8\v(1+t)}}e^{-\f{C(s_3t)^2}{8\v(1+t)}}=O(1)e^{-\f{C_h
|x|}{\v}}e^{-\f{C_h}{\v}}.
\end{array}
$
\end{itemize}
Note that we just give the pointwise estimates of $V$component
and $d_i~(i=1,2,3)$ in each region, similar estimates hold also for
the $U_1$ and $\T$ components.  In
summary, we have
\begin{equation}
|(\bar Q_{11},\bar Q_2,\bar Q_3,\bar Q_{41})|=C_{h,T}~e^{-\f{C_h
|x|}{\s}}e^{-\f{C_h}{\s}},\label{Q11}
\end{equation}
with $\s=\v^{\f15}$ and for some positive constants $C_{h,T}$ and
$C_h$  independent of $\v$.

In order to remove the non-conservative error terms
$Q^{CD}_i,(i=1,2,3,4)$ coming from the definition of the viscous contact
wave, we now introduce the following hyperbolic wave
$\vec{b}\triangleq(b_1,b_{21},b_{22},b_{23},b_3)$:
\begin{eqnarray}\label{A2.1}
\begin{cases}
b_{1t}-b_{21x}=0,\\
b_{21t}+[\bar P_v b_1+\bar P_{u_1} b_{21}+\bar P_{u_2} b_{22}+\bar P_{u_3} b_{23}+\bar P_E b_3]_x=-Q^{CD}_1,\\
b_{22t}=-Q^{CD}_2,\\
b_{23t}=-Q^{CD}_3,\\
 b_{3t}+[(\bar P\bar U_1)_v b_1+(\bar P\bar U_1)_{u_1}
b_{21}+(\bar P\bar U_1)_{u_2} b_{22}+(\bar P\bar U_1)_{u_3}
b_{23}+(\bar P\bar U_1)_E b_3]_x=-Q^{CD}_4,
\end{cases}
\end{eqnarray}
where $P=\frac{2\Theta}{3V}=P(V,U,\mathcal{E})$ and $\bar
P_v=P_v(\bar V,\bar U,\bar{\mathcal{E}})$, etc. For later use, we
denote $b_2=(b_{21},b_{22},b_{23})$. Now we want to solve this
linear hyperbolic system \eqref{A2.1} on the interval $[h,T]$. Firstly,
we diagonalize the above system. Rewrite the system \eqref{A2.1} as
\begin{equation}\label{A2.2}
\left(
\begin{array}{l}
\di b_1\\
\di b_{21}\\
\di b_{22}\\
\di b_{23}\\
\di b_3
\end{array}
\right)_t +\left[\bar A(\bar V,\bar U,\bar{\mathcal{E}})\left(
\begin{array}{l}
\di b_1\\
\di b_{21}\\
\di b_{22}\\
\di b_{23}\\
\di b_3
\end{array}
\right)\right]_x=\left(
\begin{array}{c}
\di 0\\
\di -Q^{CD}_1\\
\di -Q^{CD}_2\\
\di -Q^{CD}_3\\
\di -Q^{CD}_4\\
\end{array}
\right),
\end{equation}
where the matrix
$$
\bar A(\bar V,\bar U,\bar{\mathcal{E}})=\left(
\begin{array}{ccccc}
\di 0&-1&0&0&0\\
\di \bar P_v&\bar P_{u_1}&\bar P_{u_2}&\bar P_{u_3}&\bar P_E\\
\di 0&0&0&0&0\\
\di 0&0&0&0&0\\
\di (\bar P\bar U_1)_v&(\bar P\bar U_1)_{u_1}&(\bar P\bar
U_1)_{u_2}&(\bar P\bar U_1)_{u_3}&(\bar P\bar U_1)_E
\end{array}
\right)
$$
has three distinct eigenvalues $\bar \l_1:=\l_1(\bar V,\bar
P)<0=\bar\l_2<\l_3(\bar V,\bar P):=\bar \l_3,$~(here $\bar\l_2$
being 3-repeated eigenvalues) with the corresponding left and right
eigenvectors denoted by
$$\bar l_1,\bar l_{21},\bar l_{22},\bar l_{23},\bar l_3; \quad \bar r_1,\bar r_{21},\bar r_{22},\bar r_{23},\bar r_3.$$ It holds that
$$
\bar L\bar A\bar R={\rm
diag(\bar\l_1,\bar\l_2,\bar\l_2,\bar\l_2,\bar\l_3)}\equiv \bar
\Lambda,
$$
$$
\bar L\bar R={\rm Id.}
$$
Here $\bar L=(\bar l_1,\bar l_{21},\bar l_{22},\bar l_{23},\bar
l_3)^t, \bar R=(\bar r_1,\bar r_{21},\bar r_{22},\bar r_{23},\bar
r_3)$ with $\bar L=\bar L(\bar V, \bar U, \bar{\mathcal{E}})$ and
$\bar R=\bar R(\bar V, \bar U, \bar{\mathcal{E}})~(i=1,2,3)$ and
${\rm Id.}$  is the $5\times5$ identity matrix. Specially, we
can choose
\begin{eqnarray}\label{A2.3}
\bar l_{21}=(\bar P,-\bar U_1,0,0,1),\quad \bar
l_{22}=(0,0,1,0,0),\quad \bar l_{23}=(0,0,0,1,0).
\end{eqnarray}
Set
\begin{equation}\label{A2.4}
\vec{B}\triangleq(B_1,B_{21},B_{22},B_{23},B_3)^t=\bar
L\cdot(b_1,b_{21},b_{22},b_{23},b_3),
\end{equation}
then
\begin{equation}\label{A2.5}
(b_1,b_{21},b_{22},b_{23},b_3)^t=\bar
R\cdot(B_1,B_{21},B_{22},B_{23},B_3)^t,
\end{equation}
and $\vec{B}$ satisfies the system
\begin{equation} \label{A2.6}
\begin{array}{rr}
\di \left(
\begin{array}{l}
\di B_1\\
\di B_{21}\\
\di B_{22}\\
\di B_{23}\\
\di B_3
\end{array}
\right)_t +\left[\bar \Lambda\left(
\begin{array}{l}
\di B_1\\
\di B_{21}\\
\di B_{22}\\
\di B_{23}\\
\di B_3
\end{array}
\right)\right]_x=\bar L\left(
\begin{array}{c}
\di 0\\
\di -Q^{CD}_1\\
\di -Q^{CD}_2\\
\di -Q^{CD}_3\\
\di -Q^{CD}_4
\end{array}
\right) \di +\bar L_t\bar R\left(
\begin{array}{l}
\di B_1\\
\di B_{21}\\
\di B_{22}\\
\di B_{23}\\
\di B_3
\end{array}
\right)+\bar L_x\bar R\bar \Lambda\left(
\begin{array}{l}
\di B_1\\
\di B_{21}\\
\di B_{22}\\
\di B_{23}\\
\di B_3
\end{array}
\right).
\end{array}
\end{equation}
For simplicity,  denote
\begin{equation}
\vec{Q}^{CD}=\big(0\, -Q^{CD}_1, -Q^{CD}_2,
-Q^{CD}_3,-Q^{CD}_4\big)^t.
\end{equation}
So we obtain a diagonalized system
\begin{equation}\label{A2.7}
\left\{
\begin{array}{ll}
\di B_{1t}+(\bar \l_1B_1)_x=\bar
l_{1}\cdot\vec{Q}^{CD}+\sum_{i=1,3}(\bar l_{1t}+\bar \l_i\bar
l_{1x})\cdot
\di  \bar r_iB_{i}+\bar l_{1t}\cdot\sum_{j=1}^{3}\bar{r}_{2j}B_{2j},\\
\di B_{21t}~~ =\bar
l_{21}\cdot\vec{Q}^{CD}+\sum_{i=1,3}(\bar l_{21t}+\bar \l_i\bar
l_{21x})\cdot
 \bar r_iB_i+\bar l_{21t}\cdot\sum_{j=1}^{3}\bar{r}_{2j}B_{2j},\\
\di B_{22t}~~ =\bar l_{22}\cdot\vec{Q}^{CD},\\
\di B_{23t}~~ =\bar l_{23}\cdot\vec{Q}^{CD},\\
\di B_{3t}+(\bar \l_3B_3)_x=\bar
l_{3}\cdot\vec{Q}^{CD}+\sum_{i=1,3}(\bar l_{3t}+\bar \l_i\bar
l_{3x})\cdot
 \bar r_iB_{i} +\bar l_{3t}\cdot\sum_{j=1}^{3}\bar{r}_{2j}B_{2j}.
\end{array}\right.
\end{equation}

Now we impose the following boundary condition to the linear
hyperbolic system \eqref{A2.7} on the domain $(t,x)\in
[h,T]\times\mathbf{R}$:
\begin{equation}\label{A2.8}
(B_1,B_{21},B_{22},B_{23},B_3)(t=T,x)=0.
\end{equation}
We can solve the linear diagonalized hyperbolic system
\eqref{A2.7} under the condition \eqref{A2.8} to have the following lemma.

\begin{lemma}\label{LemmaII} There exists a positive constant
$\delta_0$ such that if the wave strength $\d\leq \d_0,$ then there
exists a positive constant $C_{h,T}$ which is independent of $\v$,
such that
\begin{eqnarray}\label{A2.9}
&&\|\f{\partial^k}{\partial
x^k}(b_1,b_{21},b_{22},b_{23},b_3)(t,\cdot)\|^2_{L^2(dx)}
+\int_{h}^{T}\|\sqrt{|U_{1x}^{S_3}|}\f{\partial^k}{\partial
x^k}(b_1,b_{21},b_{22},b_{23},b_3)(t,\cdot)\|^2_{L^2(dx)}dt\nonumber\\
&& \leq C_{h,T}~ \v^{\frac52-2k},~ k=0,1,2,3.
\end{eqnarray}
\end{lemma}

\noindent{\bf Proof:}  From the wave curves defined in
\eqref{(2.1)}-\eqref{3-shock-curve}, we know $\d^{CD}+\d^{S_3}\leq
C\d_0$. Let $N=\frac1{\delta_0}$ and $\delta\leq \delta_0 \ll 1$,
then we have
\begin{equation}\nonumber
0<c_0\leq \big(\f{\hat{V}}{v^*}\big)^{\pm N},
\big(\f{V^{S_3}}{v^*}\big)^{\pm N}\leq C_0<\infty,
\end{equation}
where $\hat{V}$ is defined in \eqref{NS-CD} and  $c_0,C_0$ are
independent of $\delta_0$.

Without loss of generality, we assume that $v^*=1$ and
$\hat{V}_x>0$. If $v^*\neq1$, then we just replace $\hat{V}^{\pm N}$
and $(V^{S_3})^{\pm N}$ by $\di (\f{\hat{V}}{v^*})^{\pm N}$ and
$\di(\f{V^{S_3}}{v^*})^{\pm N}$, respectively.

Firstly,  multiplying the equation $(\ref{A2.7})_1$ by
$[(\hat{V})^{-N}+(V^{S_3})^{-N}]B_1$, we obtain
\begin{equation}\label{A2.10}
\begin{array}{ll}
\di \bigg(\frac12[\hat{V}^{-N}+(V^{S_3})^{-N}]B^2_1\bigg)_t
+\frac12N\bigg(\hat{V}^{-N-1}\hat{V}_t+(V^{S_3})^{-N-1}V_t^{S_3}\bigg)B^2_1\\
\di \qquad\qquad\quad+(\bar \l_{1}(\bar V,\bar P)B_1)_x[\hat{V}^{-N}+(V^{S_3})^{-N}]B_1\\
\di =[\hat{V}^{-N}+(V^{S_3})^{-N}]B_1\bigg(\bar
l_{1}\cdot\vec{Q}^{CD}+\sum_{i=1,3}(\bar l_{1t}+\bar \l_i\bar
l_{1x})\cdot \bar r_iB_i +\bar
l_{1t}\cdot\sum_{j=1}^{3}\bar{r}_{2j}B_{2j}\bigg).
\end{array}
\end{equation}
Now the second and third terms on the left hand side of
\eqref{A2.10} can be estimated by
\begin{eqnarray}\label{A2.11}
&&\frac12N\bigg(\hat{V}^{-N-1}\hat{V}_t+(V^{S_3})^{-N-1}V_t^{S_3}\bigg)B^2_1
+(\bar\l_1(\bar V,\bar P)B_1)_x[\hat{V}^{-N}+(V^{S_3})^{-N}]B_1\nonumber\\
&&\leq C_{h,T}B^2_1+C(\hat{V}_x+|U_{1x}^{S_3}|)B^2_1-\frac12N|\bar\l_1(\bar V,\bar P)|\hat{V}^{-N-1}\hat{V}_xB_1^2\nonumber\\
&&\qquad\qquad-\frac12N(1+\frac{|\bar \lambda_{1}(\bar V,\bar
P)|}{s_3})(V^{S_3})^{-N-1}|U_{1x}^{S_3}|B_1^2+(\cdots)_x.
\end{eqnarray}
Then we estimate the right hand side of \eqref{A2.10} term by term.
On one hand,
\begin{equation}\label{A2.12}
|[\hat{V}^{-N}+(V^{S_3})^{-N}]B_1\bar l_{1}\cdot\vec{Q}^{CD}|\leq
CB^2_1+ C|\vec{Q}^{CD}|^2.
\end{equation}
On the other hand, we have
\begin{equation}\label{A2.13}
|[\hat{V}^{-N}+(V^{S_3})^{-N}](\bar l_{1t}+\bar \l_1\bar
l_{1x})\cdot
 \bar r_1B^2_1|
\leq C_{h,T}B^2_1 +C\hat{V}_xB^2_1+C|U_{1x}^{S_3}|B^2_1,
\end{equation}
\begin{eqnarray}\label{A2.14}
&&|[\hat{V}^{-N}+(V^{S_3})^{-N}]B_1\bar
l_{1t}\cdot\sum_{j=1}^{3}\bar{r}_{2j}B_{2j}|\nonumber\\
&&\leq
C_{h,T}|\vec{B}|^2+C|U_{1x}^{S_3}|(B^2_1+|B_{21}|^2+|B_{22}|^2+|B_{23}|^2),
\end{eqnarray}
and
\begin{eqnarray}\label{A2.15}
&&|[\hat{V}^{-N}+(V^{S_3})^{-N}](\bar l_{1t}+\bar \l_3\bar
l_{1x})\cdot
\bar r_3B_1B_3|\nonumber\\
&&\leq
C_{h,T}(B_1^2+B_3^2)+C\hat{V}_x(B_1^2+B_3^2)+C|\l^{S_3}_3-s_3||U_{1x}^{S_3}|(B^2_1+B^2_3)\nonumber\\
&&\qquad+C_{h,T}\big[e^{-\frac{c_h}{\sigma}}e^{-\frac{c_h|x|}{\sigma}}+|(d_{1x},d_{2x},d_{3x})|\big](B_1^2+B_3^2)\nonumber\\
&&\leq
C_{h,T}(B_1^2+B_3^2)+C\hat{V}_x(B_1^2+B_3^2)+C\delta|U_{1x}^{S_3}|(B^2_1+B^2_3).
\end{eqnarray}
Substituting \eqref{A2.11}--\eqref{A2.15} into \eqref{A2.10} and
choosing $N$ large enough give
\begin{equation}
\label{A2.16}
\begin{array}{ll}
\di \Big(\frac12[\hat{V}^{-N}+(V^{S_3})^{-N}]B^2_1\Big)_t
-\frac14N|\bar\l_1(\bar V,\bar P)|\hat{V}^{-N-1}|\hat{V}_x|B_1^2-\frac14N(1+\frac{|\bar\l_1(\bar V,\bar P)|}{s_3})(V^{S_3})^{-N-1}|U_{1x}^{S_3}|B_1^2\\
\di \geq
-C_{h,T}|\vec{B}|^2-C_{h,T}|\vec{Q}^{CD}|^2-C|U_{1x}^{S_3}|\sum_{j=1}^{3}|B_{2j}|^2-C|\hat{V}_x|B_3^2-C\delta|U_{1x}^{S_3}|B_3^2.
\end{array}
\end{equation}
By multiplying the equation $(\ref{A2.7})_{j+1}$ by
$(V^{S_3})^{-N}B_{2j}~ (j=1,2,3)$,  and taking the summation of the
 equations together,  we obtain
\begin{equation}\label{A2.17}
\begin{array}{ll}
\di\bigg(\frac12(V^{S_3})^{-N}\sum_{j=1}^{3}|B_{2j}|^2\bigg)_t
-\frac12N(V^{S_3})^{-N-1}|U_{1x}^{S_3}|\sum_{j=1}^{3}|B_{2j}|^2\\
\di =(V^{S_3})^{-N}\Big(\sum_{j=1}^{3}B_{2j}\bar
l_{2j}\cdot\vec{Q}^{CD}+\sum_{i=1,3}(\bar l_{21t}+\bar \l_i\bar
l_{21x})\cdot \bar  r_iB_iB_{21}+\bar
l_{21t}\sum_{j=1}^{3}\bar{r}_{2j}B_{2j}B_{21}\Big).
\end{array}
\end{equation}
It is easy to check that
\begin{equation}\label{A2.18}
|(V^{S_3})^{-N}\sum_{j=1}^{3}B_{2j}\bar l_{2j}\cdot\vec{Q}^{CD}|\leq
C\sum_{j=1}^{3}|B_{2j}|^2+ C|\vec{Q}^{CD}|^2.
\end{equation}
From the construction of viscous contact wave and \eqref{A2.3}, it
holds that $l_{21x}^{CD}=O(1)\delta^{CD}$. Then the terms on the
right hand side of \eqref{A2.17} can be estimated by
\begin{eqnarray}\label{A2.19}
&&|(V^{S_3})^{-N}(\bar l_{21t}+\bar \l_1 \bar l_{21x})\cdot
 \bar r_1B_1B_{21}|\nonumber\\
&&\leq C_{h,T}|\vec{B}|^2+C|U_{1x}^{S_3}|(B^2_1+|B_{21}|^2)
+C_{h,T}\big[e^{-\frac{c_h}{\sigma}}e^{-\frac{c_h|x|}{\sigma}}+|(d_{1x},d_{2x},d_{3x})|\big](B_1^2+|B_{21}|^2)\nonumber\\
&&\leq C_{h,T}|\vec{B}|^2+C|U_{1x}^{S_3}|(B^2_1+|B_{21}|^2).
\end{eqnarray}
Similar to \eqref{A2.15} and \eqref{A2.19}, we have
\begin{eqnarray}\label{A2.20}
&&|(V^{S_3})^{-N}(\bar l_{21t}+\bar \l_3 \bar l_{21x})\cdot \bar
r_3B_3B_{21}|\leq
C_{h,T}|\vec{B}|^2+C\delta|U_{1x}^{S_3}|(B^2_3+|B_{21}|^2),
\end{eqnarray}
and
\begin{eqnarray}\label{A2.21}
&&|(V^{S_3})^{-N}\bar
l_{21t}\cdot\sum_{j=1}^{3}\bar{r}_{2j}B_{2j}B_{21}|\leq
C_{h,T}|\vec{B}|^2+C|U_{1x}^{S_3}|\sum_{j=1}^{3}|B_{2j}|^2.
\end{eqnarray}

Substituting \eqref{A2.18}-\eqref{A2.21} into \eqref{A2.17} and
choosing $N$ large enough give
\begin{eqnarray}\label{A2.22}
&&\bigg(\frac12(V^{S_3})^{-N}\sum_{j=1}^{3}|B_{2j}|^2\bigg)_t
-\frac14N(V^{S_3})^{-N-1}|U_{1x}^{S_3}|\sum_{j=1}^{3}|B_{2j}|^2\nonumber\\
&&\geq-C_{h,T}|\vec{B}|^2-C|\vec{Q}^{CD}|^2-C\delta|U_{1x}^{S_3}|B^2_3-C|U_{1x}^{S_3}|(B^2_1+B_{21}^2).
\end{eqnarray}
Multiplying $(\ref{A2.7})_5$ by $\hat{V}^{N}B_3$ yields
\begin{equation}\label{A2.23}
\begin{array}{ll}
\di \bigg(\frac12\hat{V}^{N}B^2_3\bigg)_t+\frac12\bar
\lambda_{3x}(\bar V,\bar P)\hat{V}^{N}B_3^2
-\frac{\bar\lambda_3}{2}N \hat{V}^{N-1}|\hat{V}_x|B^2_3\\
\di =\f12N\hat{V}^{N-1}\hat{V}_t B_3^2 +\hat{V}^{N}B_3\bigg(\bar
l_{3}\cdot\vec{Q}^{CD}+\sum_{i=1,3}(\bar l_{3t}+\bar \l_i\bar
l_{3x})\cdot \bar r_iB_i +\bar
l_{3t}\cdot\sum_{j=1}^{3}\bar{r}_{2j}B_{2j}\bigg).
\end{array}
\end{equation}
The wave interaction estimations imply
\begin{equation}\label{A2.24}
\begin{array}{ll}
\di \bar \l_{3x}&\di
=\l^{R_1}_{3x}+\l^{CD}_{3x}+\l^{S_3}_{3x}+\big(\bar
\l_{3x}-\l^{R_1}_{3x}-\l^{CD}_{3x}-\l^{S_3}_{3x}\big)
\leq\l^{S_3}_{3x}+C|\hat{V}_x|+C_{h,T},
\end{array}
\end{equation}
so that
\begin{equation}\label{A2.25}
\frac12\bar \lambda_{3x}\hat{V}^{N}B_3^2\leq
\frac12\lambda_{3x}^{S_3}\hat{V}^{N}B_3^2+C|\hat{V}_x|\hat{V}^{N}B_3^2+C_{h,T}B_3^2.
\end{equation}
The other terms on the right hand side of \eqref{A2.23} can be
estimated similarly to \eqref{A2.15}-\eqref{A2.19}. Then we have
\begin{equation}\label{A2.26}
|\f12N\hat{V}^{N-1}\hat{V}_t B_3^2+\hat{V}^{N}B_3\bar
l_{3}\cdot\vec{Q}^{CD}|\leq CB^2_3+ C|\vec{Q}^{CD}|^2,
\end{equation}
\begin{eqnarray}\label{A2.27}
&&|\hat{V}^{N}(\bar l_{3t}+\bar \l_1 \bar l_{3x})\cdot
 \bar r_1B_1B_3|\nonumber\\
&&\leq C_{h,T}(B_1^2+B_3^2)+\beta|U_{1x}^{S_3}|\hat{V}^{N}B^2_3
+C_{\beta}|U_{1x}^{S_3}|B^2_1+C|\hat{V}_x|\hat{V}^{N}(B_1^2+B_3^2),
\end{eqnarray}
\begin{eqnarray}\label{A2.29}
&&|\hat{V}^{N}(\bar l_{3t}+\bar \l_3 \bar l_{3x})\cdot
 \bar r_3B_3^2|\leq
C_{h,T}\hat{V}^{N}B_3^2+C\delta|U_{1x}^{S_3}|\hat{V}^{N}B^2_3+C|\hat{V}_x|\hat{V}^{N}B^3_2,
\end{eqnarray}
and
\begin{equation}\label{A2.28}
|\hat{V}^{N}B_3\bar l_{3t}\cdot
\sum_{j=1}^{3}\bar{r}_{2j}B_{2j}|\leq
C_{h,T}|\vec{B}|^2+\beta|U_{1x}^{S_3}|\hat{V}^{N}B^2_3+C_{\beta}|U_{1x}^{S_3}|\sum_{j=1}^{3}|B_{2j}|^2.
\end{equation}
By choosing $\beta$ and $\delta_0$ small enough, substituting
\eqref{A2.25}--\eqref{A2.29} into \eqref{A2.23} gives
\begin{eqnarray}\label{A2.30}
&&\bigg(\frac12\hat{V}^{N}B^2_3\bigg)_t-\frac{|\lambda_{3x}^{S_3}|}{4}\hat{V}^{N}B_3^2
-\frac{\bar\lambda_3}{4}N \hat{V}^{N-1}|\hat{V}_x|B^2_3\nonumber\\
&&\qquad\geq
-C_{h,T}|\vec{B}|^2-C|\vec{Q}^{CD}|^2-C|\hat{V}_x|\hat{V}^{N}B^2_1
-C_{\beta}|U_{1x}^{S_3}|(B^2_1+\sum_{j=1}^{3}|B_{2j}|^2).
\end{eqnarray}
By combining \eqref{A2.16},\eqref{A2.22} and \eqref{A2.30} and
noticing that $\delta_0\ll 1$ and $N\gg1$, we can get
\begin{equation}\label{A2.31}
\begin{array}{ll}
\di
\bigg(\frac12[\hat{V}^{-N}+(V^{S_3})^{-N}]B^2_1+\frac12(V^{S_3})^{-N}\sum_{j=1}^{3}|B_{2j}|^2
+\frac12\hat{V}^{N}B^2_3\bigg)_t\\
\di -\frac14N|\bar\l_1|\hat{V}^{-N-1}|\hat{V}_x|B_1^2
-\frac14N(1+\frac{|\bar\l_1|}{s_3})\hat{V}^{-N-1}|U_{1x}^{S_3}|B_1^2\\
\di -\frac14N(V^{S_3})^{-N-1}|U_{1x}^{S_3}|\sum_{j=1}^{3}|B_{2j}|^2
-\frac{|\lambda_{3x}^{S_3}|}{4}\hat{V}^{N}B_3^2
-\frac{\bar\lambda_3}{4}N \hat{V}^{N-1}|\hat{V}_x|B_3^2\\
\di \geq -C_{h,T}|\vec{B}|^2-C_{h,T}|Q^{CD}|^2.
\end{array}
\end{equation}
Integrating \eqref{A2.31} over $[t,T]\times \mathbf{R}$ with $t\in
(h,T)$ then yields
\begin{eqnarray}\label{A2.32}
&&\int_{\mathbf{R}}|\vec{B}|^2dx+\int_{t}^{T}\int_{\mathbf{R}}\Big[|\hat{V}_x|(B_1^2+B_3^2)
+|U_x^{S_3}||\vec{B}|^2\Big]dxdt\nonumber\\
&&\leq C_{h,T}\int_{\mathbf{R}}|\vec{B}|^2+C_{h,T}~\v^{\frac52}, \ \
\forall t\in[h,T].
\end{eqnarray}
Applying Gronwall inequality to \eqref{A2.32} gives
\begin{equation}\label{A2.33}
\int_{\mathbf{R}}|\vec{B}|^2dx
+\int_{t}^{T}\int_{\mathbf{R}}\Big[|\hat{V}_x|(B_1^2+B_3^2)
+|U_{1x}^{S_3}||\vec{B}|^2\Big]dxdt\leq C_{h,T}~\v^{\frac52},\
\forall t\in[h,T].
\end{equation}
This completes the proof for the case  when $k=0$ in  Lemma
\ref{LemmaII}. The case $k=1,2,3$ can be proved similarly to the
differentiated system, and we omit the details for brevity.
$\hfill\Box$

\subsection{Superposition of Waves}
With the above preparation, finally,  the approximate superposition wave
$(V,U,\mathcal{E})(t,x)$ can be defined by
\begin{equation}
\begin{array}{l}
 \left(\begin{array}{cc} V\\ U_i \\
 \mathcal{E}
\end{array}
\right)(t, x)= \left(\begin{array}{cc}  \bar V+b_1\\ \bar U_i+b_{2i} \\
\bar{\mathcal{E}}+b_3
\end{array}
\right)(t, x),~~i=1,2,3,
 \end{array}
\label{wave}
\end{equation}
where $\mathcal{E}=\T+\f{|U|^2}{2}.$

Thus, we have
\begin{equation}
\T=\bar\T-\sum_{i=1}^3\bar U_i b_{2i}+b_3-\f{|b_2|^2}{2},
\label{Theta}
\end{equation}
where $b_2=(b_{21},b_{22},b_{23})^t$ and
$|b_2|^2=\sum_{i=1}^3b_{2i}^2.$

From the construction of the contact wave and Lemma \ref{Lemma
2.3} and by noting that $\s=\v^{\f15}$, we have the following relation
between the approximate wave pattern $(V,U,\mathcal{E},\T)(t,x)$ of
the Boltzmann equation and the inviscid superposition wave pattern
$(\wt V,\wt U,\wt{\mathcal{E}},\wt\T)(t,x)$ to the Euler equations
\begin{equation}
\begin{array}{ll}
\di |(V,U,\mathcal{E},\T)(t,x)-(\wt V,\wt U,\wt{\mathcal{E}},\wt\T)(t,x)|\\
\di\leq C\Big[|(V^{R_1},U^{R_1},\mathcal{E}^{R_1},\T^{R_1})(t,x)-(v^{r_1},u^{r_1},E^{r_1},\t^{r_1})(t,x)|+|(d_1,d_2,d_3)(t,x)|\\
\di\quad +|(V^{CD},U^{CD},\mathcal{E}^{CD},\T^{CD})(t,x)-(v^{cd},u^{cd},E^{cd},\t^{cd})(t,x)|\\
\di\quad +|(V^{S_3},U^{S_3},\mathcal{E}^{S_3},\T^{S_3})(t,x)-(v^{s_3},u^{s_3},E^{s_3},\t^{s_3})(t,x)|+|(b_1,b_{21},b_{22},b_{23},b_3)(t,x)|\Big] \\
 \di\leq
C_{h,T}\Big[\f{1}{t}\big(\s\ln(1+t)+\s|\ln\s|\big)+\f\v\s+\d^{CD}e^{-\f{cx^2}{\v(1+t)}}+\d^{S_3}e^{-c\f{\d^{S_3}|x-s_3t|}{\v}}+\v^{\f34}\Big]\\
\di \leq
C_{h,T}\Big[\v^{\f15}|\ln\v|+\d^{CD}e^{-\f{cx^2}{\v(1+t)}}+\d^{S_3}e^{-c\f{\d^{S_3}|x-s_3t|}{\v}}\Big].
\end{array}\label{difference}
\end{equation}

Moreover, the approximate wave pattern $(V,U,\mathcal{E},\T)(t,x)$ satisfies
\begin{equation}
\left\{
\begin{array}{l}
  \di V_t-U _{1x} = 0,  \\
   \di U_{1t}+P_x
    = \f43\v (\frac{\mu(\T)U_{1x}}{V}) _{x}-\int\x_1^2\Pi^{CD}_{11x}d\x-\int\x_1^2\Pi^{S_3}_{1x}d\x+\bar Q_{1x}+Q_{1x},\\
\di
U_{it}=\v(\f{\mu(\T)U_{ix}}{V})_x-\int\x_1\x_i\Pi^{CD}_{11x}d\x-\int\x_1\x_i\Pi^{S_3}_{1x}d\x+\bar
Q_{ix}+Q_{ix},~~i=2,3,
    \\
  \di \mathcal{E} _t+(P U_{1})_x
    =\v( \frac{\k(\T)\T_{x}}{ V})_x + \f43\v(\frac{ \mu(\T)U_1 U_{1x}}{
    V})_x+\sum_{i=2}^3\v(\frac{ \mu(\T)U_{i} U_{ix}}{V})_x\\
    \di\qquad\qquad -\int\x_1\f{|\x|^2}{2}\Pi^{CD}_{11x}d\x-\int\x_1\f{|\x|^2}{2}\Pi^{S_3}_{1x}d\x
     +\bar Q_{4x}+Q_{4x},
\end{array}
\right.\label{(2.41)}
\end{equation}
where $ P =p( V , \T  )$ and
\begin{equation}\begin{array}{ll}
   \di Q_1&\di=\Big[P-\bar P-(\bar P_v b_1+\bar P_{u}\cdot b_2+\bar P_E b_3)\Big]-\f43\v\Big[\f{\mu(\T)U_{1x}}{V}-\f{\mu(\bar\T)\bar U_{1x}}{\bar V}\Big],\\
   &:=\di Q_{11}+Q_{12},\\
 \di Q_i&\di=-\v\Big[\f{\mu(\T)U_{ix}}{V}-\f{\mu(\bar\T)\bar U_{ix}}{\bar V}\Big],~~i=2,3,\\
    \di Q_4&\di= \Big[PU_{1}-\bar P\bar U_{1}-\left((\bar P\bar U_1)_v b_1+(\bar P\bar U_1)_{u}\cdot b_2+(\bar
P\bar U_1)_E
b_3\right)\Big]\\[3mm]
  &\di\quad-\v\Big[(\f{\k(\T)\T_x}{V}-\f{\k(\bar\T)\bar\T_x}{\bar
V})+\f43(\f{\mu(\T)U_1U_{1x}}{V}-\f{\mu(\bar \T)\bar U_1\bar
U_{1x}}{\bar
  V})\\
  &\di\quad\qquad +\sum_{i=2}^3(\f{\mu(\T)U_iU_{ix}}{V}-\f{\mu(\bar \T)\bar U_i\bar
U_{ix}}{\bar
  V})\Big]\\
  &\di:=Q_{41}+Q_{42}.
\end{array}\label{Q1-2}
\end{equation}
Straightforward calculation shows that
\begin{equation}
(Q_{11},Q_{41})=O(1)|\vec{b}|^2.\label{Q11-21}
\end{equation}

\section{Proof of Main Result}
\setcounter{equation}{0}

With the above preparation, we will give the proof of the main theorem as follows.
For this, we will first reformulate the problem in the following subsection. The energy
estimates will then be given in the second subsection.

\subsection{Reformulation of the Problem}

We now reformulate the system by introducing  a scaling for the
independent variables. Set
\begin{equation}
y=\f x\v,~~\tau=\f t\v. \label{scaling}
\end{equation}
In the following, we will also use the notations $(v,u,\t)(\tau,y), \mb G(\tau,y,\x),
\Pi_1(\tau,y,\x)$ and $(V,U,\T)(\tau,y)$, etc., in the scaled
independent variables. Set the perturbation around the superposition
wave $(V,U,\T)(\tau,y)$ by
\begin{equation}
\begin{array}{ll}
\di(\phi,\psi,\omega,\zeta)(\tau,y)=(v-V,u-U,E-\mathcal{E},\t-\T)(\tau,y),\\[2mm]
\di \wt
{\mathbf{G}}(\tau,y,\x)=\mathbf{G}(\tau,y,\x)-\mathbf{G}^{S_3}(\tau,y,\x),\\[2mm]
\di \wt f(\tau,y,\x)=f(\tau,y,\x)-F^{S_3}(\tau,y,\x).
\end{array}\label{perturb}
\end{equation}
Under this scaling, the hydrodynamic limit problem is reduced to a
time asymptotic stability problem for the
Boltzmann equation.

 In particular, we can choose the initial value  as
\begin{equation}
(\phi,\psi,\o)(\tau=\f h\v,y)=(0,0,0),\quad \wt{\mb{G}}(\tau=\f h\v,
y,\x)=0.\label{initial-value}
\end{equation}
Introduce the anti-derivative variables
$$
(\Phi,\Psi,\bar
W)(\tau,y)=\int_{-\i}^y(\phi,\psi,\o)(\tau,y^\prime)d y^\prime.
$$



Then $(\Phi,\Psi,\bar W)(\tau,y)$ satisfies that
\begin{equation}
\left\{
\begin{array}{l}
\di \Phi_\tau-\Psi_{1y}=0,\\
\di \Psi_{1\tau}+(p-P)=\f43\Big(\frac {\mu(\t)u_{1y}}v-\f{\mu(\T)U_{1y}}{V}\Big)-\int\xi_1^2(\Pi_{1}-\Pi^{CD}_{11}-\Pi^{S_3}_1)d\xi-\bar Q_1-Q_{1},\\
\di \Psi_{i\tau}=\Big(\frac {\mu(\t)u_{iy}}v-\f{\mu(\T)U_{iy}}{V}\Big)-\int\xi_1\xi_i(\Pi_{1}-\Pi^{CD}_{11}-\Pi^{S_3}_1)d\xi-\bar Q_i-Q_i,~~i=2,3,\\
\di \bar
W_\tau+(pu_1-PU_1)=\Big(\frac{\k(\t)\theta_{y}}{v}-\f{\k(\T)\T_y}{V}\Big)+\f43\Big(\frac
{\mu(\t)u_1u_{1y}}{v}-\f{\mu(\T)U_1U_{1y}}{V}\Big)\\
\di\qquad+\sum_{i=2}^3\Big(\frac
{\mu(\t)u_iu_{iy}}v-\f{\mu(\T)U_iU_{iy}}{V}\Big)-\int\xi_1\f{|\x|^2}{2}(\Pi_{1}-\Pi^{CD}_{11}-\Pi^{S_3}_1)d\xi-\bar
Q_4-Q_{4}.
\end{array} \right.\label{new-sys}
\end{equation}
To precisely capture the dissiaption of heat conduction, we introduce another variable related to the
absolute temperature
$$
W=\bar W-U\cdot\Psi=\bar W-\sum_{i=1}^3U_i\Psi_i,
$$
then
\begin{equation}
 \z=W_y-(\f{|\Psi_y|^2}{2}-U_{y}\cdot\Psi).\label{zeta}
\end{equation}
Linearizing the system \eqref{new-sys} around the approximate wave pattern
$(V,U,\T)(\tau,y)$ implies that
\begin{equation}
\left\{
\begin{array}{ll}
\di \Phi_\tau-\Psi_{1y}=0,\\
\di
\Psi_{1\tau}-\f{Z}{V}\Phi_y+\f{2}{3V}W_y+\f{2}{3V}U_{y}\cdot\Psi-\f43\f{\mu^\prime(\T)}{V}(W_y+U_{y}\cdot\Psi)U_{1y}=\f43\f{\mu(\T)}{V}\Psi_{1yy}\\
 \di\hspace{5cm} -\int\xi_1^2(\Pi_{1}-\Pi^{CD}_{11}-\Pi^{S_3}_1)d\xi+N_1-\bar
Q_1-Q_{1},\\
\di \Psi_{i\tau}+\f{\mu(\T)U_{iy}}{V^2}\Phi_y-\f{\mu^\prime(\T)}{V}(W_y+U_{y}\cdot\Psi)U_{iy}=\f{\mu(\T)}{V}\Psi_{iyy}\\
\di\hspace{3cm} -\int\xi_1\xi_i(\Pi_{1}-\Pi^{CD}_{11}-\Pi^{S_3}_1)d\xi+N_i-\bar Q_i-Q_i, i=2,3,\\
\di
W_\tau+Z\Psi_{1y}-\sum_{i=2}^3\f{\mu(\T)U_{iy}}{V}\Psi_{iy}+U_{\tau}\cdot\Psi-\f{\k(\T)}{V}(U_{y}\cdot\Psi)_y+\f{\k(\T)}{V^2}\T_y\Phi_y\\
[3mm] \di~-\f{\k^\prime(\T)}{V}(W_y+U_{y}\cdot\Psi)\T_{y}=
\f{\k(\T)}{V}W_{yy}-\int\xi_1\f{|\x|^2}{2}(\Pi_{1}-\Pi^{CD}_{11}-\Pi^{S_3}_1)d\xi\\
\di~
+\sum_{i=1}^3U_i\int\xi_1\x_i(\Pi_{1}-\Pi^{CD}_{11}-\Pi^{S_3}_1)d\xi
+N_4-\bar Q_4+\sum_{i=1}^3U_i\bar Q_i-Q_{4}+\sum_{i=1}^3U_iQ_{i},
\end{array} \right.\label{sys}
\end{equation}
where
\begin{equation}
Z=P-\f43\f{\mu(\T)U_{1y}}{V},\label{Z}
\end{equation}
\begin{equation}
\begin{array}{ll}
N_1&\di=\f{p-P}{V}\Phi_y+\f{1}{3V}|\Psi_y|^2+\f43(\f{\mu(\t)}{v}-\f{\mu(\T)}{V})\Psi_{1yy}\\[3mm]
&\di~~
+\f43U_{1y}\Big[\f{\mu(\t)}{v}-\f{\mu(\T)}{V}+\f{\mu(\T)}{V^2}\Phi_y-\f{\mu^\prime(\T)}{V}(\zeta+\f{|\Psi_y|^2}{2})\Big]\\
&\di=O(1)\Big[|\Phi_y|^2+|\Psi_y|^2+|\zeta|^2+|\Psi_{1yy}|^2\Big],
\end{array}\label{N1}
\end{equation}
\begin{equation}
\begin{array}{ll}
N_i&\di=(\f{\mu(\t)}{v}-\f{\mu(\T)}{V})\Psi_{iyy}+U_{iy}\Big[\f{\mu(\t)}{v}-\f{\mu(\T)}{V}+\f{\mu(\T)}{V^2}\Phi_y-\f{\mu^\prime(\T)}{V}(\zeta+\f{|\Psi_y|^2}{2})\Big]\\
&\di=O(1)\Big[|\Phi_y|^2+|\Psi_y|^2+|\zeta|^2+|\Psi_{iyy}|^2\Big],~~i=2,3,
\end{array}\label{Ni}
\end{equation}
and
\begin{equation}
\begin{array}{ll}
N_4&\di=-(p-P)\Psi_{1y}-\f{\k(\T)}{V}\Psi_y\cdot\Psi_{yy}+(\f{\k(\t)}{v}-\f{\k(\T)}{V})\zeta_y\\
&\di~~+\f43(\f{\mu(\t)u_{1y}}{v}-\f{\mu(\T)U_{1y}}{V})\Psi_{1y}+\sum_{i=2}^3\Big(\f{\mu(\t)u_{iy}}{v}-\f{\mu(\T)U_{iy}}{V}\Big)\Psi_{iy}\\[3mm]
&\di~~
+\T_{y}\Big[\f{\k(\t)}{v}-\f{\k(\T)}{V}+\f{\k(\T)}{V^2}\Phi_y-\f{\k^\prime(\T)}{V}(\zeta+\f{|\Psi_y|^2}{2})\Big]\\
&\di=O(1)\Big[|\Phi_y|^2+|\Psi_y|^2+|\zeta|^2+|\Psi_{yy}|^2+|\zeta_{y}|^2\Big].
\end{array}\label{N4}
\end{equation}
We now derive the equation for the  non-fluid component $\wt{\mb{G}}
(\tau,y,\xi)$ in the scaled independent variables. From
(\ref{(1.18)}), we have
\begin{equation}
\begin{array}{ll}
\di
\wt{\mb{G}}_{\tau}-\mb{L}_\mb{M}\wt{\mb{G}}&\di=\f{u_1}{v}\wt{\mb{G}}
_y-\f{1}{v}\mb{P}_1(\xi_1\wt{\mb{G}}_y)-\Big[\f{1}{v}\mb{P}_1(\xi_1\mb{M}_y)
-\f{1}{V^{S_3}}\mb{P}^{S_3}_1(\xi_1\mb{M}^{S_3}_y)\Big]\\[4mm]
&\di\quad+2Q(\wt{\mb{G}},\mb{G}^{S_3})+Q(\wt{\mb{G}},\wt{\mb{G}})+J_1,
\end{array}
\label{t-G}
\end{equation}
where
\begin{equation}
J_1=
\big(\mb{L}_{\mb{M}}-\mb{L}_{\mb{M}^{S_3}}\big)\mb{G}^{S_3}+\big(\f
uv-\f{U^{S_3}_1}{V^{S_3}}\big)\mb{G}^{S_3}_y
-\Big[\f{1}{v}\mb{P}_1(\xi_1\mb{G}^{S_3}_y)-\f{1}{V^{S_3}}\mb{P}^{S_3}_1(\xi_1\mb{G}^{S_3}_y)\Big].
\label{J1}
\end{equation}
Let
\begin{equation}
\mb{G}^{R_1}(\tau,y,\xi)=\f{3}{2v\t}\mb{L}^{-1}_\mb{M}\{\mb{P}_1[\xi_1(\f{|\xi-u|^2}{2\t}{\T}^{R_1}_y+\xi\cdot{U}^{R_1}_{y})\mb{M}]\},
\label{G-R1}
\end{equation}
and
\begin{equation}
\wt{\mb{G}}
_1(\tau,y,\xi)=\wt{\mb{G}}(\tau,y,\xi)-\mb{G}^{R_1}(\tau,y,\xi)-\mb{G}^{CD}(\tau,y,\xi),
\label{G1}
\end{equation}
where $\mb{G}^{CD}(\tau,y,\xi)$ is defined in \eqref{G-CD}.

Then $\mb{G} _1(\tau,y,\xi)$ satisfies
\begin{equation}
\wt{\mb{G}} _{1\tau}-\mb{L}_\mb{M}\wt{\mb{G}}_1=
\f{u_1}{v}\wt{\mb{G}}
_y-\f{1}{v}\mb{P}_1(\xi_1\wt{\mb{G}}_y)+2Q(\wt{\mb{G}},\mb{G}^{S_3})+Q(\wt{\mb{G}},\wt{\mb{G}})+J_1+J_2-\mb{G}^{R_1}_{\tau}-\mb{G}^{CD}_{\tau}.
\label{G1e}
\end{equation}
with
\begin{equation}
\begin{array}{ll}
J_2&\di =-\Big[\f{1}{v}\mb{P}_1(\xi_1\mb{M}_y)
-\f{1}{V^{S_3}}\mb{P}^{S_3}_1(\xi_1\mb{M}^{S_3}_y)\\
&\di \qquad\qquad~~~-\f{3}{2v\t}\mb{P}_1\Big(\xi_1
(\f{|\xi-u|^2}{2\t}({\T}^{R_1}_y+\T^{CD}_y)+\xi\cdot ({U}^{R_1}_{y}+{U}^{CD}_{y}))\mb{M}\Big)\Big]\\
&\di =-\Big[\f{3}{2v\t}\mb{P}_1\Big(\xi_1
\big[\f{|\xi-u|^2}{2\t}(\t-{\T}^{R_1}-\T^{CD})_y+\xi\cdot
(u-U^{R_1}-U^{CD})_y\big]\mb{M}\Big)
\\
&\di\qquad\qquad~~~ -\f{3}{2V^{S_3}\T^{S_3}}\mb{P}^{S_3}_1\Big(\xi_1
\big[\f{|\xi-U^{S_3}|^2}{2\T^{S_3}}\T^{S_3}_y+\xi\cdot
U^{S_3}_{y}\big]\mb{M}^{S_3}\Big)\Big].
\end{array}
\label{J2}
\end{equation}

 Notice that in \eqref{G1} and \eqref{G1e},
$\mb{G}^{R_1}$ and $\mb{G}^{CD}$ are subtracted from $\wt{\mb{G}}$
when carrying out the lower order energy estimates because $\di
\int_{\f h\v}^{\f
T\v}\|(\T^{R_1}_y,U^{R_1}_{1y})\|_{L^2(dy)}^2d\tau$ is uniformly
bounded with respect to $\v$, while $\di\int_{\f h\v}^{\f
T\v}\|\T^{CD}_y\|_{L^2(dy)}^2d\tau$ is only of the order of $\v^{-\f12}$.
Both do not give any  decay  with respect to Knudsen number
$\v$ in the above integrals.

From (\ref{Lag-B}) and the scaling transformation
\eqref{scaling}, we have
\begin{equation}
\di f_{\tau}-\f{u_1}{v}f_y+\f{\xi_1}{v}f_y=Q(f,f). \label{Lag-B1}
\end{equation}
Thus, we have the equation for $\wt f$ defined in \eqref{perturb}
\begin{equation}
\wt f_\tau-\f{u_1}{v}\wt f_y+\f{\x_1}{v}\wt f_y=\mb{L}_\mb{M}\wt
{\mb{G}}+Q(\wt{\mb G},\wt{\mb G})+J_F, \label{tidle-f}
\end{equation}
with
\begin{equation}
J_F=(\f{u_1}{v}-\f{U^{S_3}_1}{V^{S_3}})F^{S_3}_y-(\f1v-\f{1}{V^{S_3}})\x_1F^{S_3}_y+2Q(\mb{M}-\mb{M}^{S_3},\mb{G}^{S_3})+2Q(\wt{\mb
G},\mb{G}^{S_3}).\label{JF}
\end{equation}
The estimation on the fluid and non-fluid components governed by the
above equations will be given in the next subsection.
In the following, we will state the main estimate we want to obtain
and also give the a priori estimate.

Note that to prove the main theorem in this paper, it is sufficient to
prove the following theorem on
 the Boltzmann equation
\eqref{Lag-B1} in the scaled independent variables based on
 the construction of the approximate wave pattern.

\begin{theorem}\label{Theorem 4.1} There exist a small positive
constants $\delta_1$ and a global Maxwellian
$\mb{M}_\star=\mb{M}_{[v_\star,u_\star,\theta_\star]}$ such that if
the wave strength $\delta$ satisfies  $\delta \le \delta_1$, then on
the time interval $[\f h\v, \f T\v]$ for any  $0<h<T$, there is a positive constant
$\v_1(\d,h,T)$. If the Knudsen number $\v\leq\v_1$, then the problem
(\ref{Lag-B1}) admits a family of smooth solution
$f^{\v,h}(\tau,y,\xi)$ satisfying
\begin{equation}
\begin{array}{l}
\di \sup_{\tau\in[\f h\v, \f
T\v]}\sup_{y\in\mathbf{R}}\|f^{\v,h}(\tau,y,\xi)-\mb{M}_{[V,U,\Theta]}(\tau,y,\xi)
\|_{L^2_{\xi}(\frac{1}{\sqrt{\mb{M}_\star}})}\le
C\v^{\f15}.\\
\end{array}
\label{(4.10)}
\end{equation}
\end{theorem}

Consider the reformulated system \eqref{sys} and \eqref{G1e}. Since
the local existence of solution to \eqref{sys} and \eqref{G1e} is
known, cf. \cite{Guo} and \cite{UYZ}, to prove the
existence on the time interval $[\f h\v, \f T\v]$, we only need to
close the following a priori estimate by the continuity argument. Set
\begin{equation}
\begin{array}{ll}
\di \mathcal{N}(\tau)=\sup_{\f h\v\leq \tau^\prime\leq \tau}\Bigg\{
\|(\Phi,\Psi,W)(\tau^\prime,\cdot)\|^2+\|(\phi,\psi,\zeta)(\tau^\prime,\cdot)\|_1^2+\int\int\f{|\wt{\mb{G}}_1|^2}{\mb{M}_\star}d\xi dy\\
\di\qquad +\sum_{|\a^\prime|=1} \int\int\f{|\partial^{\a^\prime}
\wt{\mb{G}} |^2}{\mb{M}_\star}d\xi
dy+\sum_{|\a|=2}\int\int\f{|\partial^\a \wt
f|^2}{\mb{M}_{\star}}d\xi dy\Bigg\}\leq \chi^2 =
\v^{\f{1}{10}},~~\forall \tau\in[\f h\v, \f T\v],
\end{array}
\label{priori}
\end{equation}
where  $\partial^\a,\partial^{\a^\prime}$ denote the derivatives
with respect to $y$ and $\tau$, and
$\mb{M}_\star$ is a global Maxwellian to be chosen.

 \begin{remark} In the paper, we simply choose the initial data for the Boltzmann equation
(\ref{Lag-B1}) as
\begin{equation}
 f^{\v,h}(\tau=\f h\v,y,\xi)=\mb{M}_{[V,U,\T]}(\f
h\v,y,\xi)+\mb{G}^{S_3}(\f h\v,y,\x), \label{(4.12)}
\end{equation}
so that
\begin{equation}
\begin{array}{ll}
\mathcal{N}(\tau)|_{\tau=\f h\v}&\di\leq  C_{h,T}~\v^{\f12}.
\end{array} \label{(4.13)}
\end{equation}
In this case, the functional measuring the perturbation
$\mathcal{N}(\tau)$ is smaller at the initial time
$\tau=\f h\v$ than the estimate
given in Theorem \ref{Theorem 4.1} in the whole time interval
 when  $\v$ is small.
\end{remark}

Note that the a priori assumption (\ref{priori}) implies that
\begin{equation}
\|(\Phi,\Psi,W)\|^2_{L_\i}+\|(\phi,\psi,\z,\o)\|^2_{L_\i}\leq
C\chi^2, \label{L-infty}
\end{equation}
and
\begin{equation}
\begin{array}{ll}
\di \|\int\f{\wt{\mb{G}}_1^2}{\mb{M}_\star}d\xi\|_{L_\i^y}&\di \leq
C\left(\int\int\f{\wt{\mb{G}} _1^2}{\mb{M}_\star}d\xi
dy\right)^{\f{1}{2}}\cdot\left(\int\int\f{|\wt{\mb{G}}
_{1y}|^2}{\mb{M}_\star}d\xi dy\right)^{\f{1}{2}}\\
&\di \leq C\left(\int\int\f{\wt{\mb{G}} _1^2}{\mb{M}_\star}d\xi
dy\right)^{\f{1}{2}}\cdot\left(\int\int\f{|\wt{\mb{G}}
_{y}|^2+|\mb{G}^{R_1}_{y}|^2+|\mb{G}^{CD}_{y}|^2}{\mb{M}_\star}d\xi
dy\right)^{\f{1}{2}}\\[2mm]
&\di \leq
C\chi\big[\chi+\|(v_y,u_y,\t_y)(\T^{R_1}_y,U^{R_1}_{1y},\T^{CD}_y,U^{CD}_{y})\|_{L^2(dy)}\\[2mm]
&\di\qquad\quad
+\|(\T^{R_1}_{yy},U^{R_1}_{1yy},\T^{CD}_{yy},U^{CD}_{yy})\|_{L^2(dy)}\big]\leq
C_{h,T}\chi^2.
\end{array}
\label{G1-infty}
\end{equation}
Furthermore,  for $|\a^\prime|=1$,
\begin{equation}
\begin{array}{ll}
 \di\|\int\f{|\partial^{\a^\prime}
\wt{\mb{G}}|^2}{\mb{M}_\star}d\xi\|_{L_\i^y}\leq
C\left(\int\int\f{|\partial^{\a^\prime} \wt{\mb{G}}
|^2}{\mb{M}_\star}d\xi
dy\right)^{\f{1}{2}}\cdot\left(\int\int\f{|\partial^{\a^\prime}
\wt{\mb{G}}
_y|^2}{\mb{M}_\star}d\xi dy\right)^{\f{1}{2}}\\
 \di \leq C\left(\int\int\f{|\partial^{\a^\prime} \wt{\mb{G}}
|^2}{\mb{M}_\star}d\xi
dy\right)^{\f{1}{2}}\cdot\left(\int\int\f{|\partial^{\a^\prime}\wt
f_y|^2
+|\partial^{\a^\prime}(\mb{M}-\mb{M}^{S_3})_y|^2}{\mb{M}_\star}d\xi
dy\right)^{\f{1}{2}}\\[5mm]
\di\leq C\chi\Big[\chi+\sum_{|\alpha|=2}\|\partial^\a(v-V^{S_3},
u-U^{S_3},\t-\T^{S_3})\|\\
\qquad+\|(v-V^{S_3},
u-U^{S_3},\t-\T^{S_3})\big(\sum_{|\alpha|=2}\partial^{\a}+\sum_{|\a^\prime|=1}\partial^{\a^\prime}\big)(V^{S_3},U^{S_3},\T^{S_3})\|\\
\di \qquad+\sum_{|\a^\prime|=1}\|\partial^{\a^\prime}(v-V^{S_3},
u-U^{S_3},\t-\T^{S_3})\cdot\partial^{\a^\prime}(V^{S_3},U^{S_3},\T^{S_3})\|\\
\di \qquad
+\sum_{|\alpha^\prime|=1}\|\partial^{\a^\prime}(v-V^{S_3},
u-U^{S_3},\t-\T^{S_3})\|_{L^4(dy)}^2\Big]\leq C_{h,T}\chi^2.
\end{array}
\label{(4.17+)}
\end{equation}
From (\ref{(1.17)}) and (\ref{(2.41)}), we have
\begin{equation}
\left\{
\begin{array}{l}
\di \p_{\tau}-\psi_{1y}=0,\\
\di \psi_{1\tau}+(p-P)_y
              =-\f{4}{3}\Big[(\f{\mu(\T)}{V}U_{1y})_y-(\f{\mu(\T^{S_3})}{V^{S_3}}U^{S_3}_{1y})_y\Big]+\int\x_1^2\Pi^{CD}_{11y}d\x\\
              \di \hspace{5cm}-\bar Q_{1y}- Q_{1y}-\int\xi_1^2\wt{\mb{G}}_yd\xi,\\
\di \psi_{i\tau}=-(\f{\mu(\T)}{V}U_{iy})_y+\int\x_1\x_i\Pi^{CD}_{11y}d\x-\bar Q_{iy}- Q_{iy}-\int\xi_1\xi_i\wt{\mb{G}}_yd\xi,~~i=2,3,\\
\di
\o_{\tau}+(pu_{1}-PU_{1})_y=-\Big[(\f{\k(\T)}{V}\T_y)_y-(\f{\k(\T^{S_3})}{V^{S_3}}\T^{S_3}_y)_y\Big]
              -\f{4}{3}\Big[(\f{\mu(\T)}{V}U_1U_{1y})_y\\
              \di\qquad\qquad-(\f{\mu(\T^{S_3})}{V^{S_3}}U^{S_3}_1U^{S_3}_{1y})_y\Big]-\sum_{i=2}^3(\f{\mu(\T)}{V}U_iU_{iy})_y+\int\x_1\f{|\x|^2}{2}\Pi^{CD}_{11y}d\x\\
              \di\qquad\qquad-\bar Q_{4y}-Q_{4y}-\f{1}{2}\int\xi_1|\xi|^2\wt{\mb{G}}_yd\xi.
\end{array}
\right. \label{(4.18)}
\end{equation}
And from \eqref{barQ1}, we get
\begin{equation}
\begin{array}{ll}
\di \bar Q_{1y}=\v[\bar
P-P^{R_1}-P^{CD}-P^{S_3}-(p_v^{R_1}d_1+p_{u_1}^{R_1}d_2+p_E^{R_1}d_3)]_x\\
\di \qquad\quad
+\v^2(\f{\mu(\bar \T)\bar U_{1x}}{\bar V}-\f{\mu(\T^{R_1})U^{R_1}_{1x}}{V^{R_1}}-\f{\mu(\T^{CD})U^{CD}_{1x}}{V^{CD}}-\f{\mu(\T^{S_3})U^{S_3}_{1x}}{V^{S_3}})_x\\
\quad~~\di= O(1)\v\Big[|(V^{R_1}_x,U^{R_1}_x,\mathcal{E}^{R_1}_x,d_{ix})(|d_i|^2,V^{S_3}-v_*,U^{S_3}-u_*,\mathcal{E}^{S_3}-E_*)|+|d_id_{ix}|\Big]\\
\di~ \quad~~~\di
+O(1)\v^2\Big[|(d_{2xx},d_{1x}d_{2x},V^{R_1}_xd_{2x},U^{R_1}_xd_{1x})|+|(U^{R_1}_{xx},U^{R_1}_xV^{R_1}_x)(d_1,V^{S_3}-v_*)|\\
\qquad~~~\di
+|(U^{S_3}_{xx},U^{S_3}_xV^{S_3}_x)(V^{R_1}-v_*,d_1)|+|(U^{S_3}_x,V^{S_3}_x)(U^{R_1}_x,V^{R_1}_x,d_{1x},d_{2x})|\Big]\\
\quad~\di:=\bar Q_{13}+\bar Q_{14},
\end{array}\label{Q1y-e}
\end{equation}
where $\bar Q_{13}$ represents the wave interaction satisfying
\begin{equation}
\begin{array}{ll}
\di \bar Q_{13}&\di=O(1)\v\Big[|(V^{R_1}_x,U^{R_1}_x,\mathcal{E}^{R_1}_x,d_{ix})(V^{S_3}-v_*,U^{S_3}-u_*,\mathcal{E}^{S_3}-E_*)|\Big]\\
&\di~~+O(1)\v^2\Big[|(U^{R_1}_{xx},U^{R_1}_xV^{R_1}_x)(d_1,V^{S_3}-v_*)|+|(U^{S_3}_{xx},U^{S_3}_xV^{S_3}_x)(V^{R_1}-v_*,d_1)|\\
&\di\qquad\qquad~~~
+|(U^{S_3}_x,V^{S_3}_x)(U^{R_1}_x,V^{R_1}_x,d_{1x},d_{2x})|\Big]\leq
C_{h,T}~e^{-\f{c|x|}{\s}}e^{-\f{C_h}{\s}},
\end{array}\label{Q13}
\end{equation}
and $\bar Q_{14}$ represents the terms related to the hyperbolic waves
$d_i~(i=1,2,3)$ satisfying
\begin{equation}
\begin{array}{ll}
\di \bar Q_{14}&\di=
O(1)\v\Big[|(V^{R_1}_x,U^{R_1}_x,\mathcal{E}^{R_1}_x,d_{ix})||d_i|^2+|d_id_{ix}|\Big]\\[3mm]
&\di~~
+O(1)\v^2\Big[|(d_{2xx},d_{1x}d_{2x},V^{R_1}_xd_{2x},U^{R_1}_xd_{1x})|\Big].
\end{array}
\label{Q14}
\end{equation}
Then we have
\begin{equation}
\int_{\f h\v}^\tau\int|\bar Q_{1y}|^2 dyd\tau\leq C\int_{\f
h\v}^\tau\int(|\bar Q_{13}|^2+|\bar Q_{14}|^2) dyd\tau \leq
C_{h,T}~\v^{\f12}.\label{barQ1-E}
\end{equation}
Similar estimates hold for $\bar Q_{iy}~(i=2,3,4).$

By \eqref{Q1-2}, straightforward calcuation gives
\begin{equation}
\begin{array}{ll}
\di Q_{1y} =\v\big[P-\bar P-(\bar p_v b_1+\bar p_u\cdot b_2+\bar p_E
b_3)\big]_x
-\f43\v^2\big(\f{\mu(\T)U_{1x}}{V}-\f{\mu(\bar \T)\bar U_{1x}}{\bar V}\big)_x\\
\qquad\di =O(1)\v\sum_{i=1}^3\Big[|b_i|^2|(\bar V_x,\bar
U_{1x},\bar{\mathcal{E}}_x)|+|b_i||b_{ix}|\Big]\\
\qquad~\di
+O(1)\v^2\sum_{i=1}^3\Big[|b_{ixx}|+|b_{ix}|^2+|b_i||(\bar
U_{xx},\bar V_x\bar U_x,\bar V_x\bar \T_x)|+|b_{ix}||(\bar V_x,\bar
U_x,\bar \T_x)|\Big].
\end{array}
\label{Q1-y}
\end{equation}
Hence,
\begin{equation}
\int_{\f h\v}^\tau\int |Q_{1y}|^2dyd\tau \leq
C_{h,T}~\v^{\f12}.\label{Q1-E}
\end{equation}
Similar estimates hold for $Q_{iy}~(i=2,3,4).$

Thus from the system \eqref{(4.18)}, we have
\begin{equation}
\|(\p_{\tau},\psi_{\tau},\o_{\tau})\|^2\leq C_{h,T}\chi^2,
\label{(4.19)}
\end{equation}
and
\begin{equation}
\|(\p_{\tau},\psi_{\tau},\z_{\tau})\|^2\leq
C\|(\p_{\tau},\psi_{\tau},\o_{\tau},U_{\tau}\cdot\psi)\|^2\leq
C_{h,T}\chi^2. \label{(4.19+)}
\end{equation}

Now we want to obtain the estimates on
$\|\partial^{\a}(\p,\psi,\zeta)\|^2$ for $|\a|=2$. For brevity, we
only calculate $|\partial_{yy}(v-V^{S_3},u-U^{S_3},\t-\T^{S_3})|$,
and the others can be estimated similarly. From \eqref{(1.6+)}
and \eqref{shock-fluid}, we have
\begin{eqnarray}
&&|\partial_{yy}(v-V^{S_3},u-U^{S_3},\t-\T^{S_3})|\nonumber\\
&&=
O(1)|(V_{yy}^{S_3},U_{yy}^{S_3},\T_{yy}^{S_3},(V_{y}^{S_3})^2,(U_{y}^{S_3})^2,(\T_{y}^{S_3})^2)|\cdot
|(v-V^{S_3},u-U^{S_3},\t-\T^{S_3})|\nonumber\\
&& \quad+O(1)|(V_{y}^{S_3},U_{y}^{S_3},\T_{y}^{S_3})|\cdot
|(v_{y}-V_{y}^{S_3},u_{y}-U_{y}^{S_3},\t_{y}-\T_{y}^{S_3})|\nonumber\\
&&\quad+O(1)|(v_{y}-V_{y}^{S_3},u_{y}-U_{y}^{S_3},\t_{y}-\T_{y}^{S_3})|^2+O(1)\sum_{i=0}^{4}\int|\varphi_{i}(\xi)|\cdot|\partial_{yy}\wt{f}|d\xi.\label{yy}
\end{eqnarray}
Therefore, we have
\begin{eqnarray}\label{alpha=2}
\|\partial^\a(\phi,\psi,\z)\|^2 &\leq&
 C\Big[\|\partial^\a(v-V^{S_3},u-U^{S_3},\t-\T^{S_3})\|^2+\|\partial^\a(V^{R_1},U^{R_1},\T^{R_1})\|^2\nonumber\\
&+&\|\partial^\a(V^{CD},U^{CD},\T^{CD})\|^2+\|\partial^\a(d_1,d_2,d_3)\|^2\nonumber\\
&+&\|\partial^\a
(b_1,b_{21},b_{22},b_{23},b_3)\|^2\Big]+C\int\int\f{|\partial^{\a}\wt{f}|^2}{\mb{M}_{\star}}d\xi
dy\leq C_{h,T}\chi^2.
\end{eqnarray}

Finally, by noticing the fact that $f=\mb{M}+\mb{G}$ and
$F^{S_3}=\mb{M}^{S_3}+\mb{G}^{S_3}$,  \eqref{alpha=2},
with $|\a|=2$ yields,
\begin{equation}
\begin{array}{ll}
\di \int\int\f{|\partial^\a\wt{\mb{G}}|^2}{\mb{M}_\star}d\xi dy\leq
C\int\int\f{|\partial^\a \wt f|^2}{\mb{M}_\star}d\xi dy+C
\int\int\f{|\partial^\a(\mb{M}-\mb{M}^{S_3})|^2}{\mb{M}_\star}d\xi
dy \leq C_{h,T}\chi^2,
\end{array}
\label{(4.25)}
\end{equation}
where in the last inequality we have used a similar argument used for
\eqref{(4.17+)}.

Before closing the a priori estimate (\ref{priori}), we list some
basic lemmas based on the celebrated H-theorem for later use. The
first lemma is from \cite{GPS}.

 \begin{lemma}\label{Lemma 4.1} There exists a positive
constant $C$ such that
$$
\int\f{\nu(|\xi|)^{-1}Q(f,g)^2}{\wt{\mb{M}}}d\xi\le
C\left\{\int\f{\nu(|\xi|)f^2}{\wt{\mb{M}}}d\xi\cdot\int\f{g^2}{\wt{\mb{M}}}d\xi+
\int\f{f^2}{\wt{\mb{M}}}d\xi\cdot\int\f{\nu(|\xi|)g^2}{\wt{\mb{M}}}d\xi\right\},
$$
where $\wt{\mb{M}}$ can be any Maxwellian so that the above
integrals are well-defined.
\end{lemma}

Based on Lemma \ref{Lemma 4.1}, the following three lemmas are taken
from \cite{Liu-Yang-Yu-Zhao}. And their proofs are straightforward by
using Cauchy inequality.

\begin{lemma}\label{Lemma 4.2} If $\t/2<\t_\star<\t$, then there exist two
positive constants $\wt\sigma=\wt\sigma(v,u,\t;\break
v_\star,u_\star,\t_\star)$ and
$\eta_0=\eta_0(v,u,\t;v_\star,u_\star,\t_\star)$ such that if
$|v-v_\star|+|u-u_\star|+|\t-\t_\star|<\eta_0$, we have for
$g(\xi)\in  \mathfrak{N}^\bot$,
$$
-\int\f{g\mb{L}_\mb{M}g}{\mb{M}_\star}d\xi\geq
\wt\sigma\int\f{\nu(|\xi|)g^2}{\mb{M}_\star}d\xi.
$$
\end{lemma}

 \begin{lemma}\label{Lemma 4.3} Under the assumptions in Lemma \ref{Lemma 4.2}, we
have  for each $g(\xi)\in  \mathfrak{N}^\bot$,
$$
 \int\f{\nu(|\xi|)}{\mb{M}}|\mb{L}_\mb{M}^{-1}g|^2d\xi
\leq \wt\sigma^{-2}\int\f{\nu(|\xi|)^{-1}g^2}{\mb{M}}d\xi,~~{\rm
and}~~
\int\f{\nu(|\xi|)}{\mb{M}_\star}|\mb{L}_\mb{M}^{-1}g|^2d\xi\le
\wt\sigma^{-2}\int\f{\nu(|\xi|)^{-1}g^2}{\mb{M}_\star}d\xi.
$$
\end{lemma}

\begin{lemma}\label{Lemma 4.4} Under the conditions in Lemma \ref{Lemma 4.2},  for any
positive constants $k$ and $\lambda$, it holds that
$$
|\int\f{g_1\mb{P}_1(|\xi|^kg_2)}{\mb{M}_\star}d\xi-\int\f{g_1|\xi|^kg_2}{\mb{M}_\star}d\xi|\le
C_{k,\lambda}\int\f{\lambda|g_1|^2+\lambda^{-1}|g_2|^2}{\mb{M}_\star}d\xi,
$$
where the constant $C_{k,\lambda}$ depends on $k$ and $\lambda$.
\end{lemma}

\subsection{Energy Estimates}

To close the a priori estimate \eqref{priori} and to prove Theorem
\ref{Theorem 4.1}, we need  the following energy estimates given in
Propositions \ref{Prop3.1} and Proposition \ref{Prop3.2}. First, the
lower order estimates to the system \eqref{sys} and \eqref{G1e} are
given in the following Proposition.

\begin{proposition} \label{Prop3.1}
Under the assumptions of Theorem \ref{Theorem 4.1}, there exist
positive constants $C$ and $C_{h,T}$  independent of $\v$
such that
\begin{equation*}
\begin{array}{ll}
\di \sup_{\f h\v\leq \tau_1\leq
\tau}\Big[\|(\Phi,\Psi,W,\Phi_y)(\tau_1,\cdot)\|^2+\int\int\f{|\wt{\mb{G}}_1|^2}{\mb{M}_\star}(\tau_1,y,\x)d\xi
dy\Big]\\
\di\quad+\int_{\f
h\v}^\tau\Big[\|\sqrt{|U^{S_3}_{1y}|}(\Psi,W)\|^2+\|(\Phi_y,\Psi_y,W_y,\z,\Psi_\tau,W_\tau)\|^2\Big]d\tau+\int_{\f
h\v}^\tau\int\int\f{\nu(|\xi|)}{\mb{M}_\star}|\wt{\mb{G}}_1|^2d\xi
dyd\tau\\
\di \leq C_{h,T}~\v\int_{\f
h\v}^\tau\|(\Psi,W)\|^2d\tau+C\sum_{|\a^\prime|=1}\int_{\f
h\v}^\tau\|\partial^{\a^\prime}(\p,\psi,\z)\|^2d\tau\\
\di\qquad+C\sum_{|\a^\prime|=1}\int_{\f
h\v}^\tau\int\int\f{\nu(|\xi|)}{\mb{M}_\star}|\partial^{\a^\prime}\wt{\mb{G}}|^2d\xi
dyd\tau+C_{h,T}~\v^{\f25}.
\end{array}
\end{equation*}
\end{proposition}

For brevity of the presentation, we also put the proof of Proposition \ref{Prop3.1} to the Appendices.

Then we perform the higher order estimates. Firstly, we apply
$\partial_y$ to the system \eqref{sys} to get the following system
for $(\phi,\psi,\zeta)$:
\begin{equation}
\left\{
\begin{array}{ll}
\di \phi_\tau-\psi_{1y}=0,\\[3mm]
\di \psi_{1\tau}-\f{Z}{V}\phi_y+\f{2}{3V}\z_y+H_1
=\big(\f43\f{\mu(\T)}{V}\psi_{1y}\big)_y-\int\xi_1^2(\Pi_{1}-\Pi^{CD}_{11}-\Pi^{S_3}_1)_yd\xi\\
\di \hspace{11cm}+N_5-\bar
Q_{1y}-Q_{1y},\\
\di \psi_{i\tau}+H_i=\big(\f{\mu(\T)}{V}\psi_{iy}\big)_y-\int\xi_1\xi_i(\Pi_{1}-\Pi^{CD}_{11}-\Pi^{S_3}_1)_yd\xi+N_{i+4}-\bar Q_{iy}-Q_{iy}, i=2,3,\\
\di \z_\tau+Z\psi_{1y}+H_4=
\big(\f{\k(\T)}{V}\z_{y}\big)_y-\int\xi_1\f{|\x|^2}{2}(\Pi_{1}-\Pi^{CD}_{11}-\Pi^{S_3}_1)_yd\xi+\sum_{i=1}^3\psi_i\int\x_1\x_i\Pi_{1y}d\x\\
\di\qquad
+\sum_{i=1}^3U_i\int\xi_1\x_i(\Pi_{1}-\Pi^{CD}_{11}-\Pi^{S_3}_1)_yd\xi+N_8-\bar
Q_{4y}-Q_{4y}+\sum_{i=1}^3U_i(\bar Q_{iy}+Q_{iy}),
\end{array} \right.\label{sys-h}
\end{equation}
where $Z$ is defined in \eqref{Z}. Here,  the linear terms are
\begin{equation}
H_1=-\f43\f{\mu^\prime(\T)}{V}U_{1y}\z_y-\big(\f{Z}{V}\big)_y\phi+(\f{2}{3V})_y\z-\big(\f{4\mu^\prime(\T)}{3V}U_{1y}\big)_y\z,
\label{H1}
\end{equation}
\begin{equation}
H_i=\f{\mu(\T)U_{iy}}{V^2}\phi_y-\f{\mu^\prime(\T)}{V}U_{iy}\z_y-\big(\f{\mu(\T)U_{iy}}{V^2}\big)_y\phi
-\big(\f{\mu^\prime(\T)}{V}U_{iy}\big)_y\z,~~i=2,3,
 \label{Hi}
\end{equation}
\begin{equation}
\begin{array}{ll}
H_4=&\di-\f{U_{1y}}{V}\big(Z\phi-\f{2}{3}\z+\f{4\mu^\prime(\T)U_{1y}}{3}\z+\f{4\mu(\T)}{3}\psi_{1y}\big)+\big(\f{\k(\T)}{V^2}\T_y\phi\big)_y\\
[3mm]
&\di-\sum_{i=2}^3\f{U_{iy}}{V}\big(-\f{\mu(\T)U_{iy}}{V}\phi+\mu^\prime(\T)U_{iy}\z+2\mu(\T)\psi_{iy}\big)-\big(\f{\k^\prime(\T)}{V}\T_{y}\z\big)_y,
\end{array}\label{H4}
\end{equation}
 and the nonlinear terms are
\begin{equation}
\begin{array}{ll}
N_5&\di=\Big[\f{p-P}{V}\phi+\f43(\f{\mu(\t)}{v}-\f{\mu(\T)}{V})\psi_{1y}
+\f43U_{1y}\Big(\f{\mu(\t)}{v}-\f{\mu(\T)}{V}+\f{\mu(\T)}{V^2}\phi-\f{\mu^\prime(\T)}{V}\zeta\Big)\Big]_y\\
&\di=O(1)\Big[|(\phi,\psi,\zeta)|^2+|(\phi_y,\psi_{1y},\z_y)|^2+|\psi_{1yy}|^2\Big],
\end{array}\label{N5}
\end{equation}
\begin{equation}
\begin{array}{ll}
N_{i+4}&\di=\Big[(\f{\mu(\t)}{v}-\f{\mu(\T)}{V})\psi_{iy}+U_{iy}\Big(\f{\mu(\t)}{v}-\f{\mu(\T)}{V}+\f{\mu(\T)}{V^2}\phi-\f{\mu^\prime(\T)}{V}\zeta\Big)\Big]_y\\
&\di=O(1)\Big[|(\phi,\psi,\zeta)|^2+|(\phi_y,\psi_{iy},\z_y)|^2+|\psi_{iyy}|^2\Big],~~i=2,3,
\end{array}\label{Ni+4}
\end{equation}
and
\begin{equation}
\begin{array}{ll}
N_8&\di=-(p-P)\psi_{1y}+\f{U_{1y}}{V}(p-P)\phi+2\big(\f{\mu(\t)}{v}-\f{\mu(\T)}{V}\big)\Big(\f43U_{1y}\psi_{1y}+\sum_{i=2}^3U_{iy}\psi_{iy}\Big)\\
&\di~~+\Big(\f{\mu(\t)}{v}-\f{\mu(\T)}{V}+\f{\mu(\T)}{V^2}\phi-\f{\mu^\prime(\T)}{V}\zeta\Big)\Big(\f43U_{1y}^2+\sum_{i=2}^3U_{iy}^2\Big)+\f{4\mu(\t)}{3v}\psi_{1y}^2\\
&\di~~
+\sum_{i=2}^3\f{\mu(\t)}{v}\psi_{iy}^2+\Big[(\f{\k(\t)}{v}-\f{\k(\T)}{V})\zeta_y+\T_{y}\Big(\f{\k(\t)}{v}-\f{\k(\T)}{V}+\f{\k(\T)}{V^2}\phi-\f{\k^\prime(\T)}{V}\zeta\Big)\Big]_y\\
&\di=O(1)\Big[|(\phi,\psi,\zeta)|^2+|(\phi_y,\psi_{y},\z_y)|^2+|\z_{yy}|^2\Big].
\end{array}\label{N8}
\end{equation}

To derive the estimate on the higher order derivatives,
applying $\partial_y$ to the system \eqref{sys-h}, gives
\begin{equation}
\left\{
\begin{array}{ll}
\di \phi_{y\tau}-\psi_{1yy}=0,\\[3mm]
\di \psi_{1y\tau}-\f{Z}{V}\phi_{yy}+\f{2}{3V}\z_{yy}+H_5
=\big(\f43\f{\mu(\T)}{V}\psi_{1y}\big)_{yy}\\[3mm]
 \di\hspace{2cm}-\int\xi_1^2(\Pi_{1}-\Pi^{CD}_{11}-\Pi^{S_3}_1)_{yy}d\xi+N_{5y}-\bar
Q_{1yy}-Q_{1yy},\\
\di \psi_{iy\tau}+H_{i+4}
=\big(\f{\mu(\T)}{V}\psi_{iy}\big)_{yy}-\int\xi_1\xi_i(\Pi_{1}-\Pi^{CD}_{11}-\Pi^{S_3}_1)_{yy}d\xi
\\[3mm]
 \di\hspace{6cm}+(N_{i+4})_y-\bar Q_{iyy}-Q_{iyy}, i=2,3,\\[3mm]
\di \z_{y\tau} +Z\psi_{1yy}+H_8=
\big(\f{\k(\T)}{V}\z_{y}\big)_{yy}-\int\xi_1\f{|\x|^2}{2}(\Pi_{1}-\Pi^{CD}_{11}-\Pi^{S_3}_1)_{yy}d\xi\\[3mm]
 \di\qquad\qquad +\sum_{i=1}^3\Big[\psi_i\int\x_1\x_i\Pi_{1y}d\x\Big]_y
+\sum_{i=1}^3\Big[U_i\int\xi_1\x_i(\Pi_{1}-\Pi^{CD}_{11}-\Pi^{S_3}_1)_yd\xi\Big]_y\\
\di\qquad\qquad +N_{8y}-\bar
Q_{4yy}-Q_{4yy}+\sum_{i=1}^3\Big[U_i(\bar Q_{iy}+Q_{iy})\Big]_y,
\end{array} \right.\label{sys-h-h}
\end{equation}
where
\begin{equation}
\begin{array}{ll}
H_5=&\di
-\f43\f{\mu^\prime(\T)}{V}U_{1y}\z_{yy}-2\big(\f{Z}{V}\big)_y\phi_y
-\big(\f{Z}{V}\big)_{yy}\phi
+(\f{4}{3V})_y\z_y\\
&\di
+(\f{2}{3V})_{yy}\z-\big(\f{8\mu^\prime(\T)}{3V}U_{1y}\big)_y\z_y
-\big(\f{4\mu^\prime(\T)}{3V}U_{1y}\big)_{yy}\z,
\end{array}
\label{H5}
\end{equation}
\begin{equation}
\begin{array}{ll}
H_{i+4}=&\di\f{\mu(\T)U_{iy}}{V^2}\phi_{yy}-\f{\mu^\prime(\T)}{V}U_{iy}\z_{yy}+2\big(\f{\mu(\T)U_{iy}}{V^2}\big)_y\phi_y+\big(\f{\mu(\T)U_{iy}}{V^2}\big)_{yy}\phi\\
&\di-2\big(\f{\mu^\prime(\T)}{V}U_{iy}\big)_y\z_y-\big(\f{\mu^\prime(\T)}{V}U_{iy}\big)_{yy}\z,~~i=2,3,
\end{array}\label{H4+i}
\end{equation}
and
\begin{equation}
\begin{array}{ll}
H_8&\di=Z_y\psi_{1y}-\Big[\f{U_{1y}}{V}\big(Z\phi-\f{2}{3}\z+\f{4\mu^\prime(\T)U_{1y}}{3}\z+\f{4\mu(\T)}{3}\psi_{1y}\big)\Big]_y+\big(\f{\k(\T)}{V^2}\T_y\phi\big)_{yy}\\
[3mm]
&\di-\sum_{i=2}^3\Big[\f{U_{iy}}{V}\big(-\f{\mu(\T)U_{iy}}{V}\phi+\mu^\prime(\T)U_{iy}\z+2\mu(\T)\psi_{iy}\big)\Big]_y-\big(\f{\k^\prime(\T)}{V}\T_{y}\z\big)_{yy}.
\end{array}\label{H8}
\end{equation}

By using the above two systems and the
equation for the non-fluid component, we can reach the following proposition
for the
higher order energy estimates.

\begin{proposition}\label{Prop3.2}
Under the assumptions of Theorem \ref{Theorem 4.1}, there exist
positive constants $C$ and $C_{h,T}$  independent of $\v$
such that
\begin{equation*}
\begin{array}{ll}
\di \sup_{\f h\v\leq \tau_1\leq
\tau}\Big[\|(\phi,\psi,\z,\phi_y,\psi_y,\z_y)(\tau_1,\cdot)\|^2+\sum_{|\a^\prime|=1}
\int\int\f{|\partial^{\a^\prime}\wt{\mb{G}}|^2}{\mb{M}_\star}(\tau_1,y,\x)d\xi
dy+\sum_{|\a|=2}\int\int \f{|\partial^\a
\wt{f}|^2}{2\mb{M}_\star}(\tau_1,y,\x)d\xi dy\Big]\\
\di \qquad+\int_{\f
h\v}^{\tau}\sum_{1\leq|\a|\leq2}\|\partial^\a(\p,\psi,\z)\|^2d\tau+\sum_{1\leq|\a|\leq2}\int_{\f
h\v}^\tau\int\int\f{\nu(|\xi|)}{\mb{M}_\star}|\partial^\a
\wt{\mb{G}}|^2d\xi dyd\tau\\
\di \leq C(\d+C_{h,T}\chi)\int_{\f h\v}^\tau\int\int
\f{\nu(|\xi|)}{\mb{M}_\star}|\wt{\mb{G}}_1|^2d\x dyd\tau+
C(\d+C_{h,T}\chi)\int_{\f h\v}^\tau\|(\p,\psi,\z)\|^2d\tau
+C_{h,T}~\v^{\f12}.
\end{array}
\end{equation*}
\end{proposition}

Again, the proof of Proposition \ref{Prop3.2} will be given in the Appendices.

By combining the above lower and higher order estimates given in Propositions
\ref{Prop3.1} and \ref{Prop3.2} and choosing the wave strength $\d$,
the bound on the a priori estimate $\chi$ and the Knudsen number
$\v$ to be suitably small, we obtain
$$
\begin{array}{ll}
\di \mathcal{N}(\tau)+\int_{\f
h\v}^\tau\Big[\sum_{0\leq|\a|\leq2}\|\partial^\a(\p,\psi,\z)\|^2+\|\sqrt{|U^{S_3}_{1y}|}(\Psi,W)\|^2\Big]d\tau
+\int_{\f h\v}^\tau\int\int\f{\nu(|\x|)|\wt{\mb{G}}
_1|^2}{\mb{M}_\star}d\xi dyd\tau\\
\di +\sum_{1\leq|\a|\leq2} \int_{\f
h\v}^\tau\int\int\f{\nu(|\x|)|\partial^{\a} \wt{\mb{G}}
|^2}{\mb{M}_\star}(\tau,y,\x)d\xi dyd\tau\leq C_{h,T}~\v^{\f25}.
\end{array}
$$
Therefore, we close the a priori assumption
\eqref{priori} and then complete the proof of Theorem \ref{Theorem 4.1}.

\section{Appendices}
\renewcommand{\theequation}{\arabic{section}.\arabic{subsection}.\arabic{equation}}
As mentioned before, since the proofs of Propositions 3.1 and 3.2
are technical and long, we put them in the following three
subsections.

\subsection{Proof of Proposition \ref{Prop3.1}}
\setcounter{equation}{0}

\underline{Proof of Proposition 3.1.} Firstly, from the fact that
$$
v_-<V^{R_1}<v_*,~~V^{CD}\in (\min\{v_*,v^*\},\max\{v_*,v^*\})~~{\rm
and}~~ v^*<V^{S_3}<v_+,
$$
we have
$$
v_--\d^{CD}-\|(d_1,b_1)\|_{L^\i}\leq V\leq
v_++\d^{CD}+\|(d_1,b_1)\|_{L^\i}.
$$
Thus
\begin{equation}
\f{v_-}{2}\leq V\leq 2v_+,\quad {\rm if}\quad \v<\!<1~~{\rm
and}~~\d^{CD}\ll1.\label{fact}
\end{equation}
From $\eqref{shock-profile}_1$ and $\eqref{shock-profile}_2$, we can
obtain
\begin{equation}
Z^{S_3}:=P^{S_3}-\f43\f{\mu(\T^{S_3})U^{S_3}_{1y}}{V^{S_3}}=(a_1-s_3^2V^{S_3})-\int\x^2_1\Pi^{S_3}_1d\x,
\label{ZS3}
\end{equation}
where $a_1=p_++s_3^2v_+=p^*+s_3^2v^*.$

Since
\begin{equation}
0<p_+\leq (a_1-s_3^2V^{S_3})\leq p^*,~~{\rm and}~~
|\int\x^2_1\Pi^{S_3}_1d\x|\leq C\d^{S_3},
\end{equation}
we have
$$
p_+-C\d^{S_3}\leq Z^{S_3}\leq p^*+C\d^{S_3},
$$
and then
\begin{equation}
\f{p_+}{2}\leq Z^{S_3}\leq 2p^*,~~{\rm if}~~\d^{S_3}\ll1.
\end{equation}
Now we estimate $Z$ defined in \eqref{Z} as follows.
\begin{equation}
\begin{array}{ll}
Z\di =Z^{S_3}+(P^{R_1}-p_*)-\f43\f{\mu(\T^{R_1})U^{R_1}_{1y}}{V^{R_1}}+(P^{CD}-p^*)-\f43\f{\mu(\T^{CD})U^{CD}_{1y}}{V^{CD}}\\
\di \qquad+(P-P^{R_1}-P^{CD}-P^{S_3}+p_*+p^*)\\
\di \qquad-\f43\Big[\f{\mu(\T)U_{1y}}{V}-\f{\mu(\T^{R_1})U^{R_1}_{1y}}{V^{R_1}}-\f{\mu(\T^{CD})U^{CD}_{1y}}{V^{CD}}-\f{\mu(\T^{S_3})U^{S_3}_{1y}}{V^{S_3}}\Big]\\
\di
\quad:=Z^{S_3}+(P^{R_1}-p_*)-\f43\f{\mu(\T^{R_1})U^{R_1}_{1y}}{V^{R_1}}+(P^{CD}-p^*)-\f43\f{\mu(\T^{CD})U^{CD}_{1y}}{V^{CD}}+Q_5,
\end{array}\label{Z1}
\end{equation}
where $Q_5$ is related to  the hyperbolic waves and  the wave
interaction terms given by
\begin{equation}
\begin{array}{ll}
Q_5
&\di
=O(1)\sum_{i=1}^3|(d_i,b_i)|+O(1)|(d_{2y},b_{2y})|+O(1)e^{-\f{C_h}{\s}}\\
&\di =O(1)(\f\v\s+\v^{\f34}+\f{\v^2}{\s^2}+e^{-\f{C_h}{\s}}),
\end{array}
\label{Q5}
\end{equation}
with $\s=\v^{\f15}.$

 From the
properties of the approximate rarefaction wave, we have
\begin{equation}
0\leq P^{R_1}-p_*\leq p_--p_*, \quad{\rm and }~~
\f43\f{\mu(\T^{R_1})U^{R_1}_{1y}}{V^{R_1}}\leq C_{h,T}~\v.
\label{ZR1}
\end{equation}
 And the
properties of the viscous contact wave imply that
\begin{equation}
P^{CD}-p^*\leq C_{h,T}~\v, \quad{\rm and }~~
\f43\f{\mu(\T^{CD})U^{CD}_{1y}}{V^{CD}}\leq C_{h,T}~\v^{\f32}.
\label{ZCD}
\end{equation}
Substituting \eqref{ZS3}, \eqref{Q5}, \eqref{ZR1} and \eqref{ZCD}
into \eqref{Z1}, we know that there exist positive constants $c$ and
$C$ such that
\begin{equation}
0<c\leq Z\leq C,\label{Z-bound}
\end{equation}
 provided $\v\ll1$ and the wave strength $\d\ll1$.

\

\underline{Step 1. Estimation on
$\|(\Phi,\Psi,W)(\tau,\cdot)\|^2$.}

\

 By
multiplying $(\ref{sys})_1$ by $\Phi$, $(\ref{sys})_2$ by
$\f{V}{Z}\Psi_1$, $(\ref{sys})_3$ by $V^{S_3}\Psi_i$,
$(\ref{sys})_4$ by $\f{2}{3Z^2}W$ respectively and adding all of them
together, we have
\begin{equation}
\begin{array}{l}
\di
I_1(\Phi,\Psi,W)_\tau+I_2(\Psi,W)+I_3(\Psi_y,W_y)=I_4(\Psi,W,\Phi_y,\Psi_y,W_y)\\
\di \qquad +\f{V\Psi_1}{Z}(N_1-\bar
Q_1-Q_1)+V^{S_3}\sum_{i=2}^3(N_i-\bar
Q_i-Q_i)\Psi_i\\
\di \qquad +\f{2W}{3Z^2}(N_4-\bar Q_4+\sum_{i=1}^3U_i\bar
Q_i-Q_4+\sum_{i=1}^3U_iQ_i)+K_1 \di+(\cdots)_y,
\end{array}\label{le1}
\end{equation}
where
\begin{equation}
\begin{array}{ll}
\di I_1(\Phi,\Psi,W)
=\f{\Phi^2}2+\f{V}{2Z}\Psi_1^2+\f{V^{S_3}}{2}\sum_{i=2}^3\Psi_i^2+\f{W^2}{3Z^2},\\
\di I_2(\Psi,W) =\Big[\f{2}{3Z}U_{1y}-\f12 (\f
VZ)_{\tau}\Big]\Psi_1^2-V^{S_3}_\tau\sum_{i=2}^3\f{\Psi_i^2}{2}-(\f{1}{3Z^2})_\tau
W^2,\\
\di I_3(\Psi_y,W_y)
=\f{4\mu(\T)}{3Z}\Psi_{1y}^2+\f{V^{S_3}\mu(\T)}{V}\sum_{i=2}^3\Psi_{iy}^2
+\f{2\k(\T)}{3Z^2 V}W_y^2,
\end{array}
\end{equation}
\begin{equation}
\begin{array}{ll} \di I_4(\Psi,W,\Phi_y,\Psi_y,W_y) =
-\f{2}{3Z^2}W\Psi_1(U_{1\tau}+Z_y)-\f{2W}{3Z^2}\sum_{i=2}^3U_{i\tau}\Psi_i\\
\di\qquad~~-\Big[\f{2\Psi_1}{3Z}+\f{\mu(\T)V^{S_3}}{V^2}\Phi_y-\f{\mu^\prime(\T)V^{S_3}}{V}(W_y+U_y\cdot\Psi)\Big]\sum_{i=2}^3U_{iy}\Psi_i\\
\di\qquad~~
-\f{4}{3}(\f{\mu(\T)}{Z})_y\Psi_{1}\Psi_{1y}-\sum_{i=2}^3(\f{V^{S_3}\mu(\T)}{V})_y\Psi_{i}\Psi_{iy}
-\f{2}{3}(\f{\k(\T)}{VZ^2})_yWW_y
\\
\di\qquad~~+\f43\f{\mu^\prime(\T)}{Z}U_{1y}\Psi_1(W_y+U_{y}\cdot\Psi)-\f{2\k(\T)}{3V^2Z^2}\T_y\Phi_yW+\f{2\k(\T)}{3VZ^2}W(U_{y}\cdot\Psi)_y
\\[3mm]
\di\qquad~~+\f{2\k^\prime(\T)}{3VZ^2}\T_yW(W_y+U_{y}\cdot\Psi)-\sum_{i=2}^3\f{2\mu(\T)U_{iy}}{3Z^2V}W\Psi_{iy}:=\sum_{i=1}^{11}I_4^i,
\end{array}\label{I4}
\end{equation}
and $K_1$ denotes the non-fluid parts given by
\begin{equation}
\begin{array}{ll}
\di K_1=-
\f{V}{Z}\Psi_1\int\xi_1^2(\Pi_1-\Pi^{CD}_{11}-\Pi_1^{S_3})d\xi-\sum_{i=2}^3V^{S_3}\Psi_i\int\xi_1\xi_i(\Pi_1-\Pi^{CD}_{11}-\Pi_1^{S_3})d\xi\\
\quad\di
+\f{2W}{3Z^2}\Big[-\int\f12\xi_1|\xi|^2(\Pi_1-\Pi^{CD}_{11}-\Pi_1^{S_3})d\xi+\sum_{i=1}^3U_i\int\xi_1\xi_i(\Pi_1-\Pi^{CD}_{11}-\Pi_1^{S_3})d\xi\Big].
\end{array} \label{K1}
\end{equation}

From \eqref{Z-bound}, we have
\begin{equation}
c(\Phi,\Psi,W)^2\leq I_1\leq C (\Phi,\Psi,W)^2,
\end{equation}
 for some positive
constants $c$ and $C$.

Now we estimate $I_i,~(i=2,3,4)$ term by term. Note that
\begin{equation}
Z^{S_3}_\tau=-s_3^2U^{S_3}_{1y}+s_3\int\x^2_1\Pi^{S_3}_{1y}d\x,
\end{equation}
and
\begin{equation}
\begin{array}{ll}
\di Z_\tau=-s_3^2U^{S_3}_{1y}+s_3\int\x^2_1\Pi^{S_3}_{1y}d\x
+Z^{CD}_\tau+Z^{R_1}_\tau+Q_6,
\end{array}\label{Z-tau-1}
\end{equation}
where $\di
Z^{R_1}=P^{R_1}-\f43\f{\mu(\T^{R_1})U^{R_1}_{1y}}{V^{R_1}}$, $\di
Z^{CD}=P^{CD}-\f43\f{\mu(\T^{CD})U^{CD}_{1y}}{V^{CD}}$ and
\begin{equation}
\begin{array}{ll}
Q_6&\di =(P-P^{R_1}-P^{CD}-P^{S_3})_\tau\\
&\di \qquad
-\f43\Big[\f{\mu(\T)U_{1y}}{V}-\f{\mu(\T^{R_1})U^{R_1}_{1y}}{V^{R_1}}-\f{\mu(\T^{CD})U^{CD}_{1y}}{V^{CD}}
-\f{\mu(\T^{S_3})U^{S_3}_{1y}}{V^{S_3}}\Big]_\tau\\[3mm]
&\di
=\f23\Big(\f{\T_\tau}{V}-\f{\T^{R_1}_\tau}{V^{R_1}}-\f{\T^{CD}_\tau}{V^{CD}}-\f{\T^{S_3}_\tau}{V^{S_3}}\Big)-\f23\Big(\f{\T
V_\tau}{V^2}-\f{\T^{R_1}V^{R_1}_\tau}{(V^{R_1})^2}-\f{\T^{CD}V^{CD}_\tau}{(V^{CD})^2}-\f{\T^{S_3}V^{S_3}_\tau}{(V^{S_3})^2}\Big)\\[3mm]
&\di\quad
-\f43\Big(\f{\mu(\T)U_{1y\tau}}{V}-\f{\mu(\T^{R_1})U^{R_1}_{1y\tau}}{V^{R_1}}-\f{\mu(\T^{CD})U^{CD}_{1y\tau}}{V^{CD}}-\f{\mu(\T^{S_3})U^{S_3}_{1y\tau}}{V^{S_3}}\Big)\\
&\di\quad
+\f43\Big[\f{\mu(\T)U_{1y}V_\tau}{V^2}-\f{\mu(\T^{R_1})U^{R_1}_{1y}V^{R_1}_\tau}{(V^{R_1})^2}-\f{\mu(\T^{CD})U^{CD}_{1y}V^{CD}_\tau}{(V^{CD})^2}-\f{\mu(\T^{S_3})U^{S_3}_{1y}V^{S_3}_\tau}{(V^{S_3})^2}\Big]\\
&\di\quad-\f43\Big[\f{\mu^\prime(\T)U_{1y}\T_\tau}{V}-\f{\mu^\prime(\T^{R_1})U^{R_1}_{1y}\T^{R_1}_\tau}{V^{R_1}}-\f{\mu^\prime(\T^{CD})U^{CD}_{1y}\T^{CD}_\tau}{V^{CD}}-\f{\mu^\prime(\T^{S_3})U^{S_3}_{1y}\T^{S_3}_\tau}{V^{S_3}}\Big]:=\sum_{i=1}^5Q_{6i}.
\end{array}\label{Q6}
\end{equation}
Then we get for the first term of $I_2$ that
\begin{equation}
\begin{array}{ll}
\di \f{2}{3Z}U_{1y}-\f12(\f VZ)_{\tau} =\f{1}{6Z}U_{1y}+\f{V}{2Z^2}Z_\tau\\[3mm]
\di =\f{U^{S_3}_{1y}}{6Z^2}(Z^{S_3}-3s_3^2V^{S_3})+\f{U^{S_3}_{1y}}{6Z^2}\big[(Z-Z^{S_3})-3s_3^2(V-V^{S_3})\big]+\f{s_3V}{2Z^2}\int\x^2_1\Pi^{S_3}_{1y}d\x\\
\di \quad+\f{1}{6Z}\big(U^{R_1}_{1y}+U^{CD}_{1y}+d_{2y}+b_{21y}\big)+\f{V}{2Z^2}\big(Z^{CD}_\tau+Z^{R_1}_\tau+\sum_{i=1}^5Q_{6i}\big)\\[3mm]
\di=\f{|U^{S_3}_{1y}|}{6Z^2}\big[-(a_1-s_3^2V^{S_3})+3s_3^2V^{S_3}+\int\x^2_1\Pi^{S_3}_1d\x\big]+\f{s_3V}{2Z^2}\int\x^2_1\Pi^{S_3}_{1y}d\x+Q_7,
\end{array}\label{Z-Psi-1}
\end{equation}
where
\begin{equation}
\begin{array}{ll}
Q_7=&\di \f{U^{S_3}_{1y}}{6Z^2}\big[(Z-Z^{S_3})-3s_3^2(V-V^{S_3})\big]+\f{1}{6Z}\big(U^{R_1}_{1y}+U^{CD}_{1y}+d_{2y}+b_{21y}\big)\\
&\di+\f{V}{2Z^2}\big(Z^{CD}_\tau+Z^{R_1}_\tau+\sum_{i=1}^5Q_{6i}\big).
\end{array}
\label{Q7}
\end{equation}
Firstly, straightforward calculation gives
$$
-(a_1-s_3^2V^{S_3})+3s_3^2V^{S_3}\geq 4p_*-C\d^{S_3}\geq 3p_*,~~{\rm
if}~~\d^{S_3}\ll1,
$$
$$
|\int\x^2_1\Pi^{S_3}_1d\x|\leq C\d^{S_3},~~{\rm
and}~~|\int\x^2_1\Pi^{S_3}_{1y}d\x|\leq C\d^{S_3}|U^{S_3}_{1y}|.
$$
And \eqref{Z-Psi-1} implies that
\begin{equation}
\Big[\f{2}{3Z}U_{1y}-\f12(\f VZ)_{\tau}\Big]\Psi_1^2\geq
C^{-1}|U^{S_3}_{1y}|\Psi_1^2+Q_7\Psi_1^2, \label{p1}
\end{equation}
provided $\d\ll1.$

By the definition of the approximate wave pattern defined in \eqref{wave} and
\eqref{Theta}, we have
\begin{equation}
\begin{array}{ll}
\di
Q_{61}=O(1)\v\Big[|(\T^{R_1}_t,U^{R_1}_{1t})||(V^{CD}-v_*,U^{CD}_1-u_{1*},V^{S_3}-v^*,U_1^{S_3}-u_{1}^*)|+|(\T^{R_1}_td_1,U^{R_1}_{1t}d_2)|\\
\di\qquad\qquad~~~+|(\T^{CD}_t,U^{CD}_{1t})||(V^{R_1}-v_*,U_1^{R_1}-u_{1*},V^{S_3}-v^*,U_1^{S_3}-u_{1}^*,d_1,d_2)|\\
\di\qquad\qquad~~~+|(\T^{S_3}_t,U^{S_3}_{1t})||(V^{R_1}-v_*,U_1^{R_1}-u_{1*},V^{CD}-v^*,U^{CD}_1-u^*_{1},d_1,d_2)|\\
\di\qquad\qquad~~~+|(d_{2t},d_{3t})|+|b_1||(\T^{R_1}_t,\T^{CD}_t,\T^{S_3}_t)|\Big]+\f{2}{3V}(b_{3}-\bar U\cdot b_2-\f{|b_2|^2}{2})_\tau\\
 \di\quad\leq C_{h,T}~\v+C|(b_2,b_{2\tau},b_{3\tau})|.
\end{array}\label{Q61}
\end{equation}
On the other hand, direct calculation gives
\begin{equation}
\begin{array}{ll}
\di Q_{63}=O(1)\v^2\Big[|U^{R_1}_{1xt}||(V^{CD}-v_*,U^{CD}_1-u_{1*},\T^{CD}-\t_*,V^{S_3}-v^*,U^{S_3}_1-u_1^*,\T^{S_3}-\t^*)|\\
\di\quad+|U^{CD}_{1xt}||(V^{R_1}-v_{*},U^{R_1}_1-u_{1*},\T^{R_1}-\t_{*},V^{S_3}-v^*,U^{S_3}_1-u_{1}^*,\T^{S_3}-\t^*,d_1,d_2,d_{3})|\\
\di\quad+|U^{S_3}_{1xt}||(V^{R_1}-v_{*},U^{R_1}_1-u_{1*},\T^{R_1}-\t_{*},V^{CD}-v^*,U^{CD}_1-u_{1}^*,\T^{CD}-\t^*,d_1,d_2,d_{3})|\\
\di\quad+|U^{R_1}_{1xt}||(d_1,d_{2},d_3)|+|(b_1,b_2,b_3)||(U^{R_1}_{1xt},U^{CD}_{1xt},U^{S_3}_{1xt})|+|d_{2xt}|\Big]-\f{4\mu(\T)}{3V}b_{21y\tau}\\
\di\leq C_{h,T}~\v+C|b_{21y\tau}|.
\end{array}\label{Q63}
\end{equation}
Similar estimates hold for $Q_{62}$, $Q_{64}$ and $Q_{65}$.

Moreover,  we have
\begin{equation}
\begin{array}{ll}
\di |Z^{R_1}_\tau|
=|(P^{R_1}-\f43\f{\mu(\T^{R_1})U^{R_1}_{1y}}{V^{R_1}})_\tau|\leq
C_{h,T}~\v,
\end{array}
\label{Z-R1-tau}
\end{equation}
and
\begin{equation}
\begin{array}{ll}
\di |Z^{CD}_\tau|
=|(P^{CD}-\f43\f{\mu(\T^{CD})U^{CD}_{1y}}{V^{CD}})_\tau|
 \leq C_{h,T}~\v^2.
\end{array}
\label{Z-CD-tau}
\end{equation}
Now from \eqref{Z1}, \eqref{Q7} and \eqref{Q61}-\eqref{Z-CD-tau}, we
have
\begin{equation}
\begin{array}{ll}
\di \int_{\f h\v}^\tau\int |Q_7|\Psi_1^2dyd\tau \leq
C_{h,T}~\v\int_{\f h\v}^\tau\int \Psi_1^2dyd\tau
+C\int_{\f h\v}^\tau\int|(b,b_y,b_\tau,b_{y\tau})||\Psi_1|^2dyd\tau\\
\di\qquad\quad \leq C_{h,T}~\v\int_{\f h\v}^\tau\int \Psi_1^2dyd\tau
+C\int_{\f h\v}^\tau\|(b,b_y,b_\tau,b_{y\tau})\|\|\Psi_1\|\|\Psi_1\|_{L^\i}d\tau\\
\di\qquad\quad \leq C_{h,T}~\v\int_{\f h\v}^\tau\int \Psi_1^2dyd\tau
+C_{h,T}~\v^{\f34}\int_{\f h\v}^\tau\|\Psi_1\|^{\f32}\|\Psi_{1y}\|^{\f12}d\tau\\
\di\qquad\quad \leq \b\int_{\f
h\v}^\tau\|\Psi_{1y}\|^2d\tau+C_{h,T,\b}~\v\int_{\f h\v}^\tau\int
\Psi_1^2dyd\tau.
\end{array}\label{Q7-Psi-1}
\end{equation}
Substituting \eqref{Q7-Psi-1} into \eqref{p1} gives that
\begin{equation}
\begin{array}{ll}
\di\int_{\f h\v}^\tau\int\Big[\f{2}{3Z}U_{1y}-\f12(\f
VZ)_{\tau}\Big]\Psi_1^2dyd\tau\\
\di \geq C^{-1}\int_{\f
h\v}^\tau\int|U^{S_3}_{1y}|\Psi_1^2dyd\tau-\b\int_{\f
h\v}^\tau\|\Psi_{1y}\|^2d\tau-C_{h,T,\b}~\v\int_{\f h\v}^\tau\int
\Psi_1^2dyd\tau,
\end{array}  \label{p2}
\end{equation}
provided $\d^{S_3}\ll1.$

On the other hand, we have
$$
-V^{S_3}_\tau\sum_{i=2}^3\f{\Psi_i^2}{2}=\sum_{i=2}^3|U^{S_3}_{1y}|\f{\Psi_i^2}{2}.
$$
Similar to \eqref{p2}, we have
$$
\begin{array}{ll}
\di\int_{\f h\v}^\tau\int\Big[-(\f{1}{3Z^2})_\tau
W^2\Big]dyd\tau=\int_{\f h\v}^\tau\int\f{2}{3Z^3}Z_\tau W^2dyd\tau\\
\di =\int_{\f
h\v}^\tau\int\f{2}{3Z^3}\Big(-s_3^2U^{S_3}_{1y}+s_3\int\x^2_1\Pi^{S_3}_{1y}d\x
+Z^{CD}_\tau+Z^{R_1}_\tau+Q_6\Big)W^2dyd\tau\\
\di \geq C^{-1}\int_{\f
h\v}^\tau\int|U^{S_3}_{1y}|W^2dyd\tau-\b\int_{\f
h\v}^\tau\|W_{y}\|^2d\tau-C_{h,T,\b}~\v\int_{\f h\v}^\tau\int
W^2dyd\tau.
\end{array}
$$
Thus we complete the estimation for $I_2$ as
\begin{equation}
\begin{array}{ll}
\di \int_{\f h\v}^\tau\int I_2(\Psi,W)dyd\tau\geq C^{-1}\int_{\f
h\v}^\tau\int|U^{S_3}_{1y}||(\Psi,W)|^2dyd\tau\\
\di \qquad\qquad\qquad\qquad-\b\int_{\f
h\v}^\tau\|(\Psi_{1y},W_{y})\|^2d\tau-C_{h,T,\b}~\v\int_{\f
h\v}^\tau\int |(\Psi,W)|^2dyd\tau.
\end{array}
\end{equation}
Direct computation yields that
\begin{equation}
I_3(\Psi_y,W_y)\geq C^{-1}|(\Psi_y,W_y)|^2.
\end{equation}
Now we estimate $I_4$ in \eqref{I4}. Note that
$$
U_{1\tau}+Z_y=-\int\x_1^2\Pi_{11y}^{CD}d\x-\int\x_1^2\Pi_{1y}^{S_3}d\x+\bar
Q_{1y}+Q_{1y}.
$$
Thus
$$
\begin{array}{ll}
I_4^1&\di=\f{2}{3Z^2}W\Psi_1\Big[\int\x_1^2\Pi_{11y}^{CD}d\x+\int\x_1^2\Pi_{1y}^{S_3}d\x-\bar
Q_{1y}-Q_{1y}\Big]\\
&\di=\f{2}{3Z^2}W\Psi_1\Big(\int\x_1^2\Pi_{11y}^{CD}d\x+\int\x_1^2\Pi_{1y}^{S_3}d\x\Big)+(\f{2W\Psi_1}{3Z^2})_y(\bar
Q_{1}+Q_{1})+(\cdots)_y\\
&\di
\leq\b|(\Psi_{1y},W_y)|^2+C_\b\d^{S_3}|U^{S_3}_{1y}||(\Psi_1,W)|^2+C_{h,T,\b}~\v|(\Psi_1,W)|^2+(\cdots)_y,
\end{array}
$$
where from now on,  $\beta $ is a small positive constant to be
determined and $C_\b$ is some positive constant depending on $\b$
but independent of $h,T$ and $\v$, while $C_{h,T,\b}$ depends on on
$h,T,\b$ but independent of $\v$.

Similarly, we have
$$
\begin{array}{ll}
I_4^2&\di=\sum_{i=2}^3\f{2W\Psi_i}{3Z^2}\Big[-\big(\f{\mu(\T)U_{iy}}{V}\big)_y+\int\x_1\x_i\Pi_{11y}^{CD}d\x+\int\x_1\x_i\Pi_{1y}^{S_3}d\x-\bar
Q_{iy}-Q_{iy}\Big]\\
&\di=\sum_{i=2}^3\f{2W\Psi_i}{3Z^2}\Big(\int\x_1\x_i\Pi_{11y}^{CD}d\x+\int\x_1\x_i\Pi_{1y}^{S_3}d\x\Big)+\sum_{i=2}^3(\f{2W\Psi_i}{3Z^2})_y\big(\f{\mu(\T)U_{iy}}{V}+\bar
Q_{i}+Q_{i}\big)+(\cdots)_y\\
&\di
\leq\b\sum_{i=2}^3|(\Psi_{iy},W_y)|^2+C_\b\d^{S_3}\sum_{i=2}^3|U^{S_3}_{1y}||(\Psi_i,W)|^2+C_{h,T,\b}~\v\sum_{i=2}^3|(\Psi_i,W)|^2+(\cdots)_y.
\end{array}
$$
Next, we can show that
\begin{equation}
\begin{array}{ll}
I_4^3&\di
=\sum_{i=2}^3U_i(\f{2\Psi_1\Psi_i}{3Z})_y-\Big[\f{\mu(\T)V^{S_3}}{V^2}\Phi_y-\f{\mu^\prime(\T)V^{S_3}}{V}(W_y+U_y\cdot\Psi)\Big]\sum_{i=2}^3U_{iy}\Psi_i+(\cdots)_y\\
&\di\leq\b|(\Phi_y,\Psi_{y},W_y)|^2+C_\b\d^{S_3}|U^{S_3}_{1y}||(\Psi,W)|^2+C_{h,T,\b}~\v|(\Psi,W)|^2+(\cdots)_y,
\end{array}
\end{equation}
where we have used the fact that $U_i=\bar
U_i+b_{i+1}=U^{CD}_i+b_{i+1}=O(1)\v^{\f12}, (i=2,3)$ and that
$U_{iy}=U^{CD}_{iy}+(b_{i+1})_y=O(1)\v^{\f34}$, $i=2,3.$

By H${\rm\ddot{o}}$lder inequality, we have
\begin{equation}
\sum_{i=4}^{11}I_4^i\leq\b|(\Phi_y,\Psi_{y},W_y)|^2+C_\b\d^{S_3}|U^{S_3}_{1y}||(\Psi,W)|^2+C_{h,T,\b}~\v|(\Psi,W)|^2.
\end{equation}
Hence, we have
\begin{equation}
\begin{array}{ll}
\di  I_4(\Psi,W,\Phi_y,\Psi_y,W_y)\\
\di \leq
\b|(\Phi_y,\Psi_y,W_y)|^2+(C_\b\d^{S_3}+C_{h,T}~\v)|U^{S_3}_{1y}||(\Psi,W)|^2+C_{h,T,\b}~\v|(\Psi,W)|^2+(\cdots)_y.
\end{array}
\end{equation}
Now we estimate the nonlinear terms $\f{V}{Z}N_1\Psi_1$ by
\begin{equation}
\begin{array}{ll}
\di \int_{\f h\v}^\tau\int\f{V}{Z}N_1\Psi_1dyd\tau
 &\di \leq
C \int_{\f
h\v}^\tau\int|\Psi_1|\Big[|\Phi_y|^2+|\Psi_y|^2+|\zeta|^2+|\Psi_{1yy}|^2\Big]dyd\tau\\
&\di \leq C\chi\int_{\f h\v}^\tau
\|(\Phi_y,\Psi_y,\zeta,\Psi_{1yy})\|^2d\tau.
\end{array}
\label{(4.17)}
\end{equation}

Similarly, we have
$$
\int_{\f h\v}^\tau\int
\sum_{i=2,3}V^{S_3}N_i\Psi_i+\f{2W}{3Z^2}N_4dyd\tau\leq
C\chi\int_{\f h\v}^\tau
\|(\Phi_y,\Psi_y,\zeta,\Psi_{yy},\zeta_y)\|^2d\tau.
$$

We now turn to  estimate the terms $-\f{V}{Z}\Psi_1 \bar Q_1$,
$-\sum_{i=2}^3 V^{S_3}\bar Q_i\Psi_i$ and $-\f{2W}{3Z^2}(\bar
Q_4-\sum_{i=1}^3U_i\bar Q_i)$. From (\ref{barQ1}) and \eqref{barQ4},
we have
\begin{equation}
\begin{array}{ll}
\di |-\f{V}{Z}\Psi_1 \bar Q_1|,|-\sum_{i=2}^3 V^{S_3}\bar
Q_i\Psi_i|,|-\f{2W}{3Z^2}(\bar
Q_4-\sum_{i=1}^3U_i\bar Q_i)|\\[5mm]
\di \leq C|(\Psi,W)||(\bar Q_{11},\bar Q_2,\bar Q_3,\bar
Q_{41})|+C|(\Psi,W)||(\bar Q_{12},\bar Q_{42})|.
\end{array}
\end{equation}
And  from (\ref{Q11}), we have
\begin{equation}
\begin{array}{ll}
\di \int_{\f h\v}^{\tau}\int_{\mathbf{R}}|(\Psi,W)||(\bar
Q_{11},\bar Q_2,\bar Q_3,\bar Q_{41})|dyd\tau
 \leq C\int_{\f h\v}^{\tau}\int_{\mathbf{R}}|(\Psi,W)|e^{\f{-c|x|}{\s}}e^{-\f{C_h}{\s}}dyd\tau\\
\di \leq \v\int_{\f h\v}^{\tau}\|(\Psi,W)\|^2d\tau
+C_{h,T}~e^{-\f{C_h}{2\s}}.
\end{array}
\end{equation}
On the other hand, from (\ref{barQ1})-\eqref{barQ4}, Lemma 2.3 and
noting that $\s=\v^{\f15}$, we have
\begin{equation}
\begin{array}{ll}
&\di \int_{\f h\v}^{\tau}\int_{\mathbf{R}}|(\Psi,W)||(\bar Q_{12},\bar Q_{42})|dyd\tau\\
&\di \leq C\int_{\f
h\v}^{\tau}\int_{\mathbf{R}}|(\Psi,W)|\Big[\sum_{i=1}^3|d_i|^2+\v|(d_{2x},d_{3x})|
+\v|(U^{R_1}_{1x},\T^{R_1}_x)||(d_1,d_2,d_3)|\Big]dyd\tau\\
&\di \leq C_{h,T}\int_{\f
h\v}^{\tau}\|(\Psi,W)\|_{L^2(dy)}\sum_{i=1}^3\Big[\|d_i\|_{L^\i}\|d_i\|_{L^2(dy)}+\v\|d_{ix}\|_{L^2(dy)}+\v
\|d_{i}\|_{L^2(dy)}\Big]d\tau\\
&\di \leq C_{h,T}(\f{\v}{\s})^{\f32}\int_{\f
h\v}^{\tau}\|(\Psi,W)\|_{L^2(dy)}d\tau
 \leq \v\int_{\f
h\v}^{\tau}\|(\Psi,W)\|^2d\tau+C_{h,T}~\v^{\f25}.
\end{array}
\end{equation}
Now we estimate the terms $-\f{V}{Z}\Psi_1 Q_1$, $-\sum_{i=2}^3
V^{S_3}Q_i\Psi_i$ and $-\f{2W}{3Z^2}(Q_4-\sum_{i=1}^3U_iQ_i)$. For this,
 from \eqref{Q1-2} and (\ref{Q11-21}), we have
\begin{equation}
\begin{array}{ll}
\di\int_{\f h\v}^{\tau}\int_{\mathbf{R}}\Big[-\f{V}{Z}\Psi_1 Q_{11}-\f{2W}{3Z^2}(Q_{41}-U_1Q_{11})\Big]dyd\tau\\
\di \leq C\int_{\f h\v}^{\tau}\int_{\mathbf{R}}|(\Psi_1,W)||(Q_{11},Q_{41})|dyd\tau \leq C\int_{\f h\v}^{\tau}\int_{\mathbf{R}}|(\Psi_1,W)|\sum_{i=1}^3|b_i|^2dyd\tau\\
\di\leq C\int_{\f h\v}^{\tau}\|(\Psi_1,W)\|_{L^2(dy)}\sum_{i=1}^3\|b_i\|_{L^\i}\|b_i\|_{L^2(dy)}d\tau\\
\di\leq C_{h,T}~\v^{\f32}\int_{\f h\v}^{\tau}\|(\Psi_1,W)\|d\tau
\leq \v\int_{\f h\v}^{\tau}\|(\Psi_1,W)\|^2d\tau+C_{h,T}~\v.
\end{array}
\end{equation}
On the other hand, from Lemma \ref{LemmaII} with the properties of
the hyperbolic wave II, we can show that
\begin{equation}
\begin{array}{ll}
\di \int_{\f h\v}^{\tau}\int_{\mathbf{R}}\Big[-\f{V}{Z}\Psi_1
Q_{12}-\sum_{i=2}^3 V^{S_3} Q_i\Psi_i -\f{2W}{3Z^2}(Q_{42}-U_1Q_{12}-\sum_{i=2}^3U_iQ_i)\Big]dyd\tau\\
\di =\int_{\f
h\v}^{\tau}\int_{\mathbf{R}}\Big[\Big(\f{4V\Psi_1}{3Z}-\f{8U_1W}{9Z^2}\Big)\Big(\f{\mu(\T)U_{1y}}{V}-\f{\mu(\bar\T)\bar
U_{1y}}{\bar
V}\Big)\\
\di
\qquad\qquad+\sum_{i=2}^3\Big(V^{S_3}\Psi_i-\f{2U_iW}{3Z^2}\Big)\Big(\f{\mu(\T)U_{iy}}{V}-\f{\mu(\bar\T)\bar
U_{iy}}{\bar
V}\Big)\\
\di
\qquad\qquad+\f{2W}{3Z^2}\Big(\f{\mu(\T)\T_{y}}{V}-\f{\mu(\bar\T)\bar
\T_{y}}{\bar
V}\Big)+\f{8W}{9Z^2}\Big(\f{\mu(\T)U_1U_{1y}}{V}-\f{\mu(\bar\T)\bar
U_1\bar U_{1y}}{\bar V}\Big)
\\
\qquad\qquad\di+\f{2W}{3Z^2}\sum_{i=2}^3\Big(\f{\mu(\T)U_iU_{iy}}{V}-\f{\mu(\bar\T)\bar
U_i\bar U_{iy}}{\bar V}\Big) \Big]dyd\tau:=\sum_{i=1}^5 I_{5i}.
\end{array}
\end{equation}
For brevity, we only compute $I_{51}$ in the following because
$I_{5i}~(i=2,3,4,5)$ can be estimated similarly.
\begin{equation}
\begin{array}{ll}
I_{51}&\di =\int_{\f h\v}^{\tau}\int_{\mathbf{R}}\Big(\f{4
V\Psi_1}{3Z}-\f{8
U_1W}{9Z^2}\Big)\Big(\f{\mu(\T)U_{1y}}{V}-\f{\mu(\bar\T)\bar
U_{1y}}{\bar V}\Big)dyd\tau\\
&\di =\int_{\f h\v}^{\tau}\int_{\mathbf{R}}\Big(\f{4
V\Psi_1}{3Z}-\f{8
U_1W}{9Z^2}\Big)\Big[\f{\mu(\T)}{V}b_{2y}+\Big(\f{\mu(\T)}{V}-\f{\mu(\bar\T)}{\bar
V}\Big)\bar U_{1y}\Big]dyd\tau\\
&\di =-\int_{\f h\v}^{\tau}\int_{\mathbf{R}}\Big[\Big(\f{4
V\Psi_1}{3Z}-\f{8
U_1W}{9Z^2}\Big)\f{\mu(\T)}{V}\Big]_yb_{2}dyd\tau\\
&\di\qquad +\int_{\f h\v}^{\tau}\int_{\mathbf{R}}\Big(\f{4
V\Psi_1}{3Z}-\f{8
U_1W}{9Z^2}\Big)\Big(\f{\mu(\T)}{V}-\f{\mu(\bar\T)}{\bar V}\Big)\bar
U_{1y}\Big]dyd\tau:=I_{51}^1+I_{51}^2.
\end{array}
\end{equation}
We can further have
\begin{equation}
\begin{array}{ll}
 \di I_{51}^1&\di =O(1)\int_{\f
 h\v}^{\tau}\int_{\mathbf{R}}|b_2|\Big[|(\Psi_{1y},W_y)|+|(V_y,U_{1y},\T_y,Z_y)||(\Psi_1,W)|\Big]dyd\tau\\
&\di \leq \b\int_{\f
 h\v}^{\tau}\Big[\|(\Psi_{1y},W_y)\|_{L^2(dy)}^2+\|\sqrt{|U^{S_3}_{1y}|}(\Psi_1,W)\|^2_{L^2(dy)}\Big]d\tau\\
 &\di \qquad
 +C_{h,T}~\v\int_{\f h\v}^{\tau}\|(\Psi_1,W)\|^2_{L^2(dy)}d\tau+C_\b\int_{\f
 h\v}^{\tau}\|(b_1,b_2,b_3)\|^2_{L^2(dy)}d\tau\\
 &\di \leq \b\int_{\f
 h\v}^{\tau}\Big[\|(\Psi_{1y},W_y)\|_{L^2(dy)}^2+\|\sqrt{|U^{S_3}_{1y}|}(\Psi_1,W)\|^2_{L^2(dy)}\Big]d\tau\\
 &\di \qquad
 +C_{h,T}~\v\int_{\f
 h\v}^{\tau}\|(\Psi_1,W)\|^2_{L^2(dy)}d\tau+C_{h,T,\b}~\v^{\f12}.
\end{array}
 \label{I511}
\end{equation}
Similarly, we have
\begin{equation}
\begin{array}{ll}
 \di I_{51}^2 =O(1)\int_{\f h\v}^{\tau}\int_{\mathbf{R}}|(\Psi_1,W)||\bar
 U_{1y}||(b_1,b_2,b_3)|dyd\tau\\
 \di \leq \b\int_{\f
 h\v}^{\tau}\|\sqrt{|U^{S_3}_{1y}|}(\Psi_1,W)\|^2d\tau
 +C_{h,T}~\v\int_{\f
 h\v}^{\tau}\|(\Psi_1,W)\|^2d\tau+C_{h,T,\b}~\v^{\f12}.
\end{array}
 \label{I512}
\end{equation}
By integrating \eqref{le1} with respect to $y$ and $\tau$, then
combining all the above estimates, and choosing $\b, \d^{S_3},\v$
and $\chi$ small enough, we have
\begin{equation}
\begin{array}{ll}
\di \|(\Phi,\Psi,W)(\tau,\cdot)\|^2+\int_{\f
h\v}^\tau\Big[\|\sqrt{|U^{S_3}_{1y}|}(\Psi,W)\|^2+\|(\Psi_y,W_y)\|^2\Big]d\tau\\
\di \leq C(\b+\chi)\int_{\f h\v}^\tau
\|\Phi_y\|^2d\tau+C\chi\int_{\f h\v}^\tau
\|(\zeta,\Psi_{yy},\zeta_y)\|^2d\tau\\
 \di\quad +C_{h,T}~\v\int_{\f
 h\v}^{\tau}\|(\Psi,W)\|^2d\tau+C\int_{\f h\v}^\tau\int
 K_1dyd\tau+C_{h,T,\b}~\v^{\f25}.
\end{array}\label{le2}
\end{equation}
Now we turn to  estimate the microscopic term $\di \int_{\f
h\v}^\tau\int K_1dyd\tau$ in \eqref{le2}. Here, we will only
estimate $\di K_{11}=:-\int_{\f h\v}^\tau\int
\f{V}{Z}\Psi_1\int\xi_1^2(\Pi_1-\Pi_{11}^{CD}-\Pi_1^{S_3})d\xi
dyd\tau$ because the other terms in $\di \int_{\f h\v}^\tau\int
K_1dyd\tau$ can be estimated similarly. Let $\mb{M}_\star$ be a
global Maxwellian with the state $(v_\star,u_\star,\t_\star)$
satisfying $\f12\t<\t_\star<\t$ and
$|v-v_\star|+|u-u_\star|+|\t-\t_\star|\le \eta_0$ so that Lemma
\ref{Lemma 4.2} holds, and Lemmas \ref{Lemma-shock} and
\ref{Lemma-shock-L} hold with $\mb{M}_0$ being replaced by
$\mb{M}_\star$. Note that the above choice of the global Maxwellian
$\mb{M}_\star$ can be obtained if the total wave strength $\d$ is
small enough. By the definition of $\Pi_1,\Pi^{CD}_{11}$ and
$\Pi_1^{S_3}$ given in \eqref{(1.21)}, \eqref{Pi-CD-1} and
\eqref{Pi-S-1} respectively, we have
\begin{equation}
\begin{array}{ll}
\di \Pi_1-\Pi^{CD}_{11}-\Pi_1^{S_3} \di
\di
=&\di {\mb{L}}_\mb{M}^{-1}\big[\wt{\mb{G}}_\tau-\f{u_1}{v}\wt{\mb{G}}_{1y}+\f{1}{v}\mb{P}_1(\x_1\wt{\mb{G}}_{1y})-Q(\wt{\mb{G}}_1,\wt{\mb{G}}_1)\big]\\
&\di
-{\mb{L}}_\mb{M}^{-1}[2Q(\wt{\mb{G}}_1,\mb{G}^{R_1}+\mb{G}^{CD}+\mb{G}^{S_3})]+J_3+J_4+J_5.
\end{array}\label{Pi-1}
\end{equation}
Here \begin{equation}
\bar{\mb{G}}^{CD}(\tau,y,\xi)=\f{3}{2V^{CD}\T^{CD}}\mb{L}^{-1}_{\mb{M}^{CD}}\big\{\mb{P}^{CD}_1[\xi_1(\f{|\xi-U^{CD}|^2}{2\T^{CD}}\T^{CD}_y+\xi\cdot{U}^{CD}_{y})\mb{M}^{CD}]\big\},
\label{bar-G-CD}
\end{equation}
\begin{equation}
\begin{array}{ll}
J_3=&\di
({\mb{L}}_\mb{M}^{-1}-{\mb{L}}_\mb{M^{S_3}}^{-1})\big[\mb{G}^{S_3}_\tau-Q(\mb{G}^{S_3},\mb{G}^{S_3})\big]
-(\f{u_1}{v}{\mb{L}}_\mb{M}^{-1}-\f{U_1^{S_3}}{V^{S_3}}{\mb{L}}_\mb{M^{S_3}}^{-1})\mb{G}^{S_3}_y
\\
&\di ~~
+(\f1v{\mb{L}}_\mb{M}^{-1}\mb{P}_1-\f{1}{V^{S_3}}{\mb{L}}_\mb{M^S}^{-1}\mb{P}^{S_3}_1)
(\x_1\mb{G}^{S_3}_y)-{\mb{L}}_\mb{M}^{-1}[2Q(\mb{G}^{S_3},\mb{G}^{R_1}+\mb{G}^{CD})],
\end{array}\label{J3}
\end{equation}
\begin{equation}
\begin{array}{ll}
J_4=&\di
-\big(\f{u_1}{v}{\mb{L}}_\mb{M}^{-1}\mb{G}^{CD}_y-\f{U^{CD}_1}{V^{CD}}{\mb{L}}_{\mb{M}^{CD}}^{-1}\bar{\mb{G}}^{CD}_y\big)-\big[
{\mb{L}}_\mb{M}^{-1}Q(\mb{G}^{CD},\mb{G}^{CD})-{\mb{L}}_{\mb{M}^{CD}}^{-1}Q(\bar{\mb{G}}^{CD},\bar{\mb{G}}^{CD})\big]
\\
&\di
~~+\big[\f1v{\mb{L}}_\mb{M}^{-1}\mb{P}_1(\x_1\mb{G}^{CD}_y)-\f{1}{V^{CD}}{\mb{L}}_{\mb{M}^{CD}}^{-1}\mb{P}^{CD}_1(\x_1\bar{\mb{G}}^{CD}_y)\big]
-{\mb{L}}_\mb{M}^{-1}[2Q(\mb{G}^{R_1},\mb{G}^{CD})],
\end{array}\label{J4}
\end{equation}
and
\begin{equation}
J_5={\mb{L}}_\mb{M}^{-1}\big[-\f{u_1}{v}\mb{G}^{R_1}_y
+\f1v{\mb{L}}_\mb{M}^{-1}\mb{P}_1(\x_1\mb{G}^{R_1}_y)-Q(\mb{G}^{R_1},\mb{G}^{R_1})\big].
\label{J5}
\end{equation}
By \eqref{Pi-1}, we have
\begin{equation}
\begin{array}{ll}
\di K_{11}= -\int_{\f h\v}^\tau\int
\f{V}{Z}\Psi_1\int\xi_1^2\mb{L}_{\mb{M}}^{-1}(\wt{\mb{G}}_{\tau})d\xi
dyd\tau+\int_{\f h\v}^\tau\int\f{u_1
V\Psi_1}{vZ}\int\xi_1^2\mb{L}_{\mb{M}}^{-1}(\wt{\mb{G}}_{1y})d\xi
dyd\tau\\
\di\quad-\int_{\f h\v}^\tau\int\f{
V\Psi_1}{vZ}\int\xi_1^2\mb{L}_\mb{M}^{-1}[\mb{P}_1(\xi_1\wt{\mb{G}}_{1y})]d\xi
dyd\tau+\int_{\f h\v}^\tau\int
\f{V\Psi_1}{Z}\int\xi_1^2\mb{L}_\mb{M}^{-1}[Q(\wt{\mb{G}}_1,\wt{\mb{G}}_1)]d\xi
dyd\tau\\
\di\quad+\int_{\f h\v}^\tau\int\f{
V\Psi_1}{vZ}\int\xi_1^2{\mb{L}}_\mb{M}^{-1}[2Q(\wt{\mb{G}}_1,\mb{G}^{R_1}+\mb{G}^{CD}+\mb{G}^{S_3})]d\xi
dyd\tau-\int_{\f h\v}^\tau\int \f{V\Psi_1}{Z}\int\xi_1^2J_3d\xi
dyd\tau\\
\di\quad-\int_{\f h\v}^\tau\int \f{V\Psi_1}{Z}\int\xi_1^2J_4d\xi
dyd\tau-\int_{\f h\v}^\tau\int \f{V\Psi_1}{Z}\int\xi_1^2J_5d\xi
dyd\tau
 :=
\sum_{i=1}^8K_{11}^{i}.
\end{array}
\label{K11}
\end{equation}
For the integral $K_{11}^{1}$, we have
\begin{equation}
\begin{array}{ll}
K_{11}^{1}&\di =-\int_{\f h\v}^\tau\int
\f{V}{Z}\Psi_1\int\xi_1^2\mb{L}_\mb{M}^{-1}(\wt{\mb{G}}_{1\tau})d\xi
dyd\tau\\
&\di\qquad-\int_{\f h\v}^\tau\int
\f{V}{Z}\Psi_1\int\xi_1^2\mb{L}_\mb{M}^{-1}(\mb{G}^{R_1}_{\tau}+\mb{G}^{CD}_{\tau})d\xi
dyd\tau =:K_{11}^{11}+K_{11}^{12}.
\end{array}
\label{K111}
\end{equation}
Note that the linearized operator $\mb{L}_\mb{M}^{-1}$ satisfies,
for any $h\in \mathfrak{N}^\bot$,
\begin{equation}
\begin{array}{l}
(\mb{L}_\mb{M}^{-1}g)_{\varsigma}=\mb{L}_\mb{M}^{-1}(g_{\varsigma})-2\mb{L}_\mb{M}^{-1}\{Q(\mb{L}_\mb{M}^{-1}g,M_{\varsigma})\},
\quad{\rm for}~~ \varsigma=\tau,y.
\end{array}
\label{facts1}
\end{equation}
Then we have
\begin{equation}
\begin{array}{l}
\di K_{11}^{11}=-\int_{\f h\v}^\tau\int
\f{V}{Z}\Psi_1\int\xi_1^2(\mb{L}_\mb{M}^{-1}\wt{\mb{G}}_1)_{\tau}d\xi
dyd\tau-2\int_{\f h\v}^\tau\int
\f{V}{Z}\Psi_1\int\xi_1^2\mb{L}_\mb{M}^{-1}\{Q(\mb{L}_\mb{M}^{-1}\wt{\mb{G}}_1,M_{\tau})\}d\xi
dyd\tau\\
\quad~ \di=-\int\int
\Big[\f{V}{Z}\Psi_1\xi_1^2\mb{L}_\mb{M}^{-1}(\wt{\mb{G}}_1)\Big](\tau,y,\x)d\xi
dy+\int_{\f h\v}^\tau\int(
\f{V}{Z}\Psi_1)_{\tau}\int\xi_1^2\mb{L}_\mb{M}^{-1}(\wt{\mb{G}}_1)d\xi
dyd\tau\\
\qquad\qquad \di-2\int_{\f h\v}^\tau\int
\f{V}{Z}\Psi_1\int\xi_1^2\mb{L}_\mb{M}^{-1}\{Q(\mb{L}_\mb{M}^{-1}\wt{\mb{G}}_1,\mb{M}_{\tau})\}d\xi
dyd\tau.
\end{array}
\label{K1111}
\end{equation}
The H\"{o}lder inequality and Lemma \ref{Lemma 4.3} yield
\begin{equation}
|\int\xi_1^2\mb{L}_\mb{M}^{-1}(\wt{\mb{G}}_1)d\xi|^2\leq
C\int\f{|\wt{\mb{G}}_1|^2}{\mb{M}_\star}d\xi. \label{K1111-1}
\end{equation}
Moreover, from Lemmas \ref{Lemma 4.1}-\ref{Lemma 4.3}, we have
\begin{equation}
\begin{array}{l}
\di
\int\xi_1^2\mb{L}_\mb{M}^{-1}\{Q(\mb{L}_\mb{M}^{-1}\wt{\mb{G}}_1,\mb{M}_{\tau})\}d\xi\leq
C\left(\int\f{\nu(|\xi|)}{\mb{M}_\star}|\mb{L}_\mb{M}^{-1}\{Q(\mb{L}_\mb{M}^{-1}\wt{\mb{G}}_1,\mb{M}_{\tau})\}|^2d\xi\right)^{\f{1}{2}}
\\ \quad \di
\leq
C\left(\int\f{\nu(|\xi|)}{\mb{M}_\star}|\mb{L}_\mb{M}^{-1}\wt{\mb{G}}_1|^2d\xi\right)^{\f{1}{2}}\cdot\left(\int\f{\nu(|\xi|)}{\mb{M}_\star}|\mb{M}_{\tau}|^2d\xi\right)^{\f{1}{2}}
\\
\quad \di\leq
C|(v_{\tau},u_{\tau},\t_{\tau})|\left(\int\f{\nu^{-1}(|\xi|)}{\mb{M}_\star}|\wt{\mb{G}}_1|^2d\xi\right)^{\f{1}{2}}.
\end{array}
\label{K1111-2}
\end{equation}
Combining (\ref{K1111})-(\ref{K1111-2}) gives
\begin{equation}
\begin{array}{ll}
\di K_{11}^{11}\leq\b\Big[\|\Psi(\tau,\cdot)\|^2+\int_{\f
h\v}^\tau\|(\Psi_{1\tau},\sqrt{|U^{S_3}_{1y}|}\Psi_1)\|^2d\tau\Big]\\
\qquad\di+C_\b
\int\int\f{|\wt{\mb{G}}_1|^2}{\mb{M}_\star}(\tau,y,\x)d\xi
dy+C_\b\int_{\f
h\v}^\tau\int\int\f{\nu(|\xi|)}{\mb{M}_\star}|\wt{\mb{G}}_1|^2d\xi
dyd\tau\\
\qquad\di +C_{h,T}~\v\int_{\f
h\v}^\tau\|\Psi_1\|^2d\tau+C_{h,T}(\v^{\f14}+\chi) \int_{\f
h\v}^\tau\|(\p_{\tau},\psi_{\tau},\z_{\tau})\|^2d\tau.
\end{array}
\label{K1111-e}
\end{equation}
On the other hand, by (\ref{G-R1}) and \eqref{G-CD}, we have
\begin{equation}
\begin{array}{l}
\di K_{11}^{12}=-\int_{\f h\v}^\tau\int
\f{V}{Z}\Psi_1\int\xi_1^2\mb{L}_\mb{M}^{-1}(\mb{G}^{R_1}_{\tau}+\mb{G}^{CD}_{\tau})d\xi
dyd\tau\\
\quad\di \leq C\int_{\f h\v}^\tau\int
|\Psi_1|\Big[|(\T^{R_1}_{y\tau},U^{R_1}_{1y\tau},\T^{CD}_{y\tau},U^{CD}_{y\tau})|
+|(\T^{R_1}_y,U^{R_1}_{1y},\T^{CD}_y,U^{CD}_{y})||(v_{\tau},u_{\tau},\t_{\tau})|\Big]dyd\tau\\
\quad \di\leq C_{h,T}~\v\int_{\f h\v}^\tau\|\Psi_1\|^2d\tau
+C(C_{h,T}~\v+\d)\int_{\f
h\v}^\tau\|(\phi_\tau,\psi_\tau,\zeta_\tau)\|^2d\tau+C_{h,T}~\v^{\f12},
\end{array}
\label{K1112}
\end{equation}
which, together with (\ref{K1111-e}), imply
\begin{equation}
\begin{array}{l}
\di K_{11}^{1}\leq\b\Big[\|\Psi(\tau,\cdot)\|^2+\int_{\f
h\v}^\tau\|(\Psi_{1\tau},\sqrt{|U^{S_3}_{1y}|}\Psi_1)\|^2d\tau\Big]+C_{h,T}~\v^{\f12}\\
\qquad\di+C_\b
\int\int\f{|\wt{\mb{G}}_1|^2}{\mb{M}_\star}(\tau,y,\x)d\xi
dy+C_\b\int_{\f
h\v}^\tau\int\int\f{\nu(|\xi|)}{\mb{M}_\star}|\wt{\mb{G}}_1|^2d\xi
dyd\tau\\
\qquad\di+C_{h,T}~\v\int_{\f
h\v}^\tau\|\Psi_1\|^2d\tau+C[\d+C_{h,T}(\v^{\f14}+\chi)] \int_{\f
h\v}^\tau\|(\p_{\tau},\psi_{\tau},\z_{\tau})\|^2d\tau.
\end{array}
\label{K111-e}
\end{equation}
We now turn to $K_{11}^{2}$. By \eqref{facts1}, we have
\begin{equation}
\begin{array}{ll}
\di K_{11}^{2}&\di= -\int_{\f
h\v}^\tau\int(\f{u_1V\Psi_1}{vZ})_y\int\x_1^2\mb{L}_{\mb{M}}^{-1}(\wt{\mb{G}}_1)
d\x dyd\tau+2\int_{\f
h\v}^\tau\int\f{u_1V\Psi_1}{vZ}\int\x_1^2\mb{L}_{\mb{M}}^{-1}\{Q(\mb{L}_{\mb{M}}^{-1}\wt{\mb{G}}_1,\mb{M}_y)\}
d\x dyd\tau\\
&\di \leq C\int_{\f
h\v}^\tau\int\Big[|\Psi_{1y}|+|\Psi_1||(V_y,Z_y)|+|\Psi_1||(v_y,u_y,\t_y)|\Big]\Big(\int\f{\nu^{-1}(|\xi|)}{\mb{M}_\star}|\wt{\mb{G}}_1|^2d\xi\Big)^{\f{1}{2}}dyd\tau\\
&\di \leq \b\int_{\f
h\v}^\tau\|(\Psi_{1y},\sqrt{|U^{S_3}_{1y}|}\Psi_1)\|^2d\tau+C_\b\int_{\f
h\v}^\tau\int\int\f{\nu(|\xi|)}{\mb{M}_\star}|\wt{\mb{G}}_1|^2d\xi
dyd\tau\\
&\di~~+C_{h,T}~\v\int_{\f
h\v}^\tau\|\Psi_1\|^2d\tau+C_{h,T}(\v^{\f14}+\chi) \int_{\f
h\v}^\tau\|(\p_{y},\psi_{y},\z_{y})\|^2d\tau.
\end{array}
\label{K112}
\end{equation}
To estimate  $K_{11}^{3}$, notice that
\begin{equation}
\mb{P}_1(\xi_1\wt{\mb{G}}_{1y})=[\mb{P}_1(\xi_1\wt{\mb{G}}_1)]_y+\sum_{j=0}^4<\xi_1\wt{\mb{G}}_1,\chi_j>\mb{P}_1(\chi_{jy}).
\label{P1-y}
\end{equation}
Then from \eqref{facts1} and \eqref{P1-y}, we have
\begin{equation}
\begin{array}{ll}
\di K_{11}^{3}&\di=\int_{\f h\v}^\tau\int(\f{
V\Psi_1}{vZ})_y\int\xi_1^2\mb{L}_\mb{M}^{-1}[\mb{P}_1(\xi_1\wt{\mb{G}}_1)]d\xi
dyd\tau\\
 &\di\quad-\int_{\f h\v}^\tau\int\f{
V\Psi_1}{vZ}\int\xi_1^2\mb{L}_\mb{M}^{-1}
[\sum_{j=0}^4<\xi_1\wt{\mb{G}}_1,\chi_j>\mb{P}_1(\chi_{jy})]d\xi
dyd\tau
\\
&\di\quad-2\int_{\f h\v}^\tau\int\f{
V\Psi_1}{vZ}\int\xi_1^2\mb{L}_\mb{M}^{-1}
\{Q(\mb{L}_\mb{M}^{-1}[\mb{P}_1(\xi_1\wt{\mb{G}}_1)],\mb{M}_y)\}d\xi
dyd\tau\\
 &\di\leq \b\int_{\f h\v}^\tau\int
(|U^{S_3}_{1y}||\Psi_1|^2+|\Psi_{1y}|^2)dyd\tau+C_\b\int_{\f
h\v}^\tau\int\int\f{\nu(|\xi|)}{\mb{M}_\star}|\wt{\mb{G}}_1|^2d\xi
dyd\tau\\
&\di \qquad+C_{h,T}~\v\int_{\f
h\v}^\tau\int|\Psi_1|^2dyd\tau+C_{h,T}(\v^{\f14}+\chi)\int_{\f
h\v}^\tau\int|(\phi_y,\psi_y,\zeta_y)|^2 dyd\tau.
\end{array}
\label{K113}
\end{equation}
On the other hand,
\begin{equation}
\begin{array}{l}
\di K_{11}^{4}\leq C\int_{\f
h\v}^\tau\int|\Psi_1|\Big(\int\f{\nu(|\xi|)}{\mb{M}_\star}|\mb{L}_\mb{M}^{-1}\{Q(\wt{\mb{G}}_1,\wt{\mb{G}}_1)\}|^2d\xi\Big)^{\f12}dyd\tau
\\
\quad\di \leq C\int_{\f
h\v}^\tau\int|\Psi_1|\int\f{\nu(|\xi|)}{\mb{M}_\star}|\wt{\mb{G}}_1|^2d\xi
dyd\tau \leq C\chi\int_{\f
h\v}^\tau\int\int\f{\nu(|\xi|)}{\mb{M}_\star}|\wt{\mb{G}}_1|^2d\xi
dyd\tau,
\end{array}
\label{K114}
\end{equation}
and
\begin{equation}
\begin{array}{l}
\di K_{11}^{5}
\di\leq C\int_{\f
h\v}^\tau\int|\Psi_1|\Big(\int\f{\nu(|\xi|)}{\mb{M}_\star}|\wt{\mb{G}}_1|^2d\xi\Big)^{\f12}
\Big[\big(\int\f{\nu(|\xi|)}{\mb{M}_\star}|\mb{G}^{S_3}|^2d\xi\big)^{\f12}+|(\T^{R_1}_y,\T^{CD}_y,U^{R_1}_{1y},U^{CD}_y)|\Big]dyd\tau
\\
\qquad\di \leq \b\int_{\f
h\v}^\tau\int|U^{S_3}_{1y}||\Psi_1|^2dyd\tau+C_\b\int_{\f
h\v}^\tau\int\int\f{\nu(|\xi|)}{\mb{M}_\star}|\wt{\mb{G}}_1|^2d\xi
dyd\tau+C_{h,T}~\v\int_{\f h\v}^\tau\int|\Psi_1|^2dyd\tau.
\end{array}
\label{K115}
\end{equation}
By \eqref{J3}, we have
\begin{equation}
\begin{array}{ll}
\di K^6_{11} &\di=-\int_{\f h\v}^\tau\int
\f{V\Psi_1}{Z}\int\xi_1^2\Big\{({\mb{L}}_\mb{M}^{-1}-{\mb{L}}_\mb{M^{S_3}}^{-1})\mb{G}^{S_3}_\tau-({\mb{L}}_\mb{M}^{-1}-{\mb{L}}_\mb{M^{S_3}}^{-1})Q(\mb{G}^{S_3},\mb{G}^{S_3})
\\
&\di\qquad~~-(\f{u_1}{v}{\mb{L}}_\mb{M}^{-1}-\f{U_1^{S_3}}{V^{S_3}}{\mb{L}}_\mb{M^{S_3}}^{-1})\mb{G}^{S_3}_y
+(\f1v{\mb{L}}_\mb{M}^{-1}\mb{P}_1-\f{1}{V^{S_3}}{\mb{L}}_\mb{M^S}^{-1}\mb{P}^{S_3}_1)
(\x_1\mb{G}^{S_3}_y)\\
&\di
\qquad~~-{\mb{L}}_\mb{M}^{-1}[2Q(\mb{G}^{S_3},\mb{G}^{R_1}+\mb{G}^{CD})]\Big\}d\xi
dyd\tau:=\sum_{i=1}^5 K_{11}^{6i}.
\end{array}
\end{equation}
For brevity, we will only consider the term $K_{11}^{61}$ because
similar estimates hold for the terms $K_{11}^{6i}~(i=2,3,4,5)$.
Direct calculation yields
\begin{equation}\label{K116}
\begin{array}{ll}
\di K_{11}^{61}
 &\di =-\int_{\f
h\v}^{\tau}\int\f{V\Psi_1}{Z}\int\xi_1^2{\mb{L}}_\mb{M}^{-1}\big[2Q\big(\mb{M}^{S_3}-\mb{M},{\mb{L}}_\mb{M^{S_3}}^{-1}
(\mb{G}^{S_3}_\tau)\big)\big]d\xi\\
&\di \leq C\int_{\f
h\v}^{\tau}\int|\Psi_1|\big(\int\f{\nu(|\x|)|\mb{M}^{S_3}-\mb{M}|^2}{\mb{M}_*}d\x\big)^\f12
\big(\int\f{\nu(|\x|)|{\mb{L}}_\mb{M^{S_3}}^{-1}(\mb{G}^{S_3}_\tau)|^2}{\mb{M}_*}d\x\big)^\f12dyd\tau\\
&\di \leq C\int_{\f
h\v}^{\tau}\int|\Psi_1||(v-V^{S_3},u-U^{S_3},\t-\T^{S_3})|
\big(\int\f{\nu^{-1}(|\x|)|\mb{G}^{S_3}_\tau|^2}{\mb{M}_*}d\x\big)^\f12dyd\tau\\
&\di \leq C\d\int_{\f
h\v}^{\tau}\int|U^{S_3}_{1y}||\Psi_1|^2dyd\tau+C\d\int_{\f
h\v}^{\tau}\int|(\Phi_y,\Psi_y, \zeta)|^2dyd\tau +C_{h,T}~\v\int_{\f
h\v}^{\tau}\int|\Psi_1|^2dyd\tau+C_{h,T}~\v^{\f12}.
\end{array}
\end{equation}
By \eqref{J4}, we have
\begin{equation*}
\begin{array}{ll}
\di K^7_{11} \di=-\int_{\f h\v}^\tau\int
\f{V\Psi_1}{Z}\int\xi_1^2\Big\{-\big(\f{u_1}{v}{\mb{L}}_\mb{M}^{-1}\mb{G}^{CD}_y-\f{U^{CD}_1}{V^{CD}}{\mb{L}}_{\mb{M}^{CD}}^{-1}\bar{\mb{G}}^{CD}_y\big)
\\\di\qquad-\big[
{\mb{L}}_\mb{M}^{-1}Q(\mb{G}^{CD},\mb{G}^{CD})-{\mb{L}}_{\mb{M}^{CD}}^{-1}Q(\bar{\mb{G}}^{CD},\bar{\mb{G}}^{CD})\big]
\\
\di \qquad
+\big[\f1v{\mb{L}}_\mb{M}^{-1}\mb{P}_1(\x_1\mb{G}^{CD}_y)-\f{1}{V^{CD}}{\mb{L}}_{\mb{M}^{CD}}^{-1}\mb{P}^{CD}_1(\x_1\bar{\mb{G}}^{CD}_y)\big]-{\mb{L}}_\mb{M}^{-1}[2Q(\mb{G}^{R_1},\mb{G}^{CD})]\Big\}d\xi
dyd\tau:=\sum_{i=1}^4 K_{11}^{7i}.
\end{array}
\end{equation*}
For illustration, we only consider $K_{11}^{71}$ as follows.
\begin{equation}
\begin{array}{ll}
K_{11}^{71} \di=\int_{\f h\v}^\tau\int
\f{V\Psi_1}{Z}\int\xi_1^2\Big[\big(\f{u_1}{v}-\f{U^{CD}_1}{V^{CD}}\big){\mb{L}}_\mb{M}^{-1}\mb{G}^{CD}_y
+\f{U^{CD}_1}{V^{CD}}\big({\mb{L}}_\mb{M}^{-1}-{\mb{L}}_{\mb{M}^{CD}}^{-1}\big)\mb{G}^{CD}_y\\
\di\qquad\quad
+\f{U^{CD}_1}{V^{CD}}{\mb{L}}_{\mb{M}^{CD}}^{-1}\big(\mb{G}^{CD}_y-\bar{\mb{G}}^{CD}_y\big)\Big]d\xi
dyd\tau\\
\di \leq C\int_{\f h\v}^\tau\int|\Psi_1||(v-V^{CD}, u-U^{CD},
\t-\T^{CD})|
\cdot\Big[|(\T^{CD}_{yy},U^{CD}_{yy})|+|(\T^{CD}_{y},U^{CD}_{y})||(v_y,u_y,\t_y)|\Big] dy d\tau\\
\di \leq C_{h,T}~\v\int_{\f
h\v}^\tau\int|\Psi_1|^2dyd\tau+C\v\int_{\f
h\v}^\tau\|(\Phi_y,\Psi_y,\z)\|^2d\tau+C_{h,T}\chi\int_{\f
h\v}^\tau\|(\phi_y,\psi_y,\z_y)\|^2d\tau+C_{h,T}~\v^{\f12}.
\end{array}
\end{equation}
Then, by \eqref{J5}, we have
\begin{equation}
\begin{array}{ll}
K_{11}^8 &\di =-\int_{\f h\v}^\tau\int
\f{V\Psi_1}{Z}\int\xi_1^2{\mb{L}}_\mb{M}^{-1}\big[-\f{u_1}{v}\mb{G}^{R_1}_y
+\f1v{\mb{L}}_\mb{M}^{-1}\mb{P}_1(\x_1\mb{G}^{R_1}_y)+Q(\mb{G}^{R_1},\mb{G}^{R_1})\big]d\xi
dyd\tau\\
&\di\leq C\int_{\f
h\v}^\tau\int|\Psi_1|\Big[|(v_y,u_y,\t_y)||(\T^{R_1}_y,U^{R_1}_y)|+|(\T^{R_1}_y,U^{R_1}_y)|^2+|(\T^{R_1}_{yy},U^{R_1}_{yy})|\Big]dyd\tau\\
&\di\leq C_{h,T}~\v\int_{\f
h\v}^\tau\int|\Psi_1|^2dyd\tau+C_{h,T}~\v\int_{\f
h\v}^\tau\|(\phi_y,\psi_y,\z_y)\|^2d\tau+C_{h,T}~\v^{\f12}.
\end{array}
\end{equation}
By collecting all the above estimates, we have
\begin{equation}
\begin{array}{l}
\di K_{11}\leq \b\Big[\|\Psi(\tau,\cdot)\|^2+\int_{\f
h\v}^\tau\|(\Psi_{1\tau},\sqrt{|U^{S_3}_{1y}|}\Psi_1)\|^2d\tau\Big]+C_\b
\int\int\f{|\wt{\mb{G}}_1|^2}{\mb{M}_\star}(\tau,y,\x)d\xi
dy\\
\qquad\di+C_\b\int_{\f
h\v}^\tau\int\int\f{\nu(|\xi|)}{\mb{M}_\star}|\wt{\mb{G}}_1|^2d\xi
dyd\tau+C(\chi+\beta+\d)\int_{\f
h\v}^\tau\|(\Phi_y,\Psi_y,\z)\|^2d\tau\\
\qquad\di+C_{h,T}~\v\int_{\f
h\v}^\tau\|\Psi_1\|^2d\tau+C[\d+C_{h,T}\chi]\sum_{|\a^\prime|=1}\int_{\f
h\v}^\tau\|\partial^{\a^\prime}(\p,\psi,\z)\|^2d\tau+C_{h,T}~\v^{\f12}.
\end{array}
\label{K11-E}
\end{equation}

Therefore, we
have
\begin{equation}
\begin{array}{l}
\di \|(\Phi,\Psi,W)(\tau,\cdot)\|^2+\int_{\f
h\v}^\tau\Big[\|\sqrt{|U^{S_3}_{1y}|}(\Psi,W)\|^2+\|(\Psi_y,W_y)\|^2\Big]d\tau\\
\di \leq C\int\int\f{|\wt{\mb{G}}_1|^2}{\mb{M}_\star}(\tau,y,\x)d\xi
dy+ C_{h,T}~\v\int_{\f h\v}^\tau\|(\Psi,W)\|^2d\tau+C \beta \int_{\f
h\v}^\tau\|(\Psi_{\tau},W_{\tau})\|^2d\tau\\
\di +C(\chi+\beta+\d)\int_{\f
h\v}^\tau\|(\Phi_y,\z)\|^2d\tau+C_\b\int_{\f
h\v}^\tau\int\int\f{\nu(|\xi|)}{\mb{M}_\star}|\wt{\mb{G}}_1|^2d\xi
dyd\tau\\
\di+C(\d+C_{h,T}\chi)\sum_{|\a^\prime|=1}\int_{\f
h\v}^\tau\|\partial^{\a^\prime}(\p,\psi,\z)\|^2d\tau+C_{h,T,\b}~\v^{\f25},
\end{array}
\label{LE}
\end{equation}
where we have used the smallness of $\delta$, $\beta$, $\v$ and
$\chi$.

\

\underline{Step 2. Estimation on $\|\Phi_y(\tau,\cdot)\|^2$.}

\

Note that the dissipation term does not contain the term
$\|\Phi_y\|^2$. To complete the lower order energy estimate, we have to
estimate $\Phi_y$. From $\eqref{sys}_2$, we have
\begin{equation}
\begin{array}{ll}
\di \f{4}{3}\f{\mu(\T)}{
V}\Phi_{y\tau}-\Psi_{1\tau}+\f{Z}{V}\Phi_y=\f2{3
V}W_y+\f{2}{3V}U_{y}\cdot\Psi\\
\di\qquad\quad
-\f43\f{\mu^\prime(\T)}{V}U_{1y}(W_y+U_{y}\cdot\Psi)-N_1+\bar
Q_1+Q_1+\int\xi_1^2(\Pi_1-\Pi_{11}^{CD}-\Pi_1^{S_3})d\xi.
\end{array}
\label{Phi}
\end{equation}
Multiplying (\ref{Phi}) by $\Phi_y$ yields
\begin{equation}
\begin{array}{l}
\di(\f{2\mu(\T)}{3 V}\Phi_y^2-\Phi_y\Psi_1)_{\tau}+\f{Z}{
V}\Phi_y^2=(\f{2\mu(\T)}{3
V})_{\tau}\Phi_y^2+\Psi_{1y}^2+\Big[\f{2}{3V}U_{y}\cdot\Psi+\f{2}{3
V}W_y\\
\di
\quad-\f43\f{\mu^\prime(\T)}{V}U_{1y}(W_y+U_{y}\cdot\Psi)-N_1+\bar
Q_1+Q_1+\int\xi_1^2(\Pi_1-\Pi_{11}^{CD}-\Pi_1^{S_3})d\xi\Big]\Phi_y+(\cdots)_y,
\end{array}
\label{Phi-1}
\end{equation}
where we have used the fact
$$
\Phi_y\Psi_{1\tau}=(\Phi_y\Psi_1)_{\tau}
-(\Phi_{\tau}\Psi_1)_y+\Psi_{1y}^2.
$$
Integrating (\ref{Phi-1}) with respect to $y$ and $\tau$ and using
H${\rm\ddot{o}}$lder inequality give
\begin{equation}
\begin{array}{l}
\di
\int\Big[\f{2\mu(\T)}{3 V}\Phi_y^2-\Phi_y\Psi_1\Big](\tau,y)dy+\int_{\f h\v}^\tau\int\f{Z}{2 V}\Phi_y^2dyd\tau\\
\quad\di\leq C_{h,T}~\v\int_{\f h\v}^\tau\|\Psi\|^2d\tau+C\int_{\f
h\v}^\tau\|(\Psi_y,W_y)\|^2d\tau+C\int_{\f
h\v}^\tau\int (N_1^2+\bar Q_1^2+Q_1^2)dyd\tau\\
\qquad\di +C\d\int_{\f
h\v}^\tau\|\sqrt{|U^{S_3}_{1y}|}\Psi\|^2d\tau+C\int_{\f
h\v}^\tau\int|\int\xi_1^2(\Pi_1-\Pi_{11}^{CD}-\Pi_1^{S_3})d\xi|^2dyd\tau.
\end{array}
\label{(4.40)}
\end{equation}
By \eqref{barQ1}, \eqref{N1} and \eqref{Q1-2} and the Cauchy
inequality, one has
\begin{equation}
\begin{array}{ll}
&\di\int_{\f h\v}^\tau\int (N_1^2+\bar Q_1^2+Q_1^2)dyd\tau\leq
C\chi\int_{\f h\v}^\tau\|(\Phi_y,\Psi_y,\zeta)\|^2d\tau+
C_{h,T}\chi\int_{\f h\v}^\tau\|\psi_{1y}\|^2d\tau+C_{h,T}~\v^{\f12}.
\end{array}
 \label{(4.41+)}
\end{equation}
Now we estimate the last term on the right hand side of
\eqref{(4.40)}. By \eqref{Pi-1}, we have
\begin{equation}
\begin{array}{ll}
\di  K_2:=\int_{\f
h\v}^\tau\int|\int\xi_1^2(\Pi_{1}-\Pi^{CD}_{11}-\Pi^{S_3}_1)d\xi|^2dyd\tau\\
 \di \quad\leq C
\int_{\f
h\v}^\tau\int|\int\xi_1^2\mb{L}_\mb{M}^{-1}(\wt{\mb{G}}_{\tau})d\xi|^2dyd\tau+C\int_{\f
h\v}^\tau\int|\int\xi_1^2\mb{L}_\mb{M}^{-1}(\f{u_1}v\wt{\mb{G}}_{1y})d\xi|^2dyd\tau\\
\di~~+C\int_{\f
h\v}^\tau\int|\int\xi_1^2\mb{L}_\mb{M}^{-1}[\f{1}{v}\mb{P}_1(\xi_1\wt{\mb{G}}
_{1y})]d\xi|^2dyd\tau +C\int_{\f
h\v}^\tau\int|\int\xi_1^2\mb{L}_\mb{M}^{-1}[Q(\wt{\mb{G}}_1,\wt{\mb{G}}_1)]d\xi|^2dyd\tau\\
\di~~ +C\int_{\f
h\v}^\tau\int|\int\xi_1^2\mb{L}_\mb{M}^{-1}[2Q(\wt{\mb{G}}_1,\mb{G}^{R_1}+\mb{G}^{CD}+\mb{G}^{S_3})]d\xi|^2dyd\tau\\
\di~~ +C\int_{\f
h\v}^\tau\int\int\xi_1^2(|J_3|^2+|J_4|^2+|J_5|^2)d\xi dyd\tau
:=\sum_{i=1}^6 K_2^i.
\end{array}\label{K2}
\end{equation}
Then we can obtain
\begin{equation}
\begin{array}{ll}
K_2^1&\di \leq C\int_{\f
h\v}^\tau\int|\int\f{\nu(|\x|)|\mb{L}_\mb{M}^{-1}\wt{\mb{G}}_{\tau}|^2}{\mb{M}_\star}
d\x\cdot \int \nu^{-1}(|\x|)\x_1^4\mb{M}_\star d\x| dyd\tau\\
&\di\leq C\int_{\f h\v}^\tau\int\int
\f{\nu^{-1}(|\x|)|\wt{\mb{G}}_{\tau}|^2}{\mb{M}_\star}d\x dyd\tau.
\end{array}\label{(4.31)}
\end{equation}
Similarly,
\begin{equation}
\begin{array}{ll}
K_2^2 &\di \leq C\int_{\f h\v}^\tau\int\int
\f{\nu^{-1}(|\x|)|\wt{\mb{G}}_{1y}|^2}{\mb{M}_\star}d\x dyd\tau\\
&\di \leq C\int_{\f h\v}^\tau\int\int
\f{\nu^{-1}(|\x|)|\wt{\mb{G}}_{y}|^2}{\mb{M}_\star}d\x
dyd\tau+C\int_{\f h\v}^\tau\int\int
\f{\nu^{-1}(|\x|)}{\mb{M}_\star}(|\mb{G}^{R_1}_{y}|^2+|\mb{G}^{CD}_{y}|^2)d\x
dyd\tau\\
&\di \leq C\int_{\f h\v}^\tau\int\int
\f{\nu^{-1}(|\x|)|\wt{\mb{G}}_{y}|^2}{\mb{M}_\star}d\x
dyd\tau+C\int_{\f
h\v}^\tau\|(\T^{R_1}_{yy},U^{R_1}_{1yy},\T^{CD}_{yy},U^{CD}_{yy})\|^2
d\tau\\
&\di\qquad +\int_{\f
h\v}^\tau\int|(v_y,u_y,\t_y)|^2\cdot|(\T^{R_1}_{y},U^{R_1}_{1y},\T^{CD}_{y},U^{CD}_{y})|^2
dyd\tau\\
&\di \leq C\int_{\f h\v}^\tau\int\int
\f{\nu^{-1}(|\x|)|\wt{\mb{G}}_{y}|^2}{\mb{M}_\star}d\x
dyd\tau+C\v\int_{\f
h\v}^\tau\|(\phi_y,\psi_y,\z_y)\|^2d\tau+C_{h,T}~\v^{\f12}.
\end{array} \label{(4.32)}
\end{equation}
Moreover,
\begin{equation}
\begin{array}{ll}
\di K_2^3 \leq C\int_{\f
h\v}^\tau\int|\int\f{\nu(|\x|)|\mb{L}_\mb{M}^{-1}[\f{1}{v}\mb{P}_1(\xi_1\wt{\mb{G}}
_{1y})]|^2}{\mb{M}_{[2v_\star,2u_\star,2\t_\star]}} d\x\cdot \int
\nu^{-1}(|\x|)\x_1^4\mb{M}_{[2v_\star,2u_\star,2\t_\star]}d\x|
dyd\tau\\
\quad\di \leq C\int_{\f h\v}^\tau\int\int
\f{\nu^{-1}(|\x|)|\f{1}{v}\mb{P}_1(\xi_1\wt{\mb{G}}
_{1y})|^2}{\mb{M}_{[2v_\star,2u_\star,2\t_\star]}}d\x dyd\tau \leq
C\int_{\f h\v}^\tau\int\int \f{\nu(|\x|)|\wt{\mb{G}}
_{1y}|^2}{\mb{M}_\star}d\x dyd\tau\\
\quad\di \leq C\int_{\f h\v}^\tau\int\int
\f{\nu(|\x|)|\wt{\mb{G}}_{y}|^2}{\mb{M}_\star}d\x
dyd\tau+C\v\int_{\f
h\v}^\tau\|(\phi_y,\psi_y,\z_y)\|^2d\tau+C_{h,T}~\v^{\f12}.
\end{array}
\label{(4.33)}
\end{equation}
From Lemma \ref{Lemma 4.1}, we have
\begin{equation}
\begin{array}{ll}
K_2^4 &\di \leq C\int_{\f
h\v}^\tau\int\int\f{\nu^{-1}(|\x|)|Q(\wt{\mb{G}}_1 ,\wt{\mb{G}}_1
)|^2}{\mb{M}_\star}d\xi dyd\tau \leq C\int_{\f h\v}^\tau\int
\int\f{\nu(|\x|)|\wt{\mb{G}}_1|^2}{\mb{M}_\star}
d\x\cdot \int \f{|\wt{\mb{G}}_1|^2}{\mb{M}_\star}d\x dyd\tau\\
&\di\leq C_{h,T}\chi\int_{\f
h\v}^\tau\int\int\f{\nu(|\x|)|\wt{\mb{G}}_1|^2}{\mb{M}_\star}d\x
dyd\tau.
\end{array}\label{(4.34)}
\end{equation}
Using similar idea as for  \eqref{K115} and \eqref{K116}, we can obtain
 \begin{equation}
\begin{array}{ll}
K_2^5&\di\leq C\int_{\f
h\v}^\tau\int\int\f{\nu^{-1}(|\x|)|Q(\wt{\mb{G}}_1,\mb{G}^{R_1}+\mb{G}^{CD}+\mb{G}^{S_3})|^2}{\mb{M}_\star}d\xi dyd\tau\\
&\di \leq C\int_{\f
h\v}^\tau\int\int\f{\nu(|\x|)|\wt{\mb{G}}_1|^2}{\mb{M}_\star}d\xi\int\f{\nu(|\x|)(|\mb{G}^{R_1}|^2+|\mb{G}^{CD}|^2+|\mb{G}^{S_3}|^2)}{\mb{M}_\star}d\x dyd\tau\\
&\di \leq C(\d+\v)\int_{\f
h\v}^\tau\int\int\f{\nu(|\x|)|\wt{\mb{G}}_1|^2}{\mb{M}_\star}d\xi
dyd\tau,
\end{array}
\end{equation}
and
\begin{equation}
\begin{array}{ll}
\di K_2^6&\di \leq C\int_{\f
h\v}^\tau\int\int\f{\nu^{-1}(|\x|)(|J_3|^2+|J_4|^2+|J_5|^2)}{\mb{M}_\star}d\xi dyd\tau\\
& \di \leq C(\d+\v)\int_{\f
h\v}^\tau\|(\Phi_y,\Psi_y,\z)\|^2d\tau+C_{h,T}~\v^{\f12}.
\end{array}
\end{equation}

 Substituting the estimations of $K_2^i~(i=1,2,\cdots,6)$ into
\eqref{K2}, we obtain
\begin{equation}
\begin{array}{ll}
\di K_2=\int_{\f
h\v}^\tau\int|\int\xi_1^2(\Pi_1-\Pi_{11}^{CD}-\Pi_1^{S_3})d\xi|^2dyd\tau\\
\di \quad \leq C(\d+\v)\int_{\f
h\v}^\tau\|(\Phi_y,\Psi_y,\z,\phi_y,\psi_y,\z_y)\|^2d\tau+C_{h,T}~\v^{\f12}\\
\quad\di+C(\d+\chi)\int_{\f
h\v}^\tau\int\int\f{\nu(|\xi|)}{\mb{M}_\star}|\wt{\mb{G}}_1|^2d\xi
dyd\tau+C\sum_{|\a^\prime|=1}\int_{\f
h\v}^\tau\int\int\f{\nu(|\xi|)}{\mb{M}_\star}|\partial^{\a^\prime}\wt{\mb{G}}|^2d\xi
dyd\tau.
\end{array}
\label{K2-E}
\end{equation}
Thus combining \eqref{(4.40)}, (\ref{(4.41+)}) and (\ref{K2-E})
yields
\begin{equation}
\begin{array}{l}
\di
\|\Phi_y(\tau,\cdot)\|^2+\int_{\f h\v}^\tau\|\Phi_y\|^2d\tau\leq C\|\Psi_1(\tau,\cdot)\|^2+C\int_{\f h\v}^\tau\|(\Psi_y,W_y,\z)\|^2d\tau\\
\quad\di+C(\delta+\chi)\int_{\f
h\v}^\tau\|(\phi_y,\psi_y,\z_y)\|^2d\tau+C\sum_{|\a^\prime|=1}\int_{\f
h\v}^\tau\int\int\f{\nu(|\xi|)}{\mb{M}_\star}|\partial^{\a^\prime}\wt{\mb{G}}|^2d\xi
dyd\tau\\
\quad\di +C\d\int_{\f
h\v}^\tau\|\sqrt{|U^{S_3}_{1y}|}\Psi\|^2d\tau+C[\d+C_{h,T}\chi]\int_{\f
h\v}^\tau\int\int\f{\nu(|\xi|)}{\mb{M}_\star}|\wt{\mb{G}}_1|^2d\xi
dyd\tau.
\end{array}
\label{Phi-y-E}
\end{equation}
\underline{Step 3. Estimation on the non-fluid component.}

The microscopic component $\wt{\mb{G}}_1$ can be estimated by using
the equation \eqref{G1e}. Multiplying \eqref{G1e} by
$\f{v\wt{\mb{G}}_1}{\mb{M}_\star}$ gives
\begin{eqnarray}\label{M.1}
&&\bigg(\frac{v\wt{\mb{G}}_1^2}{2\mb{M}_\star}\bigg)_{1\tau}-\f{v\wt{\mb{G}}_1}{\mb{M}_\star}\mb{L}_\mb{M}\wt{\mb{G}}_1
=v_{\tau}\frac{\wt{\mb{G}}_1^2}{2\mb{M}_\star}+\bigg\{\f{u_1}{v}\wt{\mb{G}}
_y-\f{1}{v}\mb{P}_1(\xi_1\wt{\mb{G}}_y)+2Q(\wt{\mb{G}},\mb{G}^{S_3})\nonumber\\
&&\qquad\qquad\qquad\qquad\qquad\qquad+Q(\wt{\mb{G}},\wt{\mb{G}})+J_1+J_2-\mb{G}^{R_1}_{\tau}-\mb{G}^{CD}_{\tau}\bigg\}\f{v\wt{\mb{G}}_1}{\mb{M}_\star}
\\
&&\quad\di=
\bigg\{-\mb{P}_1(\xi_1\wt{\mb{G}}_y)+2vQ(\wt{\mb{G}},\mb{G}^{S_3})+vQ(\wt{\mb{G}},\wt{\mb{G}})+vJ_1+vJ_2\nonumber\\
&&\qquad\qquad\qquad
-v\mb{G}^{R_1}_{\tau}-v\mb{G}^{CD}_{\tau}+u_1\mb{G}^{R_1}_y+u_1\mb{G}^{CD}_y\bigg\}
\f{\wt{\mb{G}}_1}{\mb{M}_\star}+(\cdots)_y\nonumber.
\end{eqnarray}
The Cauchy inequality implies
\begin{equation}\label{M.4}
\begin{array}{ll}
\di \int_{\f h\v}^\tau\int\int
\bigg(-v\mb{G}^{R_1}_{\tau}-v\mb{G}^{CD}_{\tau}+u_1\mb{G}^{R_1}_y+u_1\mb{G}^{CD}_y\bigg)\f{\wt{\mb{G}}_1}{\mb{M}_\star}d\xi
dyd\tau\\
\di \leq \f{\wt{\s}}{32}\int_{\f
h\v}^\tau\int\int\f{\nu(|\xi|)\wt{\mb{G}}_1^2}{\mb{M}_\star}d\xi
dyd\tau+C\int_{\f
h\v}^\tau\int\sum_{|\a|=2}|\partial^\a(\T^{CD},U^{CD},\T^{R_1},U^{R_1})|^2dyd\tau\\
\di\quad +C\int_{\f
h\v}^\tau\int|(\T^{CD}_{y},U^{CD}_y,\T^{R_1}_{y},U^{R_1}_y)|^2\Big[|(v_{\tau},u_{\tau},\t_{\tau})|^2+|(v_y,u_y,\t_y)|^2\Big]dyd\tau\\
\di\leq \f{\wt{\s}}{32}\int_{\f
h\v}^\tau\int\int\f{\nu(|\xi|)\wt{\mb{G}}_1^2}{\mb{M}_\star}d\xi
dyd\tau+C(\d+C_{h,T}\chi)\int_{\f h\v}^\tau\sum_{|\a^\prime|=1}
\|\partial^{\a^\prime}(\phi,\psi,\zeta)\|^2d\tau + C_{h,T}\v^{\f12}.
\end{array}
\end{equation}
Notice that
$\mb{P}_1(\xi_1\wt{\mb{G}}_y)=\xi_1\wt{\mb{G}}_y-\sum_{j=0}^{4}\langle\xi_1\wt{\mb{G}}_y,\chi_j\rangle\chi_j$. Then
we have
\begin{eqnarray}\label{M.5}
&&\int_{\f h\v}^\tau\int\int
\mb{P}_1(\xi_1\wt{\mb{G}}_y)\f{\wt{\mb{G}}_1}{\mb{M}_\star}d\xi
dyd\tau\nonumber\\
&&=\int_{\f
h\v}^\tau\int\int\xi_1(\wt{\mb{G}}_{1y}+{\mb{G}}^{R_1}_y+{\mb{G}}^{CD}_y)\f{\wt{\mb{G}}_1}{\mb{M}_\star}
-\sum_{j=0}^{4}\int_{\f
h\v}^\tau\int\int\langle\xi_1\wt{\mb{G}}_y,\chi_j\rangle\chi_j\f{\wt{\mb{G}}_1}{\mb{M}_\star}\nonumber\\
&&\leq \f{\wt\s}{32}\int_{\f
h\v}^\tau\int\int\f{\nu(|\xi|)\wt{\mb{G}}_1^2}{\mb{M}_\star}d\xi
dyd\tau+C\int_{\f
h\v}^\tau\int\int\f{\nu(|\xi|)\wt{\mb{G}}_y^2}{\mb{M}_\star}d\xi
dyd\tau\nonumber\\
&&\quad+C(\d+C_{h,T}\chi)\int_{\f
h\v}^\tau\|\partial_y(\phi,\psi,\zeta)\|^2d\tau+C_{h,T}\v^{\f12},
\end{eqnarray}
where we have used the fact that
\begin{eqnarray}\nonumber
|\langle\xi_1\wt{\mb{G}}_y,\chi_j\rangle|^2\leq
\int\f{\nu(|\xi|)|\wt{\mb{G}}_y|^2}{\mb{M}_\star}d\xi.
\end{eqnarray}
Lemma \ref{Lemma 4.1}-Lemma \ref{Lemma 4.4} and wave interaction
estimates imply
\begin{eqnarray}\label{M.7}
&&\int_{\f h\v}^\tau\int\int
Q(\wt{\mb{G}},\wt{\mb{G}})\f{v\wt{\mb{G}}_1}{\mb{M}_\star}d\xi
dyd\tau\nonumber\\
&&\leq C\int_{\f
h\v}^\tau\int\bigg(\int\f{\nu(|\xi|)\wt{\mb{G}}_1^2}{\mb{M}_\star}d\xi\bigg)^{\f12}
\bigg(\int\f{\nu(|\xi|)\wt{\mb{G}}^2}{\mb{M}_\star}d\xi\int\f{\wt{\mb{G}}^2}{\mb{M}_\star}d\xi\bigg)^{\f12}dyd\tau\nonumber\\
&&\leq \f{\wt\s}{32}\int_{\f
h\v}^\tau\int\int\f{\nu(|\xi|)\wt{\mb{G}}_1^2}{\mb{M}_\star}d\xi
dyd\tau+C\int_{\f
h\v}^\tau\int|(\T^{CD}_{y},U^{CD}_{y},\T^{R_1}_{y},U^{R_1}_{y})|^4
dyd\tau\nonumber\\
&&\leq \f{\wt\s}{32}\int_{\f
h\v}^\tau\int\int\f{\nu(|\xi|)\wt{\mb{G}}_1^2}{\mb{M}_\star}d\xi
dyd\tau+C_{h,T}~\v^{\f12}.
\end{eqnarray}

Moreover, Lemma \ref{LemmaII} and Cauchy's inequality and wave
interaction estimates imply
\begin{eqnarray}\label{M.9}
&&\int_{\f h\v}^\tau\int\int
(J_1+J_2)\f{v\wt{\mb{G}}_1}{\mb{M}_\star}d\xi dyd\tau \leq
\f{\wt\s}{16}\int_{\f
h\v}^\tau\int\int\f{\nu(|\xi|)\wt{\mb{G}}_1^2}{\mb{M}_\star}d\xi
dyd\tau+C\int_{\f
h\v}^\tau\|(\psi_{y},\zeta_{y})\|^2d\tau\nonumber\\
&&\qquad\qquad\qquad\qquad\qquad\qquad\qquad+C\d\int_{\f
h\v}^\tau\|(\p,\psi,\zeta)\|^2d\tau+C_{h,T}~\v^{\f12}.
\end{eqnarray}

Integrating \eqref{M.1} with respect to $\xi,y$ and $\tau$ and using
\eqref{(4.13)}, \eqref{M.4}-\eqref{M.9} and the smallness of
$\chi,\d,\v$ yield that
\begin{equation}\label{M.10}
\begin{array}{ll}
\di \int\int\f{\wt{\mb{G}} _1^2}{\mb{M}_\star}(\tau,y,\xi)d\xi
dy+\int_{\f
h\v}^\tau\int\int\f{\nu(|\xi|)|\wt{\mb{G}} _1|^2}{\mb{M}_\star}d\xi dyd\tau\\
\di \leq C\d\int_{\f h\v}^\tau\|(\p,\psi,\zeta)\|^2d\tau
+C\sum_{|\a^\prime|=1}\int_{\f
h\v}^\tau\|\partial^{\a^\prime}(\p,\psi,\z)\|^2d\tau+C\int_{\f
h\v}^\tau\int\int\f{\nu(|\xi|)}{\mb{M}_\star}|\wt{\mb{G}}_y|^2d\xi
dyd\tau +C_{h,T}~\v^\f12.
\end{array}
\end{equation}

On the other hand, from the fluid-type system  (\ref{sys}), we can
get an estimate for $\|(\Psi_{\tau},W_{\tau})\|^2$ as follows.
\begin{equation}
\begin{array}{ll}
\di \int_{\f h\v}^\tau\|(\Psi_{\tau},W_{\tau})\|^2d\tau\leq C
\int_{\f h\v}^\tau\|(\Phi_y,\Psi_y,W_y,\Psi_{yy},\z_{y})\|^2d\tau
+C\d\int_{\f h\v}^\tau\|\sqrt{|U^{S_3}_{1y}|}(\Psi,W)\|^2d\tau\\
\di\qquad+C\v\int_{\f
h\v}^\tau\|(\Psi,W)\|^2d\tau+C(\delta+C_{h,T}\chi)\int_{\f
h\v}^\tau\int\int\f{\nu(|\xi|) }{\mb{M}_\star}|\wt{\mb{G}}_1|^2d\xi
dyd\tau\\
\di\qquad+C\sum_{|\a^\prime|=1}\int_{\f
h\v}^\tau\int\int\f{\nu(|\xi|)}{\mb{M}_\star}|\partial^{\a^\prime}\wt{\mb{G}}|^2d\xi
dyd\tau +C_{h,T}~\v^{\f12}.
\end{array}
\label{(4.46)}
\end{equation}
From \eqref{zeta}, we have
\begin{equation}
\begin{array}{ll}
\di\int_{\f h\v}^\tau\|\z\|^2d\tau \leq C\int_{\f
h\v}^\tau\|W_y\|^2d\tau+C\int_{\f
h\v}^\tau\int|\Psi_y|^4dyd\tau+C\int_{\f
h\v}^\tau\int|U_y\cdot\Psi|^2dy d\tau\\
\qquad\di \leq C\int_{\f h\v}^\tau\|W_y\|^2d\tau+C\chi^2\int_{\f
h\v}^\tau\|\Psi_y\|^2d\tau+C\d\int_{\f
h\v}^\tau\|\sqrt{|U^{S_3}_{1y}|}\Psi\|^2d\tau+C_{h,T}~\v\int_{\f
h\v}^\tau\|\Psi\|^2d\tau.
\end{array}\label{Zeta-E}
\end{equation}
In summary, collecting the estimates \eqref{LE}, \eqref{Phi-y-E},
\eqref{M.10}-\eqref{Zeta-E}, we complete the proof of Proposition
\ref{Prop3.1}.

\subsection{Proof of Proposition \ref{Prop3.2}}
\setcounter{equation}{0}

\underline{Proof of Proposition 3.2.}
The proof is divided into the following five steps.
\

\underline{Step 1. Estimation on
$\|(\phi,\psi,\z)(\tau,\cdot)\|^2$.}

Similar to \eqref{le1}, we multiply $\eqref{sys-h}_1$ by $\di\phi$,
$\eqref{sys-h}_2$ by $\di \f{V}{Z}\psi_1$, $\eqref{sys-h}_3$ by
$\di\psi_i$, $\eqref{sys-h}_4$ by $\di \f{2\z}{3Z^2}$ respectively
and adding them together to have
\begin{equation}
\begin{array}{l}
\di\left(\f{\phi^2}{2}+\f{V}{2Z}\psi_1^2+\sum_{i=2}^3\f{\psi_i^2}{2}+
\f{\z^2}{3Z^2}\right)_{\tau} +\f{4\mu(\T)}{3Z}\psi_{1y}^2
+\sum_{i=2}^3\f{\mu(\T)}{V}\psi_{iy}^2+\f{2\k(\T)}{3Z^2V}\z_y^2\\
\di=I_6(\phi,\psi,\z,\psi_y,\z_y)+\f{V}{Z}\psi_1(N_5-\bar
Q_{1y}-Q_{1y})+\sum_{i=2}^3\psi_i(N_{i+4}-\bar Q_{iy}-Q_{iy})\\
\di \quad+\f{2\z}{3Z^2}\big[N_8-\bar
Q_{4y}-Q_{4y}+\sum_{i=1}^3U_i(\bar
Q_{iy}+Q_{iy})\big]+K_3+(\cdots)_y,
\end{array}
\label{He1}
\end{equation}
where
\begin{equation}
\begin{array}{ll}
\di
I_6(\phi,\psi,\z,\psi_y,\z_y)=&\di (\f{V}{2Z})_\tau\psi_1^2+(\f{1}{3Z^2})_\tau\z^2+(\f{2}{3Z})_y\psi_1\z-\f{V}{Z}\psi_1H_1-\sum_{i=2}^3\psi_iH_i\\
&\di -\f{2}{3Z^2}\z
H_4-\f{4\mu(\T)}{3V}(\f{V}{Z})_y\psi_1\psi_{1y}-\f{\k(\T)}{V}(\f{2}{3Z^2})_y\z\z_y,
\end{array}
 \label{I6}
\end{equation}
and
\begin{equation}
\begin{array}{ll}
K_3&\di =-\f{V}{Z}\psi_1\int\x_1^2(\Pi_{1}-\Pi^{CD}_{11}-\Pi^{S_3}_1)_yd\x-\sum_{i=2}^3\psi_i\int\x_1\x_i(\Pi_{1}-\Pi^{CD}_{11}-\Pi^{S_3}_1)_yd\x\\
&\di\quad + \f{2\z}{3Z^2}\Big[-\f12\int\xi_1|\xi|^2
(\Pi_{1}-\Pi^{CD}_{11}-\Pi^{S_3}_1)_yd\xi+\sum_{i=1}^3\psi_i\int\xi_1\xi_i\Pi_{1y}d\x\\
&\di\qquad~~~~~~+\sum_{i=1}^3U_i\int\xi_1\xi_i(\Pi_{1}-\Pi^{CD}_{11}-\Pi^{S_3}_1)_yd\xi\Big].
\end{array}
 \label{K3}
\end{equation}
Direct calculation shows that
\begin{equation}
\begin{array}{ll}
I_6&\di \leq
C|(V_\tau,Z_\tau,V_y,U_y,\T_y,Z_y)|\big[|(\phi_y,\psi_y,\z_y)|^2+(\phi,\psi,\z)|^2\big]\\
&\di \leq
C(\d+C_{h,T}\v^{\f12})|(\phi_y,\psi_y,\z_y)|^2+C(\d+C_{h,T}\v^{\f12})|(\phi,\psi,\z)|^2.
\end{array} \label{I6}
\end{equation}
Thus, integrating \eqref{He1} with respect to $\tau $ and $y$ and
using Cauchy inequality yield that
\begin{equation}
\begin{array}{ll}
&\di \|(\phi,\psi,\zeta)(\tau,\cdot)\|^2+\int_{\f
h\v}^\tau\|(\psi_y,\zeta_y)\|^2d\tau \leq
C(\d+C_{h,T}\v^{\f12})\int_{\f
h\v}^\tau\|\phi_y\|^2d\tau\\
&\di\quad+ C(\b+\d+C_{h,T}~\v^{\f12})\int_{\f
h\v}^\tau\|(\phi,\psi,\zeta)\|^2d\tau+C_{\b}\int_{\f
h\v}^\tau\int\sum_{i=1}^4(\bar Q_{iy}^2+Q_{iy}^2)
dyd\tau\\
&\di \quad +C\int_{\f
h\v}^\tau\int|(\psi,\z)|\big[|N_5|+\sum_{i=2}^3|N_{i+4}|+|N_8|\big]
dyd\tau+\int_{\f h\v}^\tau\int K_3 dyd\tau.
\end{array}
 \label{He2}
\end{equation}
By \eqref{barQ1-E} and \eqref{Q1-E}, we have
\begin{equation}
\int_{\f h\v}^\tau\int\sum_{i=1}^4(\bar Q_{iy}^2+Q_{iy}^2)
dyd\tau\leq C_{h,T}~\v^{\f12}.\label{bar-Q-E}
\end{equation}
From \eqref{N5}, \eqref{Ni+4} and \eqref{N8}, we have
\begin{equation}
\begin{array}{ll}
\di \int_{\f
h\v}^\tau\int|(\psi,\z)|\big[|N_5|+\sum_{i=2}^3|N_{i+4}|+|N_8|\big]
dyd\tau\leq C_{h,T}\chi\int_{\f
h\v}^\tau\|(\phi,\psi,\z,\phi_y,\psi_y,\z_y,\psi_{yy},\z_{yy})\|^2d\tau.
\end{array}\label{N5-8}
\end{equation}
Now we estimate the microscopic term $\di \int_{\f h\v}^\tau\int
K_3dyd\tau$ in \eqref{He2}. We only estimate study $\di
K_{31}:=-\int_{\f h\v}^\tau\int\f{V}{Z}
\psi_1\int\xi_1^2(\Pi_{1}-\Pi^{CD}_{11}-\Pi^{S_3}_1)_yd\xi dyd\tau$
and  $\di
K_{32}:=\f{\z}{3Z^2}\sum_{i=1}^3\psi_i\int\xi_1\xi_i\Pi_{1y}d\x $
because the other terms in $\di \int_{\f h\v}^\tau\int K_3dyd\tau$
can be estimated similarly.

For $K_{31}$, integration by parts with respect to $y$ and the
Cauchy inequality yield
\begin{equation}
\begin{array}{ll}
K_{31}&\di=\int_{\f h\v}^\tau\int(\f{V}{Z}\psi_1)_y\int\xi_1^2(\Pi_{1}-\Pi^{CD}_{11}-\Pi^{S_3}_1)d\xi dyd\tau\\
&\di \leq \b \int_{\f
h\v}^\tau\|\psi_{1y}\|^2d\tau+C(\d+C_{h,T}\v^{\f12})\int_{\f
h\v}^\tau\|\psi_{1}\|^2d\tau+C_\b\int_{\f
h\v}^\tau\int|\int\xi_1^2(\Pi_{1}-\Pi^{CD}_{11}-\Pi^{S_3}_1)d\xi|^2dyd\tau,
\end{array}\label{K31}
\end{equation}
where the last term has been estimated in \eqref{K2-E}.

For $K_{32}$, we have
\begin{equation}
\begin{array}{ll}
K_{32}&\di
 =-\sum_{i=1}^3\int_{\f h\v}^\tau\int(\f{\z\psi_i}{3Z^2})_y\int\xi_1\x_i\Pi_{1}d\xi dyd\tau\\
&\di \leq C\sum_{i=1}^3\int_{\f h\v}^\tau\int\Big(|\z_y||\psi|+|\psi_y||\z|+|Z_y||\psi||\z|\Big)\\
&\di \qquad\qquad\cdot\Big[|\int\xi_1\x_i(\Pi_{1}-\Pi^{CD}_{11}-\Pi^{S_3}_1)d\xi|+ |\int\xi_1\x_i(\Pi^{CD}_{11}+\Pi^{S_3}_1)d\xi|\Big]dyd\tau\\
&\di \leq \b \int_{\f
h\v}^\tau\|(\psi_{y},\z_y)\|^2d\tau+C_\b(\d+C_{h,T}\v^{\f12})\int_{\f
h\v}^\tau\|(\psi,\z)\|^2d\tau\\
&\qquad\di+C_\b\int_{\f
h\v}^\tau\int|\int\xi_1^2(\Pi_{1}-\Pi^{CD}_{11}-\Pi^{S_3}_1)d\xi|^2dyd\tau.
\end{array}\label{K32}
\end{equation}
Substituting \eqref{bar-Q-E}-\eqref{K32} and \eqref{K2-E}  into
\eqref{He2} and choosing $\b, \v,\d,\chi$ suitably small yield that
\begin{equation}
\begin{array}{ll}
&\di \|(\phi,\psi,\zeta)(\tau,\cdot)\|^2+\int_{\f
h\v}^\tau\|(\psi_y,\zeta_y)\|^2d\tau\leq C(\d+C_{h,T}\chi)\int_{\f
h\v}^\tau\|\phi_y\|^2d\tau\\
&\di+ C_\b(\d+C_{h,T}\chi)\int_{\f
h\v}^\tau\|(\phi,\psi,\zeta)\|^2d\tau+C_{h,T}\chi\int_{\f
h\v}^\tau\|(\phi_y,\psi_{yy},\z_{yy})\|^2d\tau+C_{h,T,\b}~\v^{\f12}\\
&\di+C(\d+\chi)\int_{\f
h\v}^\tau\int\int\f{\nu(|\xi|)}{\mb{M}_\star}|\wt{\mb{G}}_1|^2d\xi
dyd\tau+C\sum_{|\a^\prime|=1}\int_{\f
h\v}^\tau\int\int\f{\nu(|\xi|)}{\mb{M}_\star}|\partial^{\a^\prime}\wt{\mb{G}}|^2d\xi
dyd\tau.
\end{array}
 \label{(4.35)}
\end{equation}

\underline{Step 2. Estimation on $\|\phi_y(\tau,\cdot)\|^2$.}

To estimate the term $\di \int_{\f h\v}^\tau\|\p_y\|^2d\tau$,  we
firstly rewrite the equation $(\ref{sys-h})_2$ as
\begin{equation}
\begin{array}{l}
\di\f{4}{3}\f{\mu(\T)}{V}\p_{y\tau}-\psi_{1\tau}+\f{Z}{V}\phi_y=-\f{4}{3}(\f{\mu({\T})}{V})_y\psi_{1y}
+\f{2}{3V}\z_y+H_1\\
\di\qquad+\int\xi_1^2(\Pi_{1}-\Pi^{CD}_{11}-\Pi^{S_3}_1)_yd\xi
-N_5+\bar Q_{1y}+Q_{1y},
\end{array}
\label{(4.36)}
\end{equation}
by using the equation of conservation of the mass $(\ref{sys-h})_1$.
Multiplying (\ref{(4.36)}) by $\p_y$, we get
\begin{equation}
\begin{array}{l}
\di(\f{2\mu(\T)}{3V}\p_{y}^2-\p_y\psi_1)_{\tau}+(\p_{\tau}\psi_1)_y
+\f{Z}{V}\p_y^2=\psi_{1y}^2\\
\di~~ +\Big[ -\f{4}{3}(\f{\mu({\T})}{V})_y\psi_{1y}
+\f{2}{3V}\z_y+H_1+
\int\xi_1^2(\Pi_{1}-\Pi^{CD}_{11}-\Pi^{S_3}_1)_yd\xi -N_5+\bar
Q_{1y}+Q_{1y}\Big]\p_y.
\end{array}
\label{phi-y}
\end{equation}
Integrating (\ref{phi-y}) with respect to $\tau,y$ and using
\eqref{bar-Q-E}, \eqref{N5-8} and the Cauchy inequality yield
\begin{equation}
\begin{array}{l}
\di \|\p_y(\tau,\cdot)\|^2+\int_{\f h\v}^\tau\|\p_y\|^2d\tau\leq
C\|\psi_1(\tau,\cdot)\|^2+ C\int_{\f
h\v}^\tau\|(\psi_{y},\z_y)\|^2d\tau\\
\di~~ +C(\d+C_{h,T}\chi)\int_{\f
h\v}^\tau\|(\phi,\psi,\zeta)\|^2d\tau+C\chi\int_{\f
h\v}^\tau\|\psi_{1yy}\|^2d\tau\\
\di ~~+C_{h,T}~\v^{\f12}+\int_{\f h\v}^\tau\int
|\int\xi_1^2(\Pi_{1}-\Pi^{CD}_{11}-\Pi^{S_3}_1)_yd\xi|^2 dyd\tau.
\end{array}
\label{(4.38)}
\end{equation}
For the microscopic term $\di \int_{\f h\v}^\tau\int
|\int\xi_1^2(\Pi_{1}-\Pi^{CD}_{11}-\Pi^{S_3}_1)_yd\xi|^2 dyd\tau$,
by \eqref{Pi-1}, we have
\begin{equation}
\begin{array}{l}
 \di K_4:=\int_{\f
h\v}^\tau\int |\int\xi_1^2(\Pi_{1}-\Pi^{CD}_{11}-\Pi^{S_3}_1)_yd\xi|^2 dyd\tau\\
 \di \leq
C\Big[\int_{\f
h\v}^\tau\int|\int\xi_1^2(\mb{L}_\mb{M}^{-1}\wt{\mb{G}}_{\tau})_y
d\xi|^2 dyd\tau+\int_{\f
h\v}^\tau\int|\int\xi_1^2(\mb{L}_\mb{M}^{-1}\f{u_1}v\wt{\mb{G}}
_{1y})_y d\xi|^2dyd\tau\\
\di+\int_{\f
h\v}^\tau\int|\int\xi_1^2[\mb{L}_\mb{M}^{-1}\f{1}{v}\mb{P}_1(\xi_1\wt{\mb{G}}
_{1y})]_yd\xi|^2dyd\tau +\int_{\f
h\v}^\tau\int|\int\xi_1^2[\mb{L}_\mb{M}^{-1}Q(\wt{\mb{G}}_1
,\wt{\mb{G}}_1)]_yd\xi|^2dyd\tau\\
\di+\int_{\f
h\v}^\tau\int|\int\xi_1^2[\mb{L}_\mb{M}^{-1}Q(\wt{\mb{G}}_1
,\mb{G}^{R_1}+\mb{G}^{CD}+\mb{G}^{S_3})]_yd\xi|^2dyd\tau\\
\di +\int_{\f
h\v}^\tau\int|\int\xi_1^2(J_{3y}+J_{4y}+J_{5y})d\xi|^2dyd\tau\Big]:=\sum_{i=1}^6K_4^{i}.
\end{array}
\label{K4}
\end{equation}
Then we have
\begin{equation}
\begin{array}{ll}
\di K_4^1&\di \leq C\int_{\f
h\v}^\tau\int|\int\xi_1^2\mb{L}_\mb{M}^{-1}\wt{\mb{G}}_{y\tau}
d\xi|^2 dyd\tau+C\int_{\f
h\v}^\tau\int|\int\xi_1^2\mb{L}_\mb{M}^{-1}\{Q(\mb{L}_\mb{M}^{-1}\wt{\mb{G}}_\tau,\mb{M}_y)\}d\xi|^2
dyd\tau\\
&\di \leq C\sum_{|\a|=2} \int_{\f
h\v}^\tau\int\int\f{\nu^{-1}(|\xi|)}{\mb{M}_\star}|\partial^\a
\wt{\mb{G}} |^2d\xi dyd\tau+C(\d+C_{h,T}\chi)\int_{\f h\v}^\tau\int
\int\f{\nu(|\xi|)|\wt{\mb{G}}_\tau|^2}{\mb{M}_\star} d\xi dyd\tau.
\end{array}
\label{K41}
\end{equation}
Similar estimates hold for $K_4^i~(i=2,3)$. Moreover,
\begin{equation}
\begin{array}{ll}
\di K_4^4&\di \leq C\int_{\f
h\v}^\tau\int|\int\xi_1^2\mb{L}_\mb{M}^{-1}Q(\wt{\mb{G}}_1,\wt{\mb{G}}_{1y})
d\xi|^2 dyd\tau\\
&\di\quad + C\int_{\f
h\v}^\tau\int|\int\xi_1^2\mb{L}_\mb{M}^{-1}\{Q(\mb{L}_\mb{M}^{-1}Q(\wt{\mb{G}}_1,\wt{\mb{G}}_1),\mb{M}_y)\}d\xi|^2
dyd\tau\\
&\di \leq C_{h,T}\chi\int_{\f
h\v}^\tau\int\int\f{\nu(|\x|)(|\wt{\mb{G}}_1|^2+|\wt{\mb{G}}_y|^2)}{\mb{M}_\star}d\x
dy d\tau+C_{h,T}~\v^{\f12},
\end{array}
\label{K44}
\end{equation}
\begin{equation}
\begin{array}{ll}
\di K_4^5 \leq C(\d+C_{h,T}\chi)\int_{\f
h\v}^\tau\int\int\f{\nu(|\x|)(|\wt{\mb{G}}_1|^2+|\wt{\mb{G}}_y|^2)}{\mb{M}_\star}d\x
dy d\tau+C_{h,T}~\v^{\f12},
\end{array}
\label{K45}
\end{equation}
and
\begin{equation}
\begin{array}{ll}
\di K_4^6&\di \leq C\int_{\f
h\v}^\tau\int\int\f{|J_{3y}|^2+|J_{4y}|^2+|J_{5y}|^2}{\mb{M}_\star}d\xi dyd\tau\\
&\di \leq C(\d+C_{h,T}\chi)\int_{\f
h\v}^\tau\|(\phi,\psi,\z,\phi_y,\psi_y,\z_y,\phi_{yy},\psi_{yy},\z_{yy})\|^2d\tau+C_{h,T}~\v^{\f12}.
\end{array}
\label{K46}
\end{equation}
Substituting \eqref{K41}-\eqref{K46} into \eqref{K4} gives
\begin{equation}
\begin{array}{ll}
K_4&\di \leq C\sum_{|\a|=2} \int_{\f
h\v}^\tau\int\int\f{\nu^{-1}(|\xi|)}{\mb{M}_\star}|\partial^\a
\wt{\mb{G}} |^2d\xi
dyd\tau+C_{h,T}~\v^{\f12}\\
&\di \quad+C(\d+C_{h,T}\chi)\int_{\f
h\v}^\tau\|(\phi,\psi,\z,\phi_y,\psi_y,\z_y,\phi_{yy},\psi_{yy},\z_{yy})\|^2d\tau\\
&\di\quad +C(\d+C_{h,T}\chi)\int_{\f h\v}^\tau\int
\int\f{\nu(|\xi|)(|\wt{\mb{G}}_\tau|^2+|\wt{\mb{G}}_1|^2+|\wt{\mb{G}}_y|^2)}{\mb{M}_\star}
d\xi dyd\tau.
\end{array}
\label{K4-E}
\end{equation}
Thus we have
\begin{equation}
\begin{array}{l}
\di \|\p_y(\tau,\cdot)\|^2+\int_{\f h\v}^\tau\|\p_y\|^2d\tau\leq
C\|\psi_1(\tau,\cdot)\|^2+ C\int_{\f
h\v}^\tau\|(\psi_{y},\z_y)\|^2d\tau\\
\di \quad+C(\d+C_{h,T}\chi)\int_{\f
h\v}^\tau\|(\phi_{yy},\psi_{yy},\z_{yy})\|^2d\tau +C_{h,T}~\v^{\f12}\\
\di \quad +C(\d+C_{h,T}\chi)\int_{\f
h\v}^\tau\|(\phi,\psi,\zeta)\|^2d\tau+ C\sum_{|\a|=2}\int_{\f
h\v}^\tau \int\int\f{\nu^{-1}(|\xi|)}{\mb{M}_\star}|\partial^\a
\mb{G} |^2d\xi
dyd\tau\\
\di\quad +C(\d+C_{h,T}\chi)\int_{\f h\v}^\tau\int
\int\f{\nu(|\xi|)(\sum_{|\a^\prime|=1}|\partial^{\a^\prime}\wt{\mb{G}}|^2+|\wt{\mb{G}}_1|^2)}{\mb{M}_\star}
d\xi dyd\tau.
\end{array}
\label{(4.43)}
\end{equation}

We now turn to the time derivatives. To estimate
$\|(\p_{\tau},\psi_{\tau},\z_{\tau})\|^2$, we need to use the system
(\ref{(4.18)}). By multiplying $(\ref{(4.18)})_1$ by $\p_{\tau}$,
$(\ref{(4.18)})_2$ by $\psi_{1\tau}$, $(\ref{(4.18)})_3$ by
$\psi_{i\tau}~(i=2,3)$ and $(\ref{(4.18)})_4$ by $\z_{\tau}$
respectively, and adding them together, after integrating with
respect to $\tau$ and $y$, we have
\begin{equation}
\begin{array}{l}
\di \int_{\f
h\v}^\tau\|(\p_{\tau},\psi_{\tau},\z_{\tau})(\tau,\cdot)\|^2d\tau\leq
C(\d+C_{h,T}\v^{\f12})\int_{\f h\v}^\tau\|(\phi,\psi,\z)\|^2d\tau+C\int_{\f h\v}^\tau\|(\p_y,\psi_y,\z_y)\|^2d\tau\\
\di\qquad +C\int_{\f
h\v}^\tau\int\int\f{\nu(|\xi|)}{\mb{M}_\star}|\wt{\mb{G}}_y|^2d\xi
dyd\tau+C_{h,T}~\v^{\f12}.
\end{array}
\label{(4.45)}
\end{equation}

\underline{Step 3. Estimation on
$\|(\phi_y,\psi_y,\z_y)(\tau,\cdot)\|^2$.}

\

Multiplying $\eqref{sys-h-h}_1$ by $\phi_{y}$, $\eqref{sys-h-h}_2$
by $\f{V}{Z}\psi_{1y}$, $\eqref{sys-h-h}_3$ by $\psi_{iy}$, and
$\eqref{sys-h-h}_4$ by $\f{2}{3Z^2}\z_{y}$, adding them together
gives
\begin{equation}
\begin{array}{l}
\di
(\f{\phi_y^2}{2}+\f{V}{2Z}\psi_{1y}^2+\sum_{i=2}^3\f{\psi_{iy}^2}{2}+\f{\z_y^2}{3Z^2})_{\tau}
+\f{4\mu(\T)}{3Z}\psi_{1yy}^2+\sum_{i=2}^3\f{\mu(\T)}{V}\psi_{iyy}^2
+\f{\k(\T)}{3Z^2V}\z_{yy}^2\\
\di
=I_7(\phi,\psi,\z,\phi_y,\psi_y,\z_y,\psi_{yy},\z_{yy})-(N_5-\bar
Q_{1y}-Q_{1y})(\f{V\psi_{1y}}{Z})_y-\sum_{i=2}^3(N_{i+4}-\bar
Q_{iy}-Q_{iy})\psi_{iyy}
\\[0.3cm]
\di -(\f{2\z_{y}}{3Z^2})_y\big[N_8-\bar
Q_{4y}-Q_{4y}+\sum_{i=1}^3U_i(\bar
Q_{iy}+Q_{iy})\big]+K_5+(\cdots)_y,
\end{array}
\label{(4.48)}
\end{equation}
where
\begin{equation}
\begin{array}{ll}
\di I_7(\phi,\psi,\z,\phi_y,\psi_y,\z_y,\psi_{yy},\z_{yy})\\
\di\quad=(\f
V{2Z})_\tau\psi_{1y}^2+(\f{1}{3Z^2})_\tau\z_y^2+(\f{2}{3Z})_y\psi_{1y}\z_y\
-\f{V}{Z}\psi_{1y}H_5-\sum_{i=2}^3\psi_{iy}H_{i+4}\\
\di\qquad-\f{2}{3Z^2}\z_yH_8-(\f{4\mu(\T)\psi_{1y}}{3V})_y(\f{V}{Z})_y\psi_{1y}-(\f{4\mu(\T)}{3V})_y\psi_{1y}\f{V}{Z}\psi_{1yy}\\
\di\qquad-\sum_{i=2}^3(\f{\mu(\T)}{V})_y\psi_{iy}\psi_{iyy}-(\f{\k(\T)\z_y}{V})_y(\f{2}{3Z^2})_y\z_y-(\f{\k(\T)}{V})_y\z_y\f{2}{3Z^2}\z_{yy}\\[3mm]
\di \quad\leq
\b|(\phi_{yy},\psi_{yy},\z_{yy})|^2+C_\b(\d+C_{h,T}\v^{\f12})|(\phi,\psi,\z,\phi_y,\psi_y,\z_y)|^2,
\end{array}\label{I7}
\end{equation}
and
\begin{equation}
\begin{array}{ll}
K_5&\di
=(\f{V\psi_{1y}}{Z})_y\int\x_1^2(\Pi_1-\Pi^{CD}_{11}-\Pi^{S_3}_1)_yd\x
+\sum_{i=2}^3\psi_{iyy}\int\x_1\x_i(\Pi_1-\Pi^{CD}_{11}-\Pi^{S_3}_1)_yd\x\\
&\di+(\f{2\z_{y}}{3Z^2})_y\Big[\int\x_1\f{|\x|^2}{2}(\Pi_1-\Pi^{CD}_{11}-\Pi^{S_3}_1)_yd\x
-\sum_{i=1}^3\psi_i\int\xi_1\xi_i\Pi_{1y}d\xi\\
&\di\qquad\qquad~-\sum_{i=1}^3U_i\int\xi_1\x_i(\Pi_1-\Pi^{CD}_{11}-\Pi^{S_3}_1)_yd\xi\Big].
\end{array}\label{K5}
\end{equation}

Integrating (\ref{(4.48)}) with respect to $\tau,y$, and
substituting \eqref{bar-Q-E}, \eqref{N5-8} and \eqref{K4-E} into
(\ref{(4.48)}) and \eqref{K5}, choosing $\b,\chi,\v$ small enough,
we have
\begin{equation}
\begin{array}{l}
\di \|(\phi_y,\psi_{y},\z_y)(\tau,\cdot)\|^2+\int_{\f
h\v}^\tau\|(\psi_{yy},\z_{yy})\|^2d\tau
\\
\di \leq \b\int_{\f
h\v}^\tau\|\phi_{yy}\|^2d\tau+C_\b(\d+C_{h,T}\v^{\f12})\int_{\f
h\v}^\tau\|(\phi,\psi,\z,\phi_y,\psi_y,\z_y)\|^2d\tau
+C_{h,T}~\v^{\f12}\\
\quad\di +C_{h,T}\chi\int_{\f
h\v}^\tau\int\int\f{\nu(|\xi|)}{\mb{M}_\star}|\wt{\mb{G}} _1|^2d\xi
dyd\tau+C\sum_{|\a|=2}\int_{\f
h\v}^\tau\int\int\f{\nu(|\xi|)}{\mb{M}_\star}|\partial^\a
\wt{\mb{G}} |^2d\xi
dyd\tau\\
\quad\di +C_{h,T}\chi\sum_{|\a^{\prime}|=1}\int_{\f
h\v}^\tau\int\int\f{\nu(|\xi|)}{\mb{M}_\star}|\partial^{\a^\prime}\wt{\mb{G}}
|^2d\xi dyd\tau.
\end{array}
\label{(4.49)}
\end{equation}

Again, to recover  $\|\p_{yy}\|^2$ in the dissipation rate, applying
$\partial_y$ to $(\ref{(4.18)})_2$, we get
\begin{equation}
\psi_{1y\tau}+(p-P)_{yy}
=-\f{4}{3}\Big[\f{\mu(\T)}{V}U_{1y}-\f{\mu(\T^{S_3})}{V^{S_3}}U^{S_3}_{1y}\Big]_{yy}-\int\xi_1\x_i\Pi^{CD}_{11yy}d\x-\bar
Q_{1yy}- Q_{1yy}-\int\xi_1^2\wt{\mb{G}}_{yy}d\xi. \label{(4.50)}
\end{equation}
Note that
\begin{equation}
(p-P)_{yy}=-\f{p}{v}\p_{yy}+\f{2}{3v}\z_{yy}-\f1v(p-P)V_{yy}
-\f{\p}{v}P_{yy}-\f{2v_y}{v}(p-P)_y-\f{2P_y}{v}\p_y. \label{(4.51)}
\end{equation}
Multiplying (\ref{(4.50)}) by $-\p_{yy}$ and using (\ref{(4.51)})
imply
\begin{equation}
\begin{array}{l}
\di -\int\psi_{1y}\p_{yy}(\tau,y)dy+\int_{\f
h\v}^\tau\int\f{p}{2v}\p_{yy}^2dyd\tau\\
\di\quad\leq C(\d+C_{h,T}\chi)\int_{\f
h\v}^\tau\|(\p,\psi,\z,\p_y,\psi_y,\z_y)\|^2d\tau
\\
\qquad\di+C\int_{\f
h\v}^\tau\|(\psi_{1yy},\z_{yy})\|^2d\tau+C_{h,T}~\v^{\f12}+C
\sum_{|\a|=2}\int_{\f
h\v}^\tau\int\int\f{\nu(|\xi|)}{\mb{M}_\star}|\partial^\a
\wt{\mb{G}} |^2d\xi dyd\tau,
\end{array}
\label{(4.52)}
\end{equation}
where we have used the fact that
\begin{equation}
\int_{\f h\v}^\tau\int(|\bar Q_{1yy}|^2+|Q_{1yy}|^2) dy d\tau\leq
C_{h,T}~\v^{\f12}.
\end{equation}

To estimate $\|(\p_{y\tau},\psi_{y\tau},\z_{y\tau})\|^2$ and
$\|(\p_{\tau\tau},\psi_{\tau\tau},\z_{\tau\tau})\|^2$, we use the
system (\ref{(4.18)}) again. Applying $\partial_y$ to
(\ref{(4.18)}), and multiplying the four equations of (\ref{(4.18)})
by $\p_{y\tau}$, $\psi_{1y\tau}$, $\psi_{iy\tau}$ $(i= 2,3)$,
$\z_{y\tau}$ respectively, then adding them together and integrating
with respect to $\tau$ and $y$, we have
\begin{equation}
\begin{array}{l}
\di \int_{\f
h\v}^\tau\|(\p_{y\tau},\psi_{y\tau},\z_{y\tau})\|^2d\tau\leq
 C\int_{\f h\v}^\tau\|(\p_{yy},\psi_{yy},\z_{yy})\|^2d\tau+C_{h,T}~\v^{\f12}\\
\quad\di+C_{h,T}\chi\int_{\f
h\v}^\tau\|(\p_y,\psi_y,\z_y)\|^2d\tau+C(\d+C_{h,T}\v^{\f12})\int_{\f
h\v}^\tau\|(\phi,\psi,\z)\|^2d\tau\\
\quad\di+C\sum_{|\a|=2}\int_{\f
h\v}^\tau\int\int\f{\nu(|\xi|)}{\mb{M}_\star}|\partial^\a
 \wt{\mb{G}}|^2d\xi dyd\tau.
\end{array}
\label{(4.53)}
\end{equation}
Similarly, we can obtain
\begin{equation}
\begin{array}{l}
\di \int_{\f
h\v}^\tau\|(\p_{\tau\tau},\psi_{\tau\tau},\z_{\tau\tau})\|^2d\tau\leq
 C\int_{\f h\v}^\tau\|(\p_{y\tau},\psi_{y\tau},\z_{y\tau})\|^2d\tau+C_{h,T}~\v^{\f12}\\
\quad\di+C_{h,T}\chi\sum_{|\a^\prime|=1}\int_{\f
h\v}^\tau\|\partial^{\a^\prime}(\p,\psi,\z)\|^2d\tau+C(\d+C_{h,T}\v^{\f12})\int_{\f
h\v}^\tau\|(\phi,\psi,\z)\|^2d\tau\\
\quad\di+C\sum_{|\a|=2}\int_{\f
h\v}^\tau\int\int\f{\nu(|\xi|)}{\mb{M}_\star}|\partial^\a
 \wt{\mb{G}}|^2d\xi dyd\tau.
\end{array}
\label{(4.54)}
\end{equation}
A suitable linear combination of \eqref{(4.49)}-\eqref{(4.54)} gives
\begin{equation}
\begin{array}{l}
\di\|(\phi_y,\psi_{y},\z_y)(\tau,\cdot)\|^2-\int\psi_{1y}\phi_{yy}
dy
+\sum_{|\a|=2}\int_{\f h\v}^\tau\|\partial^\a(\p,\psi,\z)\|^2d\tau\\
\di \leq C\sum_{|\a|=2}\int_{\f
h\v}^\tau\int\int\f{\nu(|\xi|)}{\mb{M}_\star}|\partial^\a
\wt{\mb{G}}|^2d\xi dyd\tau+C\sum_{|\a^\prime|=1}\int_{\f
h\v}^\tau\int\int\f{\nu(|\xi|)}{\mb{M}_\star}|\partial^{\a^\prime}
\wt{\mb{G}}|^2d\xi
dyd\tau\\
\di\quad+C_{h,T}\chi\int_{\f
h\v}^\tau\int\int\f{\nu(|\xi|)}{\mb{M}_\star}|\wt{\mb{G}}_1|^2d\xi
dyd\tau+C(\d+C_{h,T}\v^{\f12})\int_{\f
h\v}^\tau\|(\phi,\psi,\z)\|^2d\tau
\\
\di\quad +C_{h,T}\chi\sum_{|\a^\prime|=1}\int_{\f
h\v}^\tau\|\partial^{\a^\prime}(\p,\psi,\z)\|^2d\tau +C_{h,T}
~\v^{\f12}.
\end{array}
\label{(4.55)}
\end{equation}

\underline{Step 4: Estimation on the non-fluid component:}

\

To close the a priori estimate, we also need to estimate the
derivatives on the non-fluid component $\wt{\mb{G}} $, i.e.,
$\partial^\a \wt{\mb{G}}, (|\a|=1,2)$. For this, from \eqref{t-G},
we obtain
\begin{equation}\label{M.11}
v\wt{\mb{G}}_{\tau}-v\mb{L}_\mb{M}\wt{\mb{G}}=u_1\wt{\mb{G}}
_y-\mb{P}_1(\xi_1\wt{\mb{G}}_y)-v\Big[\f{1}{v}\mb{P}_1(\xi_1\mb{M}_y)
-\f{1}{V^{S_3}}\mb{P}^{S_3}_1(\xi_1\mb{M}^{S_3}_y)\Big]+2vQ(\wt{\mb{G}},\mb{G}^{S_3})+v
Q(\wt{\mb{G}},\wt{\mb{G}})+vJ_1.
\end{equation}
Applying $\partial_y$ on \eqref{M.11}, we have
\begin{eqnarray}\label{M.12}
&&v\wt{\mb{G}}_{y\tau}-(v\mb{L}_\mb{M}\wt{\mb{G}})_y=\bigg\{u_1\wt{\mb{G}}
_y-\mb{P}_1(\xi_1\wt{\mb{G}}_y)
-v\Big[\f{1}{v}\mb{P}_1(\xi_1\mb{M}_y)
-\f{1}{V^{S_3}}\mb{P}^{S_3}_1(\xi_1\mb{M}^{S_3}_y)\Big]\nonumber\\
&&\qquad\qquad\qquad\qquad\qquad+2vQ(\wt{\mb{G}},\mb{G}^{S_3})+v
Q(\wt{\mb{G}},\wt{\mb{G}})+vJ_1\bigg\}_y-v_y\wt{\mb{G}}_{\tau}.
\end{eqnarray}
Multiplying \eqref{M.12} by $\f{\wt{\mb{G}}_y}{\mb{M}_\star}$ gives
\begin{equation}
\begin{array}{ll}
\di
\bigg(\f{v\wt{\mb{G}}_y^2}{2\mb{M}_\star}\bigg)_{\tau}-\f{\wt{\mb{G}}_y}{\mb{M}_\star}(v\mb{L}_\mb{M}\wt{\mb{G}})_y=\bigg\{u_1\wt{\mb{G}}
_y-\mb{P}_1(\xi_1\wt{\mb{G}}_y)
-v\Big[\f{1}{v}\mb{P}_1(\xi_1\mb{M}_y)
-\f{1}{V^{S_3}}\mb{P}^{S_3}_1(\xi_1\mb{M}^{S_3}_y)\Big]\\
\qquad\qquad\qquad\di +2vQ(\wt{\mb{G}},\mb{G}^{S_3})+v
Q(\wt{\mb{G}},\wt{\mb{G}})+vJ_1\bigg\}_y\f{\wt{\mb{G}}_y}{\mb{M}_\star}-v_y\wt{\mb{G}}_{\tau}\f{\wt{\mb{G}}_y}{\mb{M}_\star}
+v_{\tau}\f{\wt{\mb{G}}_y^2}{2\mb{M}_\star}.
\end{array}\label{M.13}
\end{equation}
Then Lemma \ref{Lemma 4.1}-Lemma \ref{Lemma 4.4} and Cauchy inequality
give
\begin{eqnarray}\label{M.14}
&&-\int_{\f
h\v}^\tau\int\int\f{\wt{\mb{G}}_y}{\mb{M}_\star}(v\mb{L}_\mb{M}\wt{\mb{G}})_yd\xi
dyd\tau\nonumber\\
&&=-\int_{\f
h\v}^\tau\int\int\f{\wt{\mb{G}}_y}{\mb{M}_\star}\bigg(v\mb{L}_\mb{M}\wt{\mb{G}}_y
+2v_yQ(\mb{M},\wt{\mb{G}})+2vQ(M_y,\wt{\mb{G}})\bigg)d\xi
dyd\tau\nonumber\\
&&\geq \f{7\wt\s}{8}\int_{\f
h\v}^\tau\int\int\f{\nu(|\xi|)\wt{\mb{G}}_y^2}{\mb{M}_\star}d\xi
dyd\tau-C(\d+C_{h,T}\chi)\int_{\f
h\v}^\tau\int\int\f{\nu(|\xi|)\wt{\mb{G}}_1^2}{\mb{M}_\star}d\xi
dyd\tau\nonumber\\
&&\quad -C\v\int_{\f
h\v}^\tau\|(\p_y,\psi_y,\zeta_y)\|^2d\tau-C_{h,T}~\v^{\f12},
\end{eqnarray}
and
\begin{eqnarray}\label{M.15}
&&\int_{\f
h\v}^\tau\int\int-v_y\wt{\mb{G}}_{\tau}\f{\wt{\mb{G}}_y}{\mb{M}_\star}
+v_{\tau}\f{\wt{\mb{G}}_y^2}{2\mb{M}_\star}+(u_1\wt{\mb{G}}_y-\mb{P}_1(\xi_1\wt{\mb{G}}_y))_y\f{\wt{\mb{G}}_y}{2\mb{M}_\star}d\xi
dyd\tau\nonumber\\
&&\leq \f{\wt\s}{32}\int_{\f
h\v}^\tau\int\int\f{\nu(|\xi|)\wt{\mb{G}}_y^2}{\mb{M}_\star}d\xi
dyd\tau+C(\d+C_{h,T}\chi)\int_{\f
h\v}^\tau\int\int\f{\nu(|\xi|)}{\mb{M}_\star}\wt{\mb{G}}_{\tau}^2d\xi
dyd\tau\nonumber\\
&&\qquad\qquad+C\int_{\f
h\v}^\tau\int\int\f{\nu(|\xi|)\wt{\mb{G}}_{yy}^2}{\mb{M}_\star}d\xi
dyd\tau.
\end{eqnarray}
Note that
\begin{eqnarray}\label{M.19}
&&J_6\triangleq -v\Big[\f{1}{v}\mb{P}_1(\xi_1\mb{M}_y)
-\f{1}{V^{S_3}}\mb{P}^{S_3}_1(\xi_1\mb{M}^{S_3}_y)\Big]\nonumber\\
&&=-\f{3}{2\t}\mb{P}_1\bigg[\xi_1\mb{M}\Big(\f{|\xi-u|^2}{2\t}(\t_y-\T^{S_3}_y)+\xi\cdot
(u_y-U^{S_3}_y)\Big)\bigg]-v\bigg[\f{3}{2v\t}\mb{P}_1[\xi_1\mb{M}(\f{|\xi-u|^2}{2\t}\T^{S_3}_y+\xi\cdot
U^{S_3}_y)]\nonumber\\
&&\qquad-\f{3}{2V^{S_3}\T^{S_3}}\mb{P}^{S_3}_1[\xi_1\mb{M^{S_3}}(\f{|\xi-U^{S_3}|^2}{2\T^{S_3}}\T^{S_3}_y+\xi\cdot
U^{S_3}_y)]\bigg].
\end{eqnarray}
Then Lemma \ref{Lemma 4.1}-Lemma \ref{Lemma 4.4} and wave
interaction estimates imply that
\begin{eqnarray}\label{M.23}
&&\int_{\f h\v}^\tau\int\int
(vJ_{1}+J_6)_y\f{\wt{\mb{G}}_y}{2\mb{M}_\star}d\xi dyd\tau\nonumber\\
&& \leq  \f{\wt\s}{32}\int_{\f h\v}^\tau\int\int
\f{\nu(|\xi|)\wt{\mb{G}}_y^2}{\mb{M}_\star}d\xi dyd\tau
+C(\d+C_{h,T}\chi)\int_{\f h\v}^\tau\int\int
\f{\nu(|\xi|)\wt{\mb{G}}_1^2}{\mb{M}_\star}d\xi dyd\tau\nonumber\\
&&\qquad+C(\d+C_{h,T}\chi)\int_{\f
h\v}^\tau\|(\p,\psi,\zeta)\|_{H^1(dy)}^2d\tau+C_{h,T}\v^{\f12},
\end{eqnarray}
and
\begin{eqnarray}\label{M.24}
&&\int_{\f h\v}^\tau\int\int
\big(vQ(\wt{\mb{G}},\wt{\mb{G}})\big)_y\f{\wt{\mb{G}}_y}{2\mb{M}_\star}d\xi dyd\tau\nonumber\\
&&=\int_{\f h\v}^\tau\int\int
\big(v_yQ(\wt{\mb{G}},\wt{\mb{G}})+2vQ(\wt{\mb{G}},\wt{\mb{G}}_y)\big)\f{\wt{\mb{G}}_y}{2\mb{M}_\star}d\xi dyd\tau\nonumber\\
&&\leq C_{h,T}\v^{\f12} +\f{\wt\s}{32}\int_{\f h\v}^\tau\int\int
\f{\nu(|\xi|)\wt{\mb{G}}_y^2}{\mb{M}_\star}d\xi
dyd\tau+C(\chi+\d)\int_{\f h\v}^\tau\int\int
\f{\nu(|\xi|)\wt{\mb{G}}_1^2}{\mb{M}_\star}d\xi dyd\tau.
\end{eqnarray}

Thus, integrating \eqref{M.13} with respect to  $\x,y$ and $\tau$
and using \eqref{M.14}, \eqref{M.15} and \eqref{M.23}, \eqref{M.24},
we obtain
\begin{eqnarray}\label{M.25}
&&\int\int\f{|\wt{\mb{G}} _y|^2}{2\mb{M}_\star}(\tau,y,\x)d\xi
dy+\int_{\f
h\v}^\tau\int\int\f{\nu(|\xi|)}{\mb{M}_\star}|\wt{\mb{G}} _y|^2d\xi
dyd\tau\nonumber\\
&& \leq C_{h,T}\v^{\f12} +C\int_{\f
h\v}^\tau\int\int\f{\nu(|\xi|)}{\mb{M}_\star}|\wt{\mb{G}}
_{yy}|^2d\xi dyd\tau+C\int_{\f
h\v}^\tau\|(\p_{yy},\z_{yy})\|^2d\tau\nonumber\\
&&\quad +C(\d+C_{h,T}\chi)\int_{\f
h\v}^\tau\int\int\f{\nu(|\xi|)}{\mb{M}_\star}(|\wt{\mb{G}}_1|^2+|\wt{\mb{G}}_\tau|^2)d\xi
dyd\tau\nonumber\\
&&\qquad+C(\d+C_{h,T}\chi)\int_{\f
h\v}^\tau\big[\|(\p,\psi,\z)\|^2+\|(\p,\psi,\z)_y\|^2\big]d\tau.
\end{eqnarray}

Similarly,
\begin{eqnarray}\label{M.26}
&&\int\int\f{|\wt{\mb{G}} _\tau|^2}{2\mb{M}_\star}(\tau,y,\x)d\xi
dy+\int_{\f
h\v}^\tau\int\int\f{\nu(|\xi|)}{\mb{M}_\star}|\wt{\mb{G}}
_\tau|^2d\xi
dyd\tau\nonumber\\
&& \leq C_{h,T}\v^{\f12} +C\int_{\f
h\v}^\tau\int\int\f{\nu(|\xi|)}{\mb{M}_\star}|\wt{\mb{G}}
_{y\tau}|^2d\xi dyd\tau+C\int_{\f
h\v}^\tau\|(\p_{y\tau},\z_{y\tau})\|^2d\tau\nonumber\\
&&\quad +C(\d+C_{h,T}\chi)\int_{\f
h\v}^\tau\int\int\f{\nu(|\xi|)}{\mb{M}_\star}(|\wt{\mb{G}}_1|^2+|\wt{\mb{G}}_y|^2)d\xi
dyd\tau\nonumber\\
&&\qquad+CC(\d+C_{h,T}\chi)\int_{\f
h\v}^\tau\big[\|(\p,\psi,\z)\|^2+\|(\p,\psi,\z)_\tau\|^2\big]d\tau.
\end{eqnarray}

\underline{Step 5: Highest order estimates:}

 Finally, we estimate the highest
order derivatives, that is, $\di -\int\psi_{1y}\p_{yy}dy$ and\\ $\di
\int_{\f h\v}^\tau\int\int \f{\nu(|\xi|)|\partial^\a \wt{\mb{G}}
|^2}{\mb{M}_\star}d\xi dyd\tau$ with $|\a|=2$ in (\ref{(4.55)}). To
do so, it is sufficient to study  $\di \int\int \f{|\partial^\a
\wt{f}|^2}{\mb{M}_\star}d\xi dy~(|\a|=2)$ in view of \eqref{(4.25)}
and \eqref{alpha=2}. Using the same idea in \cite{Huang-Wang-Yang},
we obtain the estimation for the highest order derivative terms,
i.e,
\begin{eqnarray}\label{M.45}
&&\sum_{|\a|=2}\int\int \f{|\partial^\a
\wt{f}|^2}{2\mb{M}_\star}(\tau,y,\x)d\xi dy+\sum_{|\a|=2}\int_{\f
h\v}^\tau\int\int\f{\nu(|\xi|)}{\mb{M}_\star}|\partial^\a
\wt{\mb{G}}|^2d\xi dyd\tau\nonumber\\
&&\leq C(\eta_0+\d+C_{h,T}\chi)\int_{\f h\v}^\tau\int\int
\f{\nu(|\xi|)}{\mb{M}_\star}(|\wt{\mb{G}}_y|^2+|\wt{\mb{G}}_\tau|^2+|\wt{\mb{G}}_1|^2)d\x
dyd\tau\nonumber\\
&&\qquad\quad+ C(\eta_0+\d+C_{h,T}\chi)\sum_{|\a|=0}^{2}\int_{\f
h\v}^\tau\|\partial^{\a}(\p,\psi,\z)\|^2d\tau +C_{h,T}~\v^{\f12},
\end{eqnarray}
where $\eta_0$ is defined in Lemma \ref{Lemma 4.2}.

Noting that
$$
\begin{array}{ll}
\di -\int\psi_{1y}\phi_{yy} dy&\di \leq \b \|\psi_{1y}\|^2
+C_\b\|\phi_{yy}\|^2\\
&\di \leq \b \|\psi_{1y}\|^2 +C_\b\sum_{|\a|=2}\int\int
\f{|\partial^\a \wt{f}|^2}{2\mb{M}_\star}(\tau,y,\x)d\xi
dy+C_{h,T,\b}~\v^{\f12},
\end{array}
$$
and combining the estimates \eqref{(4.55)}, \eqref{M.25},
\eqref{M.26} and \eqref{M.45},
 we complete the proof of Proposition \ref{Prop3.2}.

\end{document}